\documentstyle[12pt,twoside]{article}
\newcounter{THNO}[subsection]
\renewcommand{\theTHNO}{\arabic{section}.\arabic{subsection}.\arabic{THNO} }
\def\TH#1{\refstepcounter{THNO}\par\vspace{0.9cm}
\par\noindent\begingroup \sl
\leftskip=3em\hspace{-3em}{\sc\theTHNO\ #1\ }}
\def\ETH{\par\endgroup}
\def\PR{\par{\sc Proof. }}
\author{S.~Kuleshov.}
\title{Exceptional and Rigid Sheaves on
        Surfaces with Anticanonical Class without Base Components.}
\begin{document}
%% FOLLOWING LINE CANNOT BE BROKEN BEFORE 80 CHAR
%%%%%%%%%%%%%%%%%%%%%%%%%%%%%%%%%%%%%%%%%%%%%%%%%%%%%%%%%%%%%%%%%%%%%%%%%%%%%%%%%%
%%%%%%%%%%%%%%%%CUT%%%%%%%%%ALGFONTS.TEX%%%%%%%%%%%%%%%%%%%%%%%%%%%%%%%%%
\font\twelvegtc=eufm10 scaled 1200
\font\tengtc=eufm10
\font\ninegtc=eufm9
\font\sevengtc=eufm7
\font\fivegtc=eufm5
\newfam\gtcfam
\def\gtc{\fam\gtcfam\twelvegtc}
\textfont\gtcfam=\twelvegtc
\scriptfont\gtcfam=\ninegtc
\scriptscriptfont\gtcfam=\sevengtc
\font\twelveBBB=msbm10 scaled 1200
\font\tenBBB=msbm10
\font\sevenBBB=msbm7
\newfam\BBBfam
\def\BBB{\fam\BBBfam\twelveBBB}
\textfont\BBBfam=\twelveBBB
\scriptfont\BBBfam=\tenBBB
\scriptscriptfont\BBBfam=\sevenBBB
\def\QQ{{\BBB Q}}   \def\ZZ{{\BBB Z}}    \def\NN{{\BBB N}}
\def\RR{{\BBB R}}   \def\CC{{\BBB C}}    \def\FF{{\BBB F}}
\def\AA{{\BBB A}}   \def\BB{{\BBB B}}    \def\HH{{\BBB H}}
\def\PP{{\BBB P}}
\mathchardef\kkk="707C \def\kk{{\BBB \kkk}}
\def\gtcp{{\gtc p}}  \def\gtcq{{\gtc q}}  \def\gtcm{{\gtc m}}
\def\gtcA{{\gtc A}}  \def\gtcB{{\gtc B}}  \def\gtcC{{\gtc C}}
\def\gtcG{{\gtc G}}  \def\gtcg{{\gtc g}}  \def\gtch{{\gtc h}}
\def\gtcd{{\gtc d}}
%
%%%%%%%%%%%%%%%%%%%%%%%%%%%%%%%%%%%%%%%%%%%%%%%%%%%% СПЕЦСИМВОЛЫ
%
\def\Ddots{\mathinner{\mkern1mu\raise1pt\hbox{.}\mkern2mu\raise4pt\hbox
                  {.}\mkern2mu\raise7pt\vbox{\kern7pt\hbox{.}}\mkern1mu}}
\def\bydef{\stackrel{\rm def}{=}}
\def\EV#1,#2{\langle #1,#2\rangle}
\def\empset{\not\hspace{-.5ex}\mbox{\Large o}}
\def\deg{{\rm deg}}   \def\im{{\rm im}}     \def\rk{{\rm rk}}
\def\ann{{\rm Ann}}   \def\ev{{\rm ev}}       \def\hom{{\rm Hom}}
\def\aff{{\rm aff}}   \def\Aff{{\rm Aff}}     \def\conv{{\rm conv}}
\def\Vect{{\rm Vect}} \def\codim{{\rm codim}} \def\tr{{\rm tr}}
\newlength{\gvl}
\def\gv#1{\settowidth{\gvl}{$#1$}\stackrel
         {\hspace{4pt}\vcenter{\hrule width\gvl}\hspace{-6pt}\rightarrow}
         {{#1}}}
%
%%%%%%%%%%%%%%%%%%%%%%%%%%%%%%%%%%%%%%%%%%%%%%%%%%%%%% ВЕКТОРОПОДОБНЫЕ
%
\def\VEC#1,#2{({#1}_1,{#1}_2,\dots,{#1}_{#2})}
\def\CVC#1,#2{({#1}^1,{#1}^2,\dots,{#1}^{#2})}
\def\OVEC#1,#2{({#1}_0,{#1}_1,{#1}_2,\dots,{#1}_{#2})}
\def\VVEC#1,#2{\left(\begin{array}{c}{#1}^1\\
                                     {#1}^2\\
                                     \vdots\\
                                     {#1}^{#2}\end{array}\right)}
\def\VCVC#1,#2{\left(\begin{array}{c}{#1}_1\\
                                     {#1}_2\\
                                     \vdots\\
                                     {#1}_{#2}\end{array}\right)}
\def\LC#1,#2,#3{{#1}_1#2_1+{#1}_2#2_2+\cdots+{#1}_{#3}#2_{#3}}
\def\SET#1,#2{\{{#1}_1,{#1}_2,\dots,{#1}_{#2}\}}
\def\FAM#1,#2{\ {#1}_1,{#1}_2,\dots,{#1}_{#2}\ }
%% FOLLOWING LINE CANNOT BE BROKEN BEFORE 80 CHAR
\def\POL#1,#2,#3{{{#1}}_0+{{#1}}_1{#2}+{{#1}}_2{#2}^2+\cdots+{{#1}}_{#3}{#2}^{#3}}
\def\SER#1,#2{{{#1}}_0+{{#1}}_1{#2}+{{#1}}_2{#2}^2+{{#1}}_3{#2}^3+\cdots}
%
%%%%%%%%%%%%%%%%%%%%%%%%%%%%%%%%%%%%%%%%%%%%%%%%%%%%%% ДИАГРАММЫ ЮНГА
%
\def\openrow#1#2#3{\setbox0=\vbox{\hbox
    {\vrule height#2 width#3\kern#2\vrule height#2 width0pt}\hrule height#3}
    \hbox{\leaders\copy0\hskip#1\wd0\vrule width#3}}
\def\row#1#2#3{\vbox{\hrule height#3\openrow{#1}{#2}{#3}}}
\def\Yr#1{\row{#1}{1.5ex}{.1ex}}
\def\yr#1{\openrow{#1}{1.5ex}{.1ex}}
\def\DY#1\endDY{\baselineskip=1ex\lineskip=0pt\lineskiplimit=0pt{\vcenter
    {\Yr#1}}}
\def\openclm#1#2#3{\setbox0=\vbox{\hrule height#3\hbox
    {\vrule width0pt\kern#2\vrule width#3 height#2}}\vtop
    {\leaders\copy0\vskip#1\ht0\hrule height#3}}
\def\clm#1#2#3{\hbox{\vrule width#3\openclm{#1}{#2}{#3}}}
\def\Yc#1{\clm{#1}{1.5ex}{.1ex}}
\def\yc#1{\openclm{#1}{1.5ex}{.1ex}}
\def\CDY#1\endCDY{{\vcenter{\hbox{\Yc#1}}}}
%

%% FOLLOWING LINE CANNOT BE BROKEN BEFORE 80 CHAR
%%%%%%%%%%%%%%%%%%%%%%%%%%%%%%%%%%%%%%%%%%%%%%%%%%%%%%%%%%%%%%%%%%%%%%%%%%%%%%%%%%
\maketitle
\begin{abstract}
     The paper consists of three parts. In the first of them different
     kinds stability are discussed. In particular, the stability concept
     with respect to nef divisor is introduced.
     A structure of rigid and superrigid  vector bundles on
     smooth projective  surfaces with nef anticanonical class is studied
     in the second part. We prove that any superrigid bundle has
        a unique  exceptional filtration.
     In the last part we give constructible description of exceptional bundles
     on these surfaces.
\end{abstract}
%% FOLLOWING LINE CANNOT BE BROKEN BEFORE 80 CHAR
%%%%%%%%%%%%%%%%%%%%%%%%%%%%%%%%%%%%%%%%%%%%%%%%%%%%%%%%%%%%%%%%%%%%%%%%%%%%%%%%%%%%%
%%%%%%%%%%%%%%%%%%%%%%%%%%%%%CUT%%%%%%%%%%%%%%%%%%%%%%e0.tex%%%%
\section*{Introduction.}

\addcontentsline{toc}{section}{Introduction.}

   This paper contains a generalisation of the theory of rigid
       (${\rm Ext}^1(E,E)=0$)
   and exceptional (${\rm Ext}^0(E,E)={\BBB C},\; {\rm Ext}^i(E,E)=0$ for
   $i>0$) sheaves on Del Pezzo surfaces , which was described in
   \cite{OK}.
   The objects of this paper  are rigid and exceptional sheaves
   on smooth projective surfaces  $S$ such that $-K_S$ is nef.

   If $K_S^2>0$ then these surfaces may be obtained from ${\BBB P}^2$
   by successively blowing up at most 8 points and are natural
   extension of Del Pezzo surface class.

   At first the exceptional sheaves appeared in   \cite{DP} for
   the description the possible Chern classes which a stable bundle on
   ${\BBB P}^2$. Besides in \cite{D} was proved that any rigid bundle
   on projective plane is a direct sum of exceptional ones.

   I proved the same fact for all Del Pezzo surfaces (\cite{OK}).
   But if $-K_S$ is nef then there exists indecomposable and not simple
   rigid bundles, though a structure of them is described in terms of
   exceptional collections. The proof of this statement is a goal of
   the second  part. The information about superrigid bundles
   gives convenient method for researching exceptional sheaves.

   The theory of exceptional bundles on Del Pezzo surfaces
   uses  stability with respect to  anticanonical class. In this paper
   surfaces have the nef anticanonical class. A question is arose:
   is there a sufficient stability notion with respect to nef divisor,
   and which slope axioms are sufficient for constructing the
   having meaning stability theory?
   For example, when does Garder--Narasimhan
   filtration exist? An answer is the first part's subject .

   Finally in the last part of this work we prove a
   constructibility of exceptional bundles on smooth projective surfaces
   $S$ over   $\CC$ with nef anticanonical class and $K_S^2>0$. Here the
   constructibility means that any exceptional bundle can be obtained
   from finite fixed collection of exceptional sheaves by finite procedure.

   The author is grateful to professor A.~N.~Rudakov and to contestants of
   his seminaire for useful discussions. Besides, the author is grateful
   to Sorros fond for a finance support.

\section* {Notations.}

\addcontentsline{toc}{section}{Notations.}

   Let
$X$ be complete algebraic manifold over  $\CC$;

$r(F),c_1(F),c_2(F),...$ the rank and Chern classes of a coherent sheaf
$F$ on $X$;

${\cal O}_X$ or  ${\cal O}$ the trivial line bundle on $X$;

${\cal O}_Y$ (for a closed submanifold  $Y$  in $X$)  the structure sheaf

of submanifold  $Y$,
 which we  consider sometimes as the sheaf on $X$;

${\cal O}_X(D)$ or ${\cal O}(D)$  the line bundle, which corresponds to
 a divisor $D$;

$K_X$ the  canonical class of $X$;

$F(D)$ the tensor product of $F$   and ${\cal O}(D)$;

$F^*$  the dual sheaf, that is the sheaf of local homomorphisms
${\cal H}om_{{\cal O}_X}(F,{\cal O}_X)$;

${\rm Hom}(E,F)$  the space of global maps from $E$ to $F$;

$h^i(E,F)$  dimension of space ${\rm Ext}^i(E,F)$;

$\chi (E,F)$  the Euler characteristic of any two  sheaves, which equals
$\sum (-1)^ih^i(E,F)$;

$\chi (E)$   the Euler characteristic of sheaf, which equals
 $\ h^i ({\cal O}_X,E)$;

    Direct sum  $\bigoplus\limits_{i=1}^{k} F_i$, where $\forall i\
    F_i=F$, we denote as $kF$ or $V\otimes F$ ( $V$ is a vector space
    over $\CC$ and $\dim V=k$).

    We shall identify a bundle with the sheaf of its local sections.
    Sometimes we shall place a long cohomology sequence associated with
    an exact triple into a table. For example, application of functor
    ${\rm Ext}^{\cdot}(F,\cdot)$ to the exact triple

  $$0\longrightarrow A\longrightarrow B\longrightarrow C\longrightarrow 0$$
    gives
    $$\begin{array}{|c|ccccc|}
\hline
  k&{\rm Ext}^k(F,A)&\rightarrow&{\rm Ext}^k(F,B)&\rightarrow&
  {\rm Ext}^k(F,C)\\
\hline
  0&*&&?&&*\\
  1&0&&?&&0\\
  2&*&&?&&*\\
\hline
  \end{array}$$
    This table calculates ${\rm Ext}^1(F,B)$ . From the table it
     follows  that
    ${\rm Ext}^1(F,B)=0$.

%% FOLLOWING LINE CANNOT BE BROKEN BEFORE 80 CHAR
%%%%%%%%%%%%%%%%%%%%%%%%%%%%%%%%%%%%%%%%%%%%%%%%%%%%%%%%%%%%%%%%%%%%%%%%%%%%%%%%%%%%
\section{ Axioms of Stability.}
\subsection{Definitions and Simple Properties.}

    The Gieseker and Mumford--Takemoto stabilities are
     well known.
    Recently the vector stability with respect to
    a collection of polarizations notion  arisen  (\cite{RS}). All these
    theories have
    a slope and similar properties of stable and semistable sheaves.
    In this section we introduce some slope axioms and obtain basic
    properties of stable sheaves.

    {\sc Definition.} Let $\gamma $ be a  function from the set of
    torsion free sheaves on a complete algebraic manifold $X$ over
     ${\BBB C}$ to ${\BBB R}^n$ with the lexicographic order.
     Assume $\gamma $ satisfy the axioms:

    {\sc  SLOPE.1.} For any exact triple of torsion free sheaves
  $$ 0\longrightarrow F\longrightarrow E\longrightarrow G\longrightarrow 0$$
   the following equivalences hold true:
$$\gamma (F)<\gamma (E)\qquad\Longleftrightarrow \qquad
       \gamma (E)<\gamma (G),$$
$$\gamma (F)>\gamma (E)\qquad\Longleftrightarrow \qquad
       \gamma (E)>\gamma (G),$$
$$\gamma (F)=\gamma (E)\qquad\Longleftrightarrow \qquad
       \gamma (E)=\gamma (G);$$
    {\sc  SLOPE.2.} For any two sheaves without torsion
     $F\subset E$  from $r(F)=r(E)$ it  follows that
     $$\gamma (F)\le \gamma (E).$$

    Then we say that $\gamma$ is  a {\it  slope function\/}
    and $\gamma(E)$ is $\gamma$-{\it slope} of $E$ or simply {\it slope}
     of $E$    if no confusion is likely.
     If $\gamma(E)\in \RR ^n$ then  $\gamma$ is called a
    {\it vector slope function\/}.

    {\sc Definition.} A torsion free sheaf $E$ on algebraic manifold
    $X$ is said to be $\gamma $-{\it (semi)stable\/} or simply
    {\it (semi)stable} if for any its subsheaf $F$ with $r(F)<r(E)$
        the following inequality holds true
         $$\gamma (F)<\gamma (E)$$
\centerline{($\gamma (F)\le\gamma (E)$ for semistable).}
    A subsheaf, which contradicts  (semi)stability  is called
    {\it destabilizing\/}.

    \TH{Remark.}\label{l111} 1. If rank of torsion free  sheaves equals
    1 then they are stable
    wiht respect to any slope, because they have not rank zero subsheaves.

    2. By virtue of the lexicographic order on $\RR ^n$
    the function $\gamma =(\gamma _1,\gamma _2,\ldots,\gamma _n)$
    is the slope iff any   $\gamma _i$ satisfies the slope axioms.

    3. For slopes
     $\gamma =(\gamma _1,\gamma _2,\ldots,\gamma _n)$  and
     $\gamma '=(\gamma _1,\gamma _2,\ldots,\gamma _n,\gamma _{n+1},\ldots,
     \gamma _m)$ the following statements are true

     \hspace{0.5cm} a) a $\gamma $-stable sheaf is
                     $\gamma '$-stable;

     \hspace{0.5cm} b) a $\gamma '$-semistable sheaf is
                     $\gamma $-semistable.

\ETH

    \TH{Lemma}\label{l112} A torsion free sheaf  $E$ on a manifold $X$ is
      (semi)stable if for any its subsheaf $F$ such that
      $E/F$ has not torsion
      $$ \gamma (F)<\gamma (E)\qquad (\gamma (F)\le\gamma (E)).$$
\ETH

    \PR  Let $F$  be a subsheaf of $E$ with  $r(F)<r(E)$.
    Denote by
    $T$ the torsion of the sheaf $G=E/F$.  Consider the following
    commutative diagram:
    $$\begin{array}{ccccccccc}
&&0&&&&&&\\
&&\uparrow &&&&&&\\
&&T&&&&0&&\\
&&\uparrow &&&&\uparrow &&\\
0&\longrightarrow&F'&\longrightarrow&E&\longrightarrow&G'&\longrightarrow&0\\
&&\uparrow &&\llap{$\scriptstyle id$}\uparrow &&\uparrow &&\\
 0&\longrightarrow&F&\longrightarrow&E&\longrightarrow&G&\longrightarrow&0\\
    &&\uparrow &&&&\uparrow &&\\
    &&0&&&&T&&\\
    &&&&&&\uparrow &&\\
    &&&&&&0&&\\
      \end{array}.$$

   Here $G'$ is a nonzero  torsion free sheaf and   $r(F)=r(F')$.

    Since $E$ has not torsion, we see that $F'$ is also a  torsion free
    sheaf.  Therefore it follows from {\sc SLOPE.2.} that $\gamma(F)
    \le\gamma (F')$. By the lemma condition, $\gamma (F')<\gamma (E)$. Thus
     $\gamma (F)<\gamma (E)$. This inequality proves the fact that $E$ is
     stable.

\TH{Lemma.}\label{l113} A torsion free sheaf $E$ is stable (semistable)
    if and only if the slope of any its torsion free quotient    $G$
    satisfy the inequality:
    $$\qquad \gamma (E)<\gamma (G)\qquad\qquad (\gamma (E)\le\gamma (G)).$$
\ETH
   The proof follows from \ref{l112}  and    {\sc SLOPE.1}.

\vspace{2ex}

    For studying stability properties, we need the following:
\TH{Remark.}\label{l114} Taking into account the definition
  of stability and  {\sc SLOPE.2} , we obtain that
  for any pair of torsion free sheaves $F\subset E$ ,
  from  semistability of   $E$ follows the inequality  $\gamma (F)\le
  \gamma (E)$ (without a condition on the ranks). If $E$ is
  stable then
  $$\gamma(F)=\gamma(E)\quad\Longrightarrow\quad r(E)=r(F).$$

  Similarly, if $G$ is a quotient of a semistable sheaf $E\quad$ then
    $\gamma(E)\le\gamma (G)$. If $E$ is stable then
  $$\gamma(G)=\gamma(E)\quad\Longrightarrow\quad E=G.$$
\ETH

\TH{Lemma.}\label{l115} Let $E$ and $F$ be semistable sheaves. Suppose
   $\gamma (E)>\gamma (F)$; then $${\rm Hom}(E,F)=0.$$

\ETH

\PR  Assume that there exists a nonzero morphism  $\varphi :E\longrightarrow
F$.
   Denote $\im\varphi $
   by  $I$. Since $I$ is a qoutient of $E$ and a
   subsheaf of  $F$, it follows from \ref{l114} that
   $$\gamma (E)\le\gamma (I);\qquad \gamma (I)\le\gamma (F).$$
   But by assumption, $\gamma(E)>\gamma(F)$.
   This contradiction proves the lemma.

\TH{Lemma.}\label{l116} Let $E$ and $F$ be semistable sheaves
    with  $\gamma(E)=\gamma(F)$, let
   $\varphi :E\longrightarrow F$ be a nonzero morphism. Then

  a)  $E$ is stable $\quad\Longrightarrow\qquad\varphi $ is an injection;

  b) $F$ is stable $\quad\Longrightarrow\qquad\varphi $ is an epimorphism
  in general point.  \ETH

\PR By definition, put $I=\im\varphi $. As above,
   $$\gamma (E)\le\gamma (I);\qquad \gamma (I)\le\gamma (F),$$
   and by lemma condition,
    $$\gamma (E)=\gamma (I)=\gamma (F).$$

a) Assume that $E$ is stable. Using \ref{l114}, we get $E=I$, i.e.
    $\varphi $ is a monomorphism.

b) If $F$ is stable then $r(I)=r(F)$, i.e.  ${\rm coker}\varphi $ is
   a torsion sheaf.

\TH{Lemma.}\label{l117} A stable sheaf is simple, that is
    any its endomorphism has the form $\lambda\cdot{\rm id}$.
\ETH

\PR Let us take any nonzero $\varphi\!\in\!{\rm End}(E)$.
    Using the lemma \ref{l116}, we obtain that $\varphi $ is  isomorphism
    on an open set.
    Let $\lambda$ denote a eigenvalue of this isomorphism over some point
    $x\!\in\!X$. If the map
    $(\varphi -\lambda {\rm id}_E)\in {\rm End}(E)$ is nontrivial
    then it is injection.
    Thus there exists the monomorphism of determinants
    $$\Lambda^{r(E)}(\varphi -\lambda {\rm id}_E) :\det (E)\longrightarrow
    \det (E).$$
    Since the  map $(\varphi -\lambda {\rm id}_E)$ degenerates at
    the point  $x$, the corresponding determinant map at $x$ is equal to
    zero.  Hence, $\Lambda^{r(E)}(\varphi -\lambda {\rm id}_E)\equiv 0$.
    Therefore  $(\varphi -\lambda {\rm id}_E)$ has a nontrivial kernel.
    It is impossible, i.e.
    $(\varphi -\lambda {\rm id}_E)\equiv 0 .$  This completes the proof.

\TH{Lemma.}\label{l118} Let
$$0\longrightarrow E\stackrel{i}{\longrightarrow }G\longrightarrow F
   \longrightarrow 0,$$
   be an exact sequence of torsion free sheaves such that
   $\gamma(E)=\gamma(G)=\gamma(F)$. Then
 $G$ is semistable if and only if both $E$ and $F$ are semistable.
   In particular, for any finitedimensional vector space
   $V$ over $\CC$ and a divisor  $D$ on $X$ the bundle
   $V\otimes {\cal O}_X(D)$  is semistable.
\ETH

\PR Assume that $G$ is semistable. Since $E\subset G$, we see that
    any subsheaf $E'$ of $E$ is also a subsheaf of $G$. Hence,
     $\gamma (E')\leq \gamma (G)=\gamma (E)$,
     i.e. the sheaf $E$ is also semistable. Using
     (\ref{l112}), it can similarly by checked that the sheaf
      $F$ is also semistable.

    Now suppose the sheaves $E$ and $F$ are semistable and $G$ has
    torsion free quotient $$\varphi :G\longrightarrow Q\longrightarrow 0$$
    such that $r(Q)<r(G)$ and $\gamma (G)>\gamma (Q)$ (\ref{l113}). If
    $Q$ is not semistable then consider its quotient $Q'$ such that
    $r(Q')<r(Q)$ and $\gamma (Q)>\gamma (Q')$ , ets... Taking into account
    \ref{l111}, we can assume without loss of generality that $Q$ is
    semistable.  Thus it follows from $$\gamma (E)=\gamma (G)>\gamma (Q)$$
   that $\varphi \circ i$ is the zero map (see \ref{l115}). On the other
   hand, for the same reason, $${\rm Hom}(F,Q)=0$$.

   This contradiction concludes the proof.

\vspace{2cm}

\subsection{Harder-Narasimhan Filtration.}

    The aim of this section is the construction for a torsion free sheaf the
    well known canonical filtration which becomes trivial when
    the sheaf is semistable.
    Let us remember the main definition and notations.

   The record  $Gr(E)=(G_n,G_{n-1},\ldots,G_1)$ means that the sheaf
   $E$ has a filtration:
$$ 0=E_{n+1}\subset E_{n}\subset\cdots\subset E_2\subset E_1=E$$
   and $E_i/E_{i+1}=G_i$. The sheaves  $E_i$ are called  {\it terms of
   filtration\/} and    $G_i$  are {\it quotients of filtration\/}.
   Note that $G_n=E_n$ (since $E_{n+1}=0$).

   {\sc Definition.} A filtration of a torsion free sheaf $E$
       $$\qquad Gr(E)=(G_n,G_{n-1},\ldots,G_1)$$
   is called {\it Harder-Narasimhan filtration\/} if all quotients $G_i$
   are semistable and their slopes satisfy inequalities:
  $$\gamma (G_{i+1})>\gamma (G_i)\qquad (i=1,2,...,n-1).$$

\vspace{2ex}

  To construct this filtration, we need another one slope axiom and
  several lemmas.

\TH{Lemma.}\label{l121} Let $E$ be a sheaf without torsion on $X$ and
${\cal G}$ the set of all torsion free quotients of $E$. Then there is
  $\gamma_0$ such that $\gamma(G)\geq\gamma_0$ for each $G\in {\cal G}$.
\ETH

\PR Let us choose an ample divisor $A$ on $X$. Then it follows from
   the Serre theorem (\cite{CH}) that there exists a natural number
    $n$ such that  the sheaf $E(nA)$ is generated by global sections.
    Hence we have the exact sequence:
    $$0\longrightarrow F\longrightarrow {\rm H}^0(E(nA))\otimes {\cal O}
    \stackrel{can}{\longrightarrow }E(nA)\longrightarrow 0.$$
   Therefore   $E$ is a quotient of the semistable bundle
   $${\rm H}^0(E(nA))\otimes {\cal O}(-nA).$$
   If $G$ is a torsion free quotient of $E$ then there exists an
   epimorphism $${\rm H}^0(E(nA))\otimes {\cal O}(-nA)\longrightarrow
 G\longrightarrow 0.$$ Now the lemma follows from \ref{l113}.

\TH{Lemma.}\label{l122} Suppose a slope function $\gamma$
  satisfies the axiom:

  {\sc SLOPE.3.} Let  $\gamma _0$  be a value of the function
  $\gamma $ and $M=\{G_1,G_2,G_3,\ldots\quad\}$ an ordered set of sheaves
  without torsion with $r(G_i)\leq r$ for any $i$. Then the condition
  $$\gamma (G_i)>\gamma (G_{i+1})\geq \gamma_0\qquad\forall i=1,2,3,\ldots$$
  implies that  $M$ is finite.

  Then each torsion free sheaf $E$ has  the quotient
  $G$ with the minimal slope
  $\gamma(G)$. That is, for another torsion free quotient $Q$ of
  $E$, we have: $\gamma (Q)\ge\gamma (G)$.
\ETH

  The proof follows from {SLOPE.3} and the previous lemma.

\TH{Proposition.}\label{l123} If a slope function $\gamma $ satisfies
  the axioms {\sc SLOPE.1} - {\sc SLOPE.3} then any torsion free sheaf $E$
  has the Harder-Narasimhan filtration
  $$ Gr(E)=(G_n,G_{n-1},\ldots,G_1).$$
  Moreover, if $ Gr(E)=(G'_m,G'_{m-1},\ldots,G'_1)$ is another filtration
 with semistable quotients and the inequalities
  $\gamma(G'_i)>\gamma(G'_{i-1})$ for each $i=2,3,...,m$ then
   $m=n$ and $G'_i=G_i  \qquad\forall i$.
\ETH

  {\sc Proof of existence.}
  A semistable sheaf has  trivial filtration. Suppose $E$ is not
  semistable. Denote by $G_1$ the torsion free quotient of
  $E_1=E$ with the minimal $\gamma$-slope and the maximal rank.
  That is for another quotient  $Q$  without torsion of
   $E$ we have  $\gamma (Q)\geq \gamma (G_1)$, and the equality
   $\gamma (Q)= \gamma (G_1)$ implies   $r(Q)\le r(G_1)$.
   Let  $E_2$ be the corresponding subsheaf in $E$:
  $$0\longrightarrow E_2\longrightarrow E_1\longrightarrow G_1
  \longrightarrow 0.$$
  If $E_2$ is not semistable then let us choose the torsion free
  quotient $G_2$ of $E_2$ with the minimal $\gamma$-slope and the maximal
  rank.  Denote by $E_3$ the corresponding subsheaf in $E_2$,  etc.  Note
   that all $G_i$ are semistable by construction.  Let us check the
  inequality $\gamma (G_i)<\gamma (G_{i+1})$  with aid of the following
 commutative diagram:  $$\begin{array}{ccccccccc} &&0&&&&&&\\
    &&\uparrow&&&&&&\\
    &&G_{i+1}&&&&0&&\\
    &&\uparrow&&&&\uparrow&&\\
    0&\longrightarrow&E_{i+1}&\longrightarrow&E_i&\longrightarrow &G_i&
    \longrightarrow&0\\
    &&\uparrow&&\llap{$\scriptstyle id_E$}\uparrow&&\uparrow& &\\
    0&\longrightarrow&E_{i+2}&\longrightarrow&E_i&\longrightarrow &
    Q&\longrightarrow&0\\
    &&\uparrow&&&&\uparrow&&\\
    &&0&&&&G_{i+1}&&\\
    &&&&&&\uparrow&&\\
    &&&&&&0&&
   \end{array}
.$$
   We have that  $Q$ is torsion free quotient of $E_i$.
   It follows from $r(Q)>r(G_i)$ that  $\gamma (Q)>\gamma (G_i)$.
   Finally, if we recall the axiom {\sc SLOPE.1},   we get
 $\gamma (G_{i+1})>\gamma (G_i)$.  This concludes the  proof of existence.

 \TH{Lemma.}\label{l124} Let  $E,F$ be sheaves on $X$ and
  $Gr(E)=(G_n,G_{n-1},\ldots,G_1)$ a filtration of $E$. Then

  a) ${\rm Ext}^k(G_i,F)=0\quad\forall i\qquad\Longrightarrow\qquad
  {\rm Ext}^k(E,F)=0$;

  b) ${\rm Ext}^k(F,G_i)=0\quad\forall i\qquad\Longrightarrow\qquad
  {\rm Ext}^k(F,E)=0$.
\ETH

\PR
      Let us prove the second statement. The first one can be
      checked similarly.
      Consider the exact sequences:
    $$0\longrightarrow E_{i+1}\longrightarrow E_i\longrightarrow G_i
    \longrightarrow 0,$$
     where $E_i$ are the terms of the filtration.
     It follows from the corresponding long cohomology sequences that the
     triples
      $${\rm Ext}^k(F,E_{i+1})\longrightarrow {\rm Ext}^k(F,E_i)
      \longrightarrow {\rm Ext}^k(F,G_i)$$
      are exact  $\forall i=n-1,...,1$. Since $E_n=G_n$, by lemma
      conditions we have
      $${\rm Ext}^k(F,E_n)={\rm Ext}^k(F,E_{n-1})=\cdots=
      {\rm Ext}^k(F,E_1)={\rm Ext}^k(F,E)=0.$$
      This completes the proof.

\TH{Corollary.}\label{l125} Let  $Gr(E)=(G_n,G_{n-1},\ldots,G_1)$ be the
      Harder-Narasimhan filtration of a sheaf $E$ and let $F$ be
      a semistable sheaf. Then

     a)  $\gamma (F)<\gamma (G_1)\qquad\Longrightarrow\qquad
                                                {\rm Hom}(E,F)=0$;

     b) $\gamma (F)>\gamma (G_n)\qquad\Longrightarrow\qquad{\rm Hom}(F,E)=0$.
\ETH
    The proof follows easily from the lemmas  \ref{l124}, \ref{l115},
    and the definition of Harder-Narasimhan filtration.

\TH{Lemma.}\label{l126} If a sheaf $E$ has a filtration
   $Gr(E)=(G_n,G_{n-1},\ldots,G_1)$ then
   $G_n$ is a subsheaf of $E$ and $Gr(E/G_n)=(G_{n-1},\ldots,G_1)$.
\ETH

\PR Since the last quotient of the filtration consists with its last term,
    we get $G_n\subset E$ .
    Now this lemma is immediate if we consider the following
    commutative diagram:
   $$ \begin{array}{ccccccccc}
   &&0&&0&&0&&\\
   &&\uparrow&&\uparrow&&\uparrow&&\\
   0&\longrightarrow &E_{i+1}/G_n&\longrightarrow &E_i/G_n&\longrightarrow &
   G_i&\longrightarrow &0\\
   &&\uparrow&&\uparrow&&\uparrow&&\\
   0&\longrightarrow &E_{i+1}&\longrightarrow &E_i&\longrightarrow &
   G_i&\longrightarrow &0\\
   &&\uparrow&&\uparrow&&\uparrow&&\\
   0&\longrightarrow &G_n&\longrightarrow &G_n&\longrightarrow &0&&\\
   &&\uparrow&&\uparrow&&&&\\
   &&0&&0&&&&\\
      \end{array}
.$$

\vspace{1ex}

 {\sc Proof of the uniqueness of Harder-Narasimhan filtration.}

    Let $$Gr(E)=(G_n,G_{n-1},\ldots,G_1)=(G'_m,G'_{m-1},\ldots,G_1)$$
    be two Harder-Narasimhan filtrations.
    Suppose   $\gamma (G_1)\not=\gamma (G'_1)$.
    For example, $\gamma (G_1)>\gamma (G'_1)$.
    Then it follows from  the corollary \ref{l125} and  semistability of
     $G'_1$ that  ${\rm Hom}(E,G'_1)=0$.
    This contradicts an existence of an epimorphism: $E\longrightarrow
    G'_1 \longrightarrow 0$. In the same way, the equality $\gamma
    (G_n)=\gamma (G'_m)$ is proved.

   Denote by  $E'_i$ the terms of the second filtration. Let us show by
   induction on $i$ that   $G_n$ is a subsheaf in $E'_i$. For $i=1$,
   there is nothing to prove.

   By the induction hypothesis, we have the following commutative diagram:
   $$\begin{array}{ccccccccc}
   0&\longrightarrow &E'_{i+1}&\longrightarrow &E'_i&\longrightarrow &
   G'_i&\longrightarrow &0\\
   &&\uparrow&&\uparrow&&\llap{$\scriptstyle \varphi_i$}\uparrow&&\\
   &&0&\longrightarrow &G_n&\longrightarrow &
   G_n&\longrightarrow &0\\
   &&&&\uparrow&&&&\\
   &&&&0&&&&
      \end{array}
   .$$
   By the lemma about a snake,  $\ker \varphi _i\subset E'_{i+1}$.
   On the other hand, the slopes of semistable sheaves $G_n,G'_m$ and $G'_i$
   satisfy conditions:
   $\gamma (G_n)=\gamma (G'_m)>\gamma (G'_i)$ if $i<m$.
   Hence, $\varphi _i=0$ and $\ker \varphi _i=G_n$ (\ref{l115}).

   Thus, $G_n\subset E'_{i+1}$ for $i<m$. In particular, $G_n\subset G'_m$.

   In the same way, we obtain that  $G'_m\subset G_n$. Therefore,
   $G'_m=G_n$.

   It follows from the lemma \ref{l126} that
   $$Gr(E/G_n)=(G_{n-1},\ldots,G_1)=(G'_{m-1},\ldots,G'_1).$$
   Moreover, these are the  Harder-Narasimhan filtrations of $E/G_n$.
   Now the uniqueness of  the  Harder-Narasimhan filtration
   follows easily by induction on a rank of $E$.

\vspace{2cm}

\subsection{ Examples of Slopes and Kinds of Stability.}

   The reason of slope axioms are the following well known slopes.

  The slope of bundles on a curve: $\; \mu (E)=\frac{\deg E}{r(E)}$,
  where   ${\deg E}$ is the degree of the  bundle's determinant;

  the Mumford-Takemoto slope with respect to an ample divisor $A$
  on   $n$-dimensional manifold $X$:
   $\; \mu_A(E)=\frac{c_1(E)\cdot A^{n-1}}{r(E)}$;

  the Gieseker slope relative to an ample divisor $A$:
   $\; \gamma _A(E,n)=\frac{\chi(E(nA))}{r(E)}$.

  Let us check that these slopes and the slope  $\mu_H(E)=
  \frac{c_1(E)\cdot H^{n-1}}{r(E)}$, where $H$ is nef,
   really satisfy the slope
  axioms. By definition, a divisor $A$ is  {\it nef\/}
  if the number
$D\cdot A^{n-1} $ is greater then or equals 0 for any effective divisor
   $D$ on$X$.

  We see that all examples of slopes, except for
  $\gamma _A$,
  have the form   $\gamma=d/r$, where $d$ is an additive function of
   $K_0(X)$ to $\ZZ$ and $r$ is the  rank function.

\TH{Lemma.}\label{l131} The slope function $\gamma=d/r$
  defined before satisfies the axioms
   {\sc SLOPE.1} and  {\sc SLOPE.3}.

\ETH
\PR For any exact triple of torsion free sheaves
  $$0\longrightarrow F\longrightarrow E\longrightarrow G\longrightarrow 0$$
  we have that
  $\gamma (E)=\frac {d(F)+d(G)}{r(F)+r(G)}$. Note that the sign of
  the determinant
  $$\begin{array}{|cc|}
   d(F)&d(G)\\
   r(F)&r(G)
   \end{array}$$
   corresponds to the comparison sign between the fractions:
   $\frac {d(F)}{r(F)}\qquad\frac {d(G)}{r(G)}$. Besides,
  $$\begin{array}{|cc|}
   d(F)&d(G)\\
   r(F)&r(G)
   \end{array}\quad=\quad
   \begin{array}{|cc|}
   (d(F)+d(G))&d(G)\\
   (r(F)+r(G))&r(G)
   \end{array}.$$
  This implies that  $\gamma$ satisfies the first axiom.

  For checking the axiom {\sc SLOPE.3}, note that
 $|\gamma (G_1)-\gamma (G_2)|\geq 1/r^2$ if the ranks of torsion free
  sheaves $G_1$ and $G_2$ are less than or equal to $r$.

  \TH{Corollary.}\label{l132} Let  $\gamma =(\gamma _0,\gamma _1,
     \ldots,\gamma _n)$ be a vector function of $K_0(X)$ such that each
      $\gamma _i $ has the form $d_i/r$, where $d_i$ is an additive
      function of $K_0(X)$ to $\ZZ$ and $r$ is the rank function.
      If values of  $\gamma $ are lexicographic compared then $\gamma$
      satisfies the axioms
     {\sc SLOPE.1} and  {\sc SLOPE.3}.
\ETH

\vspace{2ex}

  As for the Gieseker slope  $\gamma _A$, it is a polynomial of the degree
  $\dim X$ with rational coefficients.
  So far as the inequality
   $\gamma _A(E,n)>\gamma _A(F,n)$
  holds true if it holds for sufficiently large $n$, then the comparison
  $\gamma_A$-slopes is equivalent to lexicographic ordering of
  the coefficients of  the polynomials.

 The Hilbert polynomial $\chi(E(nA))$ is an additive function. Hence
  the Geaseker slope satisfies the axiom {\sc SLOPE.1} ( see the proof of
 lemma \ref{l131}).

  For checking {\sc SLOPE.3} note that by
  Hirzebruch-Riemann-Roch theorem (see. \cite{CH}), the Euler
  characteristic of sheaf on a smooth manifold can be calculated in the
  following way:  \begin{equation}\label{equ1} \chi(E)=\deg(ch(E)\cdot
      td(T_X))_n, \end{equation} where

  $\deg(...)_n$ means a degree $n$ component in the cohomology ring of
   $X$ ($H^*(X,\QQ)$);

  $T_X$ is the tangent bundle of  $X$;

  $$ch(E)=r+c_1+\frac{1}{2}(c_1^2-2c_2)+\frac{1}{6}(c_1^3-3c_1c_2+3c_3)+
  \frac{1}{24}(c_1^4-4c_1^2c_2+4c_1c_3+2c_2^2-4c_4)+\cdots; $$

  $$td(E)=1+c_1/2+\frac{1}{12}(c_1^2+c_2)+\frac{1}{24}(c_1c_2)-
 \frac{1}{720}(c_1^4-4c_1^2c_2-3c_2^2-c_2c_3+c_4)+\cdots .$$
 ($c_i$ are the Chern classes of a sheaf  $E$).

  This yields that the denominators of the coefficients of
  the Hilbert polynomial
   $\chi(E(nA))$ do not depend of $E$.  After some proof modifications
   of the lemma \ref{l131}, it is easily shown that the Gieseker slope
    $\gamma _A(E)$ satisfies the axiom {\sc SLOPE.3}.

  All examples of slopes satisfy the axiom {\sc SLOPE.2}
  in the different degree.

\TH{Lemma.}\label{l133} a) For any pair of torsion free sheaves
   $F\subset E$ with the same rank on a manifold $X$
   and any nef divisor $H$ the following inequality holds
   $ \mu_H(F)\le \mu_H(E)$.
   Moreover, in this case the equality $\mu_H(F)=\mu_H(E)$ is possible
   only if
    $${\rm codim\,supp}(E/F)\ge 1.$$
    Provided a slope function satisfies this reduction of axiom {\sc SLOPE.2},
   we shall call it the {\it weak\/} slope;

   b)  for any pair of torsion free sheaves
   $F\subset E$ with the same rank on a manifold $X$
   and any ample divisor $A$ the following inequality holds
   $ \mu_A(F)\le \mu_A(E)$.
   Moreover, in this case the equality is possible
   only if
    $${\rm codim\,supp}(E/F)\ge 2.$$
   Provided a slope function satisfies    this reduction of axiom {\sc
SLOPE.2},
   we shall call it the {\it Mumford-Takemoto\/} slope;

   c)  for any pair of torsion free sheaves
   $F\subset E$ with the same rank on a manifold $X$
   and any ample divisor $A$ the following inequality holds
   $ \gamma_A(F)\le \gamma_A(E)$.
   Moreover, in this case the equality of slopes is equivalent to
    $$E=F.$$
   Provided a slope function satisfies this reduction of axiom {\sc SLOPE.2},
   we shall call it the {\it Gieseker\/} slope;

   d)  the slope $\mu$ of bundles on a curve is the Gieseker slope.

\ETH

\PR  The number $c_1(E)\cdot D^{n-1}$, which is determined by a sheaf $E$
  on a $n$-dimensional manifold and a divisor
   $D$, is called the degree of sheaf with respect to $D$ and is denoted by
  $\deg_D(E)$.

   Since the ranks of sheaves  $E$ and  $F$ coincide, we see that the
   comparison of their slopes is equivalent to the comparison of the degree
   $\deg_D$ and the quotient $Q=E/F$
   has the zero rank. Hence
    $c_1(Q)=c_1(E)-c_1(F)$
  is effective or zero divisor.

   By the definition of nef divisor,
    $\deg_H(Q)\ge 0$.
    This proves the first statement of lemma.

   If $A$ is ample and  $c_1(Q)\not= 0$ then
   the Nakai-Moyshezon criterion
   (\cite{CH}) implies that $\deg_A(Q)>0$. This yields the second
   statement of lemma.

   If $A$ is ample then by the Serre theorem
   (\cite{CH}) for any nonzero sheaf  $Q$ and for sufficiently large
 $n$ $\chi(Q(nA))>0$. Therefore the third point of lemma  also
 holds.

   Finally, the degree of the effective divisor  $c_1(Q)$ on a curve
   is nonnegative and it is equal to zero only if
   $c_1(Q)=0$. This completes the proof.

\vspace{2ex}

  The more precise conditions of {\sc SLOPE.2}
  allow to formulate the following statement that is stronger the lemma
  \ref{l116}

\TH{Lemma.}\label{l134} a) Let $E$ and $F$ be a semistable sheaves
  with respect to Mumford-Takemoto slope $\gamma$,
  $\gamma (E)=\gamma (F),$ and  $F$ stable. Then
  the cokernal of any nonzero morphism $\varphi :E\longrightarrow F$
  has a support $C$ such that ${\rm codim} C \ge 2$.
  In particular,   $\varphi $ is an epimorphism if  $E$ is locally free.

  b) Let $E$ and $F$ be semistable sheaves with respect to Gieseker slope
  $\gamma$, $\gamma (E)=\gamma (F)$
  and  $F$ stable. Then any nonzero map of $E$ to $F$ is an epimorphism.
\ETH

  This lemma can be proved in the same way as
  \ref{l116}. Nevertheless, let us recall that
   ${\rm Ext}^1(Q,E)=0$ if $E$ is locally free and
  ${\rm codim\ supp}(Q)\ge 2$.

\vspace{2ex}

  Let us remark that the slope   $\gamma=
  \mu_A(E)=(c_1(E)\cdot A^{n-1})/r(E)$
  has the following property:

 {\sc SLOPE.4.} For any torsion free sheaf $E$ and
 a divisor $D$  the equalities
   $$\gamma (E^*)=-\gamma (E),\qquad
   \gamma (E(D))=\gamma (E)+\gamma ({\cal O}(D))$$ are true.

 \TH{Lemma.}\label{l135} Assume that slope function
   $\gamma$ satisfies the axiom   {\sc SLOPE.4}; then a torsion free sheaf
  $E$ is (semi)stable if and only if  $E(D)$ is (semi)stable; and
  the $\gamma$-(semi)stability of a reflexive sheaf $E$
  ($E^(**)=E$) is equivalent to the
  $\gamma$-semistability of the dual sheaf $E^*$.
\ETH

\PR Consider any subsheaf $F$ in $E(D)$. Hence, $F(-D)\subset E$ and
 $\gamma (F(-D))\le\gamma (E)$ if $E$ is semistable. Therefore,
$$\gamma (F)=\gamma (F(-D))+\gamma ({\cal O}(D))\le\gamma(E)+
\gamma({\cal O}(D))=\gamma (E(D)).$$
 That is, the sheaf   $E(D)$ is semistable.

  Now suppose that $E$ is reflexive. It is sufficient to prove that the
  semistability of
  $E^*$ implies the semistability of $E$.
  Denote by  $G$  any torsion free quotient of $E$. Hence,
  $G^*\subset E^*$.
  Therefore,
  $\gamma(G^*)\leq \gamma(E^*)$. Using {\sc SLOPE.4},
  we obtain that $\gamma (G)\geq \gamma (E)$. This implies the
  semistability of $E$.

\vspace{1ex}

   It is useful to remark that besides the canonical Harder-Narasimhan
   filtration, each semistable sheaf has Jordan-Holder filtration:

\TH{Proposition.} Any  $\gamma $-semistable sheaf  $E$ has the filtration
    $$Gr(E)=(G_n,G_{n-1},\ldots,G_1)$$ with stable quotients and the
    equalities:
    $\gamma (G_i)=\gamma (E)$.
\ETH

   We do not prove this proposition. But we constructs more
   exact filtration.

\TH{Proposition.}\label{l137} Any sheaf $E$ semistable with respect to
   Gieseker slope (see \ref{l123}) has the filtration with
   isotypic quotients:
   $$Gr(E)=(G_n,G_{n-1},\ldots,G_1),$$
   where each of $G_i$  has the filtration with isomorphic quotients:
   $$Gr(G_i)=(Q_i,Q_i,\ldots,Q_i)\qquad(\gamma (Q_i)=\gamma (G_i)=
   \gamma (E)).$$
   Moreover, this filtration can be constructed in the such way that
 $${\rm Hom}(E_i,G_{i-1})={\rm Hom}(G_i,G_{i-1})={\rm Hom}(G_{i-1},G_i)=0,$$
 where $E_i$ are the filtration terms.
\ETH

\PR  For a stable sheaf this filtration is trivial. If a sheaf
     $E=E_1$ is semistable then it has destabilizing torsion free quotient
     $Q\quad(\gamma (Q)=\gamma (E_1))$. From all such quotients let us
     choose a sheaf $Q_1$ with the minimal rank. By choice, it is stable.
     Let  $E^1_1$ be the corresponding subsheaf.
     It follows from the exact sequence:
     $$0\longrightarrow E^1_1\longrightarrow E_1\longrightarrow Q_1
     \longrightarrow 0,$$        and the equality
      $\gamma (Q_1)=\gamma (E_1)$
     that $E^1_1$ is semistable and
     $\gamma (E^1_1)=\gamma (E_1)$  (see {\sc SLOPE.1},
      \ref{l118}).
     If  ${\rm Hom}(E^1_1,Q_1)=0$ then  $E_2=E^1_1$ is the second term
     of filtration and
    $G_1=Q_1$ is the first quotient of it.

     Conversely, there exists an epimorphism;
     $E^1_1\longrightarrow Q_1\longrightarrow 0$ (\ref{l134}).
     Denote by $E^2_1$ the kernel of this epimorphism.

     Continuing this procedure, we obtain the semistable subsheaf
     $E^k_1$ such that
     $${\rm Hom}(E^k_1,Q_1)=0.$$ By definition, put
     $E_2=E^k_1$ and  $G_1=E_1/E_2$.
     From this construction it follows that  $G_1$ and $E_2$ are semistable,
     $\gamma (E_2)=\gamma (G_1)=\gamma (E_1),\qquad
     Gr(G_1)=(Q_1,Q_1,\ldots,Q_1)$
     and  ${\rm Hom}(E_2,Q_1)=0$. Now using the lemma \ref{l124}, we
     obtain, ${\rm Hom}(E_2,G_1)=0$.

     By the induction hypothesis, we can assume that $E_2$ has the
     filtration with isotypic quotients: $Gr(E_2)=(G_n,G_{n-1},\ldots,
     G_2)$. Let us show that the filtration
     $$Gr(E)=(G_n,G_{n-1},\ldots,G_2,G_1)$$ satisfies the proposition
     conditions.

    It remains to check that
         ${\rm Hom}(G_2,G_1)={\rm Hom}(G_1,G_2)=0$. Since
     ${\rm Hom}(E_2,G_1)=0$ and there exists an epimorphism
     $E_2\longrightarrow G_2\longrightarrow 0$,
     the equality ${\rm Hom}(G_2,G_1)=0$ is trivial.

     Suppose there exists a nonzero morphism
     $G_1\longrightarrow G_2$. Let us recall that
     $Gr(G_i)=(Q_i,Q_i,\ldots,Q_i)$ .Therefore by  \ref{l124},
     ${\rm Hom}(Q_1,G_2)\not= 0$. Hence , there is a nonzero map
     $\varphi :Q_1\longrightarrow Q_2$. It follows from
     (\ref{l116} and \ref{l134}) that $\varphi$ is an isomorphism.
     It implies  that ${\rm Hom}(Q_2,G_1)\not= 0$.
     But $Q_2$ is a quotient of
     $G_2$, and $Q_1$ is a subsheaf of  $G_1$. Thus,
     ${\rm Hom}(G_2,G_1)\not= 0$. This contradiction concludes the proof.

%% FOLLOWING LINE CANNOT BE BROKEN BEFORE 80 CHAR
%%%%%%%%%%%%%%%%%%%%%%%%%%%%%%%%%%%%%%%%%%%%%%%%%%%%%%%%%%%%%%%%%%%%%%%%%%%%%%%%%%%%
\section{Rigid Sheaves.}
\subsection{Preliminary information.}

    We shall study sheaves on a smooth projective surface
    $S$ over $\CC$ such that $h^1({\cal O}_S)=0$ and the anticanonical
    class $H=-K_S$ has not base components.
    Note that from the last condition it follows that
     $H$ is nef.
    In reality, suppose the cup product $H\cdot C$ is negative for some
    curve $C$; then $H$ and $C$ have a common base component.

  It is known, if $S$ is a smooth projective surface over an algebraically
  closed field with nef anticanonical class then we have one of the
  following cases:
  \begin{enumerate}
  \item
  $K_S=0$;
  \item
  $S\cong {\BBB P}({\cal O}_{{\BBB P}^1} \oplus {\cal O}_{{\BBB P}^1}(2))$;
  \item
  $S\cong {\BBB P}(F)$, where $F$ is a rank 2 vector bundle on an elliptic
  curve which is an extension of degree zero line bundles;
  \item
  $S\cong  {\BBB P}^2$  or $ {\BBB P}^1\times {\BBB P}^1$;
  \item
  $S$ is obtained from  ${\BBB P}^2$ by successively blowing up at most
  nine points.
  \end{enumerate}

    I want to note that this class of surfaces contains the surfaces such
    that they are obtained from the singular Fano surfaces (taken from
    the notes of Batyrev in the paper \cite{B}) by blowing up singular
    points.

    Let us recall the general facts, which will be needed in this text.

\TH{Theorem.}\label{l211} (The Riemann-Roch formula for surfaces.)
        The Euler characteristic of two coherent sheaves
        $E$ and $F$ on a smooth projective surface $X$
        is calculated by the following formula:
    $$\chi(E,F)=r(E)r(F)\Bigl(\chi({\cal O}_X)+\frac{1}{2}(\mu_H(F)-\mu_H(E))
    +q(F)+q(E)-\frac{(c_1(E)\cdot c_1(F)}{r(E)r(F)}\Bigr),$$
$$\mbox{\rm where}\qquad \mu_H(E)=\frac{1}{r(E)}(-K_X\cdot c_1(E)),\quad
    q(E)=\frac{c^2_1(E)-2c_2(E)}{2r(E)}.$$
\ETH

    This theorem can be proved by direct calculation.
    (see (\ref{equ1})).

\vspace{2ex}

   Note that in our case we have  $\chi({\cal O}_S)=1$.

\TH{Corollary.}\label{l212} Let $E,F$ be a sheaves on a smooth
   projective regular ($h^1({\cal O}_S)=0$)
   surface with  $\chi(E,E)=\chi(F,F)=1$ ; then
   $$q(E)=\frac{1}{2}\Bigl(\frac{c^2_1(E)+1}{r^2(E)}-1\Bigr)\qquad
   \mbox{\rm and}$$
   $$\chi(E,F)=\frac{r(E)r(F)}{2}\Bigl(\mu_H(F)-\mu_H(E)+
   \frac{1}{r^2(E)}+\frac{1}{r^2(F)}+(\frac{c_1(F)}{r(F)}-
   \frac{c_1(E)}{r(E)})^2\Bigr).$$
\ETH

\TH{Theorem.}(The Serre duality.) For any coherent sheaves
  $E$ and $F$ on a smooth projective surface $X$ the following equality
$${\rm Ext}^{k}(E,F)^*\cong {\rm Ext}^{2-k}(F,E(K_X))$$
  holds.
\ETH

   The proof is contained in \cite{CH}.

\vspace{1ex}

\TH{Lemma.}\label{l214} (Mukai.) Let $X$ be a smooth projective surface.

1. For any torsion free sheaf $E$ on $X$ we have
$$h^1(E,E)\geq h^1(E^{**},E^{**})+2length(E^{**}/E).$$
2.a) Suppose the sheaves $G_1$ and $G_2$ on $X$ from the exact sequence
$$0\longrightarrow G_2\longrightarrow E\longrightarrow G_1
\longrightarrow 0,$$
  satisfy the conditions:
 ${\rm Hom}(G_2,G_1)={\rm Ext}^2(G_1,G_2)=0$; then
$$h^1(E,E)\geq h^1(G_1,G_1)+h^1(G_2,G_2).$$

b) If besides, $h^1(E,E)=0$; then
$$h^0(E,E)=h^0(G_1,G_1)+h^0(G_2,G_2)+\chi(G_1,G_2)$$
$$h^2(E,E)=h^2(G_1,G_1)+h^2(G_2,G_2)+\chi(G_2,G_1).$$
\ETH
  This lemma follows from the spectral sequence  associated with
  the exact triple.
  Besides, the proof is contained in \cite{MUK} and
 \cite{OK}.

%% FOLLOWING LINE CANNOT BE BROKEN BEFORE 80 CHAR
%%%%%%%%%%%%%%%%%%%%%%%%%%%%%%%%%%%%%%%%%%%%%%%%%%%%%%%%%%%%%%%%%%%%%%%%%%%%%%%%%%%%%
\subsection{Exceptional sheaves.}

{\sc Definition.} A sheaf $E$ on a manifold $X$ is called  {\it rigid\/}
   whenever $${\rm Ext}^1(E,E)=0.$$

\vspace{2ex}

The most trivial rigid sheaves are exceptional.

\vspace{2ex}

{\sc Definition.} A sheaf  $E$ on a manifold is called
    {\it exceptional\/},
  provided ${\rm Ext}^0(E,E)=\CC$ and ${\rm Ext}^i(E,E)=0\quad\forall i>0.$

\vspace{2ex}

  Using methods of S.~Mukai (\cite{MUK}),
  A.~Gorodentsev (\cite{AG}),
  D.~Orlov  (\cite{OK}) and S.~Zube (\cite{ZN}) we get the starting
  information  about a structure of rigid and exceptional sheaves.

\TH{Lemma.}\label{l221} A rigid sheaf without torsion on a smooth projective
   surface is locally free.
\ETH

  This lemma follows from Mukai lemma (\ref{l214}).

\vspace{2ex}

   Let us recall that  $S$ is a smooth projective surface over
   $\CC$ with the anticanonical class $H=-K_S$ without base components

   Let  $G$ be a sheaf on $S$. Denote by $TG$ its torsion subsheaf and by
   $T^0G$ the subsheaf in $TG$ such that $T^1G=TG/T^0G$ has not a torsion
   with 0-dimensional support.

\TH{Lemma.}\label{l222}(Gorodentsev-Orlov.) Any sheaves  $G$ and $F$
    on the surface $S$ satisfy the following conditions:

   a) the inequality
   $$h^0(F,G)\ge h^2(G,F)$$
   holds whenever the support of
   $T^0G$ has not common points with the base set of anticanonical linear
   system  $|H|$;

    b) the inequality
   $$h^0(G,G)>h^2(G,G)$$
    holds provided there exists a curve   $D\in |H|$ such that
    $D\cap{\rm supp}G\not=\emptyset$. In particular, this inequality is
    satisfied whenever     $r(G)>0$.
\ETH

\PR  a) For some section $\varphi \!\in\!H^0({\cal O}_S(H))$ let us consider
     the exact sequence:
$$0\longrightarrow {\cal O}_S\stackrel{\varphi }{\longrightarrow}
{\cal O}_S(H)\longrightarrow {\cal O}_S(H)|_D\longrightarrow 0.$$
     We apply the functor ${\cal H}om (\cdot,G)$ to it to obtain
$$0\longrightarrow{\cal H}om({\cal O}_S(H)|_D,G)\longrightarrow G(K_S)
\longrightarrow
    G\longrightarrow {\cal E}xt^1({\cal O}_S(H)|_D,G)\longrightarrow 0$$
\begin{center}
(${\cal H}om({\cal O}_S(H),G)\cong {\cal O}_S(H)^*\otimes G=G(-H)$ and
$H=-K_S$).
\end{center}
   Since the quotient  $G'=G/TG$ has not torsion, we get
    ${\cal H}om({\cal O}_S(H)|_D,G')=0$. Hence it follows from the exact
    sequence
$$0\longrightarrow TG\longrightarrow G\longrightarrow G'\longrightarrow 0$$
   that
$${\cal H}om({\cal O}_S(H)|_D,TG)={\cal H}om({\cal O}_S(H)|_D,G).$$

   The support of $T^1G$ is a curve. By the lemma conditions, the support
   of  $T^0G$ has the empty intersection with the base set of
    $|H|$. On the other hand,  $|H|$ has not base component.
    Thus there exists a curve   $D\!\in\!|H|$,
    such that
$$\dim(D\cap {\rm supp}T^1G)< 1,\qquad\qquad D\cap {\rm supp}T^0G=
\emptyset.$$
    Therefore,
$${\cal H}om({\cal O}_S(H)|_D,T^1G)={\cal H}om({\cal O}_S(H)|_D,T^0G)=0.$$
    From the exact sequence
$$0\longrightarrow T^0G\longrightarrow TG\longrightarrow T^1G
\longrightarrow 0$$
  and the previous equalities it follows that
$${\cal H}om({\cal O}_S(H)|_D,G)={\cal H}om({\cal O}_S(H)|_D,TG)=0,$$
   that is we have the exact triple:
$$0\longrightarrow G(K_S)\longrightarrow
    G\longrightarrow {\cal E}xt^1({\cal O}_S(H)|_D,G)\longrightarrow 0.$$
   It implies that  ${\rm Hom}(F,G(K_S))\subset {\rm Hom}(F,G)$  and
$h^0(F,G(K_S))\leq h^0(F,G)$.

  Now the first lemma statement follows from the Serre duality.

\vspace{1ex}

b) Note that ${\cal E}xt^1({\cal O}(H)|_D,G)\cong G\otimes {\cal O}_D$.
   Hence the previous exact sequence has the form:
  $$0\longrightarrow G(K_S)\longrightarrow G\longrightarrow
                                   G\otimes {\cal O}_D\longrightarrow 0.$$
   By assumption, there exists a curve  $D$ such that
  $G\otimes {\cal O}_D\not= 0$.
  Therefore, the map ${\rm Hom}(G,G)\longrightarrow {\rm Hom}(G,G
  \otimes {\cal O}_D)$
   is nonzero. Thus we obtain the strong inequality:
   $$h^0(G,G)>h^2(G,G).$$

\TH{Corollary.}\label{l223} Let  $G$ be a rigid sheaf on $S$;
   then  its torsion subsheaf $TG$ and torsion free quotient
    $G'=G/TG$ are rigid sheaves. Moreover,
    $T^0G=0$.
\ETH

\PR Obviously, ${\rm Hom}(TG,G')=0$ and $T^0G'=0$.
    By the previous lemma,    $h^0(TG,G')\geq h^2(G',TG)$.
    We apply the Mukai lemma (\ref{l214}) to the exact triple
 $$0\longrightarrow TG\longrightarrow G\longrightarrow G'\longrightarrow 0$$
    to obtain
    $$h^1(G,G)\geq h^1(TG,TG)+h^1(G',G'),$$
    i.e. $TG$ and $G'$ are rigid as well.

    Suppose $T^0G\not= 0$. Applying the Mukai lemma to the exact sequence:
    $$0\longrightarrow T^0G\longrightarrow TG\longrightarrow T^1G
                                                      \longrightarrow 0,$$
    we see that $T^0G$ is rigid.
    This contradicts the following equality:
         $$h^1(T^0G,T^0G)=2length(T^0G)\not= 0.$$

\TH{Lemma.}\label{l224} Suppose $E$ is an exceptional torsion sheaf on $S$;
  then $c_1^2(E)=-1$. Furthermore,

   either $E={\cal O}_e(d)$ , where  $e$ is  some irreducible rational curve
   with   $e^2=-1$

   or one of the components of the support of $E$ has zero cup product with
   $K_S$.
\ETH

\PR Using the Riemann-Roch formula for surfaces (\ref{l211}), we get
$$\chi(E,E)=r^2+(r-1)c_1^2-2rc_2.$$
   On the other hand,  since $E$ is exceptional, we have $\chi(E,E)=1$.
  Besides, from the lemma conditions it follows that $r(E)=0$. Therefore,
 $c_1^2(E)=-1$.

  Denote by  $C$ the support of the sheaf $E$. Since $E$ is rigid, we
  obtain $T^0E=0$ (\ref{l223}). Suppose $C$ is irreducible; then
   $E$ is a vector bundle on it.
    Let us denote the rank of this bundle by $r_C$.
   Hence, $c_1(E)=r_CC$ and  $c_1(E)^2=r_C^2C^2=-1$. Thus, $r_C=1$.

   By adjunction
  $(2g-2=C\cdot(C+K_S))$,
   $$2g=C\cdot K_S+1\leq 1\quad$$
  ($-K_S$ is nef). That is $C$ is a rational curve $e$
 with  $e^2=-1$, and $E={\cal O}_e(d)$ (as a line bundle on a
  projective line).

  Now suppose that $C$ is reducible and $C_0$ is one of its irreducible
  component. Consider the sequence of restriction to  $C_0$:
$$0\longrightarrow E_1\longrightarrow E\longrightarrow E_0
    \longrightarrow 0,$$
  where ${\rm supp} E_0=C_0$ and ${\rm supp} E_1=C_1$. Since $T^0E=0$, we
  get  $T^0E_0=T^0E_1=0$.
  Hence, ${\rm Hom}(E_1,E_0)={\rm Ext}^2(E_0,E_1)=0$.

  Now the application of the Mukai lemma (\ref{l214}) yields
  $$h^0(E,E)=h^0(E_0,E_0)+h^0(E_1,E_1)+\chi(E_0,E_1),$$
  $$h^2(E,E)=h^2(E_0,E_0)+h^2(E_1,E_1)+\chi(E_1,E_0).$$
   Since $E$ is exceptional, we get
   $$1=h^0(E_0,E_0)-h^2(E_0,E_0)+\chi(E_0,E_1)+
       h^0(E_1,E_1)-h^2(E_1,E_1)-\chi(E_1,E_0).$$
    By the Riemann-Roch theorem for exceptional sheaves (\ref{l212}),
    $$\chi(E_0,E_1)-\chi(E_1,E_0)=
    H\cdot \Bigl(r(E_0)c_1(E_1)-r(E_1)c_1(E_0)\Bigr)=0$$
    ( $r(E_0)=r(E_1)=0$).
        Finally, we obtain
          $$1=h^0(E_0,E_0)-h^2(E_0,E_0)+h^0(E_1,E_1)-h^2(E_1,E_1).$$

   To conclude the proof, it remains to note that the last equality
   is possible only if either
   $C_0$ or $C_1$ has the zero cup product with $H$ .

\TH{Lemma.}\label{l225} Suppose $E$ is an exceptional sheaf on $S$; then
  the support of its torsion has zero cup product with $K_S$.
\ETH

\PR It follows from corollary \ref{l223} that $TE$ and $E'=E/TE$ are rigid as
  well. The application Mukai lemma to the exact sequence
  $$0\longrightarrow TE\longrightarrow E\longrightarrow E'\longrightarrow
  0$$
  yields
  $$h^0(E,E)=h^0(TE,TE)+h^0(E',E')+\chi(E',TE),$$
  $$h^2(E,E)=h^2(TE,TE)+h^2(E',E')+\chi(TE,E').$$
  Using the equality  $h^0(E,E)-h^2(E,E)=1$, we obtain
  $$ 1=h^0(TE,TE)-h^2(TE,TE)+h^0(E',E')-h^2(E',E')+
  \chi(E',TE)-\chi(TE,E').$$
  By lemma \ref{l222}, $h^0(E',E')-h^2(E',E')\geq 1$ and
                     $h^0(TE,TE)-h^2(TE,TE)\geq 0$.
  Thus we have
  $$0\geq \chi(E',TE)-\chi(TE,E')=H\cdot(r(E')c_1(TE)-r(TE)c_1(E'))=
  r(E')H\cdot c_1(TE).$$
  On the other hand, since $H\;$ is nef, we get
  $\;H\cdot c_1(TE)\geq 0.$
  This completes the proof.

\vspace{2ex}

  Combining \ref{l221}, \ref{l224} and  \ref{l225} , we can formulate the
  following proposition:

\TH{Proposition.}\label{l226} Suppose $E$ is an exceptional
 sheaf on $S$; then
  we have one of the following cases:

  1) $E$ is locally free;

  2) $E$ has a torsion subsheaf such that $({\rm supp}TE)\cdot K_S=0$;

  3) $E\cong {\cal O}_e(d)$ for some rational curve $e$ with $e^2=-1$;

  4) $r(E)=0$ and the support of $E$ contains an irreducible component $C_0$
  such that $C_0\cdot K_S=0$.
   \ETH

\TH{Corollary.}\label{l227}(Orlov.) If $-K_S$ is ample ($S$ is the
  Del Pezzo surface) then exceptional sheaf on $S$ either is locally free
  or has the form ${\cal O}_e(d)$
for some rational curve $e$ with $e^2=-1$.
\ETH

\vspace{2ex}
   Now let us prove the stability of exceptional bundles on $S$ with
   respect to the anticanonical class $H=-K_S$.

\TH{Lemma.}\label{l228}(S.~Zube) Let $D$ be a smooth elliptic curve from
  $|H|$ and let $E$ be an exceptional bundle on $S$. Then the restriction
of $E$ to $D$ ( $E'=E|_D$ )is a simple bundle, i.e.
  ${\rm Ext}^0(E',E')=\CC$.
\ETH
\PR Consider the exact sequence
$$0\longrightarrow E^*\otimes E(K_S)\longrightarrow
E^*\otimes E\longrightarrow (E^*\otimes E)|_D\longrightarrow 0.$$
  By Serre duality,
  $${\rm Ext}^k(E,E(K_S))^*\cong {\rm Ext}^{2-k}(E,E).$$
  Since $E$ is exceptional, we obtain
$${\rm Ext}^0(E,E)=\CC,\quad{\rm Ext}^1(E,E)={\rm Ext}^2(E,E)=0.$$
  Therefore,  the cohomology table associated with the exact
  sequence has the form:
$$\begin{array}{|c|ccccc|}
\hline
k&{\rm Ext}^k(E,E(K_S))&
\rightarrow&{\rm Ext}^k(E,E)&
\rightarrow&{\rm Ext}^k(E',E')\\
\hline
0& 0 &&\CC&&?  \\
1& 0 && 0 &&?  \\
2& \CC && 0 &&?\\
\hline
\end{array}\; . $$
   It implies that ${\rm Ext}^0(E',E')=\CC.$

\TH{Corollary.}\label{l229} Any exceptional bundle $E$ on $S$ is stable
   with respect to the slope $\mu_H=(H\cdot c_1(E))/r(E),$ where $H=-K_S$.
\ETH

\PR  By the Zube lemma the restriction  of $E$ to an elliptic curve
     $D\in |-K_S|$ is simple. It is known that  simple bundles
     on an elliptic curve are stable with respect to the slope
 $\mu(E)=\frac{\deg(E)}{r(E)}$. Besides,
  $\mu_H(E)=\mu(E')$, where $E'=E|_D$.  Now suppose $F$ is a subsheaf of $E$
  such that $r(F)<r(E)$ and $\mu_H(F)>\mu_H(E)$. Without loss of generality
  we can assume that $F'=F|_D$ is locally free. Thus, $\mu(F')>\mu(E')$.
  This contradicts  stability of $E'$. The corollary is proved.

%% FOLLOWING LINE CANNOT BE BROKEN BEFORE 80 CHAR
%%%%%%%%%%%%%%%%%%%%%%%%%%%%%%%%%%%%%%%%%%%%%%%%%%%%%%%%%%%%%%%%%%%%%%%%%%%%%%%%%%%%%%%%
\subsection{Exceptional Collections.}

  The main results about rigid and superrigid sheaves are formulated in
  terms of exceptional collections. The aim of this section is to study
  these collections on the surface $S$.

{\sc Definition.} An ordered collection $\VEC E,n$  of exceptional sheaves
  is called
  {\it exceptional\/} whenever
   $$
  {\rm Ext}^k(E_i,E_j)=0\qquad\mbox{for }i>j\qquad\mbox{ and}\quad
  \forall k=0,1,2.$$
   An exceptional collection $(E,F)$ is an
{\it exceptional pair\/}.

\vspace{1ex}

   By definition, an ordered collection is exceptional if and only if
   each its pair is exceptional. Thus we shall study exceptional pairs on
   $S$.

   Suppose $(E,F)$ is an exceptional pair on  Del Pezzo surface.
   It is known that  then we have one of the following cases:

   a pair $(E,F)$ has the type {\it hom}  (or in other words $(E,F)$
   is a  $hom$-pair), that is
   $${\rm Ext}^i(E,F)=0\quad\mbox{for }\
 i=1,2 \qquad\mbox{and}\quad {\rm Hom}(E,F)\not= 0;$$

   a pair $(E,F)$ has the type {\it ext}  (or in other words $(E,F)$
   is a  $ext$-pair), that is
      $${\rm Ext}^i(E,F)=0\quad\mbox{for }\
 i=0,2 \qquad\mbox{and}\quad {\rm Ext}^1(E,F)\not= 0;$$

   a pair $(E,F)$ has the type {\it zero} (or in other words $(E,F)$
   is a  zero-pair), that is
      $${\rm Ext}^i(E,F)=0\quad\mbox{for }\
 i=0,1,2.$$

  There appear exceptional pairs of a new type in our surfaces.

{\sc Definition.} An exceptional pair $(E,F)$ is called {\it singular\/}
    if  $${\rm Ext}^i(E,F)\not= 0\quad\mbox{for}\; i=0,1\quad\mbox{ and}
\qquad{\rm Ext}^2(E,F)=0.$$

\TH{Proposition.}\label{l231} Let $(E,F)$ be an exceptional pair of bundles
   on the surface $S$; then we have one of the following cases:
$$\begin{array}{clcc}
\mbox{a)}&
$(E,F)$\ \mbox{ is the {\rm hom}-pair}&
 \qquad\Longleftrightarrow&\qquad \mu_H(E)<\mu_H(F);\\
\mbox{b)}& $(E,F)$\ \mbox{ is the {\rm ext}-pair}
&\qquad\Longleftrightarrow&\qquad \mu_H(E)>\mu_H(F);\\
\mbox{c)}& $(E,F)$\ \mbox{ is singular or the zero-pair}
 &\qquad\Longleftrightarrow&
\qquad \mu_H(E)=\mu_H(F).
\end{array}$$
\ETH

\PR Consider the restriction sequence to a smooth elliptic curve:
$$0\longrightarrow E^*\otimes F(K_S)\longrightarrow E^*\otimes F
    \longrightarrow (E^*\otimes F)|_D\longrightarrow 0.$$
   Denote $E|_D$ and $F|_D$ by $E'$ and $F'$. Combining Serre duality and
   the definition of exceptional pairs, we get
 ${\rm Ext}^i(E,F(K_S))^*\cong
   {\rm Ext}^{2-i}(F,E)=0.$  Hence the cohomology sequence associated
   with this exact sequence has the form:
   $$\begin{array}{|c|ccccc|}
   \hline
   k&{\rm Ext}^k(E,F(K_S))&\rightarrow
    &{\rm Ext}^k(E,F)&\rightarrow
    &{\rm Ext}^k(E',F')\\
   \hline
   0& 0 && * && * \\
   1& 0 && * && * \\
   2& 0 && * && *\\
   \hline
   \end{array}\; . $$
    That is,
   $${\rm Ext}^i(E,F)\cong{\rm Ext}^i(E',F')\qquad\forall i.$$
   Since $E'$ and $F'$ are  stable bundles on the elliptic curve
   (see the proof of lemma \ref{l228}), we obtain that  only one of the
   spaces ${\rm Ext}^0(E',F')$   and
  ${\rm Ext}^1(E',F')$  is nonzero whenever
   $$\mu_H(E')\not= \mu_H(F')\qquad\qquad(\mu(E')=\frac{\deg(E')}{r(E')}=
   \mu_H(E)).$$
  Moreover, ${\rm Ext}^0(E',F')\not= 0$  iff $\chi(E',F')>0$ and
           ${\rm Ext}^1(E',F')\not= 0$  iff $\chi(E',F')<0$.
   In this case $\chi(E',F')$  is the Euler characteristic of two sheaves
   on an elliptic curve  (\cite{K}):
 $$\chi(E',F')=r(E')r(F')\Bigl(\mu(F')-\mu(E')\Bigr).$$
 Finally, in any case we have ${\rm Ext}^2(E',F')=0$. This completes the
 proof.

\TH{Lemma.}\label{l232} Let $(E,F)$ be an exceptional pair of bundles
on $S$ with $\mu_H(E)=\mu_H(F)$. Let $C$ be
   $c_1(F)-c_1(E)$. Then:
\begin{enumerate}
\item $r(E)=r(F)$.
\item    $C^2=-2$ and $K_S\cdot C=0$.
\item  Suppose $(E,F)$ is a singular pair; then
\begin{enumerate}
\item  $C$ is a connected curve;
\item  ${\rm Ext}^0(E,F)={\rm Ext}^1(E,F)=\CC$;
\item  there exists an exact sequence
$$0\longrightarrow E\longrightarrow F\longrightarrow Q\longrightarrow 0,$$
     where $Q$ is a torsion sheaf with $c_1(Q)=C$.
\end{enumerate}
\end{enumerate}
\ETH

\PR By the definition of an exceptional pair, $\chi(F,E)=0$.
    Substituting the discrete invariants of  $E$ and $F$ in the
    Riemann-Roch formula for exceptional sheaves   (\ref{l212}), we get
$$0=\frac{1}{r^2(E)}+\frac{1}{r^2(F)}+\Bigl(\frac{c_1(F)}{r(F)}-
                                   \frac{c_1(E)}{r(E)}\Bigr)^2.$$
    From lemma \ref{l228} it follows that the restriction of an exceptional
    bundle to the elliptic curve $D\!\in\!|H|$ is a simple bundle. Moreover,
     $$\mu_H(E)=\mu(E|_D).$$

   If $L$ is a simple bundle on an elliptic curve; then $r(L)$ and $\deg (L)$
   are coprime (\cite{A}). Hence it follows from the equality
    $\mu(E|_D)=\mu(F|_D)$ that $r(E)=r(F)=r$.

     Whence we obtain
$$0=\frac{2}{r^2}+\frac{1}{r^2}\Bigl(c_1(F)-
                                   c_1(E)\Bigr)^2.$$
    This means that $$C^2=(c_1(F)-c_1(E))^2=-2.$$

   On the other hand, $$ \mu_H(E)=\frac{c_1(E)\cdot H}{r(E)}=\mu_H(F)=
                              \frac{c_1(F)\cdot H}{r(F)}.$$
  Therefore, $C\cdot H=0$, i.e., $C\cdot K_S=0$.
  This concludes the proof of the first and the second lemma statements.

3. Let $(E,F)$ be a singular pair, i.e., there exists a nonzero map
   $\varphi :E\longrightarrow F$. Since exceptional bundles on $S$
   are  $\mu_H$-stable, it follows from lemma \ref{l116} that
    $\varphi $ is an injection. Moreover, the cokernal of
    $\varphi$ has the zero rank. By definition, put
   $Q={\rm coker}\varphi$.
   Since the first Chern class is an additive function, we get
   $c_1(Q)=c_1(F)-c_1(E)=C$.

   Consider the restriction sequence:
   $$0\longrightarrow E^*\otimes F(K_S)\longrightarrow E^*\otimes F
   \longrightarrow (E^*\otimes F)|_D\longrightarrow 0.$$
   We have the following isomorphisms:
   $${\rm Hom}(E,F)\cong {\rm Hom}(E',F');\quad
   {\rm Ext}^1(E,F)\cong {\rm Ext}^1(E',F');\quad
   {\rm Ext}^2(E,F)=0,$$
  where $E'=E|_D$ and $F'=F|_D$.

   By assumption,  ${\rm Hom}(E,F)\not= 0$. Hence there exists a nonzero map
  $\varphi ':E'\longrightarrow F'$. Since $E'$ and $F'$ are stable
  bundles on a curve and $\mu(E')=\mu(F')$, we see that $\varphi'$ is
  an isomorphism. Further, stable bundles are simple (\ref{l117}), and the
   canonical class of an elliptic curve is trivial. It  follows from Serre
   duality that   ${\rm Ext}^1(E',F')=\CC$. Thus we have
   ${\rm Ext}^1(E,F)\cong{\rm Ext}^0(E,F)=\CC$.

  Now we show that $Q$ is simple.
  Let us write cohomology tables associated with the exact sequence
  $$0\longrightarrow E\longrightarrow F\longrightarrow Q\longrightarrow 0.$$

  $$\begin{array}{|c|ccccc|}
  \hline
  k&{\rm Ext}^k(E,E)&\rightarrow
   &{\rm Ext}^k(E,F)&\rightarrow
   &{\rm Ext}^k(E,Q)\\
  \hline
  0& \CC && \CC && ? \\
  1& 0   && \CC && ? \\
  2&  0  &&  0  && ? \\
  \hline
  \end{array}\; ;$$
 $$\begin{array}{|c|ccccc|}
  \hline
  k&{\rm Ext}^k(F,E)&\rightarrow
   &{\rm Ext}^k(F,F)&\rightarrow
   &{\rm Ext}^k(F,Q)\\
  \hline
  0& 0 && \CC && ? \\
  1& 0 &&  0  && ? \\
  2& 0 &&  0  && ?\\
  \hline
  \end{array}\; ; $$
  $$\begin{array}{|c|ccccc|}
  \hline
  k&{\rm Ext}^k(Q,Q)&\rightarrow
   &{\rm Ext}^k(F,Q)&\rightarrow
   &{\rm Ext}^k(E,Q)\\
  \hline
  0& ? && \CC && 0 \\
  1& ? &&  0  && \CC \\
  2& ? &&  0  && 0 \\
  \hline
  \end{array}\; .$$
  From the last table it follows that the quotient $Q$ is simple. Hence
  $C={\rm supp}Q$ is connected. In reality, the endomorphism group
  of $Q$ contains
   projectors unless ${\rm supp}Q$  is connected.

\vspace{2ex}

   The following lemma completes the classification of exceptional
   pairs of bundles on $S$.

\TH{Lemma.}\label{l233}  Suppose $(E,F)$ is an exceptional pair of
  bundles on $S$ with
  $\mu_H(E)=\mu_H(F)$ and  $C=c_1(F)-c_1(E)$ is a connected -2-curve
  ($C^2=-2$); then $(E,F)$ is singular.
\ETH

\vspace{1ex}

  To prove this statement, we need the following remark.

\TH{Remark.}\label{l234} Let $C$ be an irreducible smooth curve on the
 surface $S$ with
 $C^2=-2$ and $C\cdot H=0$. Then $C$ is rational and there are $x,y\in
 {\BBB N}$ and $d\in\ZZ$ such that
 $E'=E|_C=x{\cal O}_C(d)\oplus y{\cal O}_C(d+1)$
\ETH

\PR Let us show that $${\rm Ext}^2(E,E(-C))=0.$$
   By Serre duality, ${\rm Ext}^2(E,E(-C))^*\cong {\rm Hom}(E,E(K_S+C))$.

  Suppose, $K^2_S>0$; then $\mu_H(E)>\mu_H(E(K_S+C))$.
  Since exceptional bundles on $S$ are $\mu_H$-stable, it follows from
 \ref{l115} that ${\rm Hom}(E,E(K_S+C))=0$.

  Suppose, $K^2_S=0$; then $\mu_H(E)=\mu_H(E(K_S+C))$.
  This implies that any nonzero map
 $\varphi:E\longrightarrow E(K_S+C)$ is  injective. Denote by $Q$ the
 cokernal of $\varphi$. It is clear that the support of $Q$ is a curve and
 $c_1(Q)=r\cdot (K_S+C),$ where $r=r(E).$ Since the linear system $|C|$
 is one-dimensional and $|-K_S|$ has not base components, the divisor
 $r\cdot (K_S+C)$ cannot be effective.

 Thus, ${\rm Ext}^2(E,E(-C))=0$. Consider the restriction sequence to $C$:
$$0\longrightarrow E^*\otimes E(-C)\longrightarrow E^*\otimes E
    \longrightarrow (E^*\otimes E)|_C\longrightarrow 0.$$
 We have,
  $${\rm Ext}^1(E,E)\longrightarrow {\rm Ext}^1(E',E')\longrightarrow 0.$$
  By the definition of exceptional sheaves, ${\rm Ext}^1(E,E)=0$.
  This yields that $E'$ is a rigid bundle on the curve $C$. Finally, by
  adjunction,
 $$2g-2=C\cdot (C+K_S).$$ Therefore the curve $C$ is rational and $E'$ is a
 direct sum of line bundles. Now the remark follows from Bott formula
 (\cite{BOTT}).

\vspace{3ex}

{\sc Proof of lemma \ref{l233}}
  It is known that any irreducible component of a -2-curve $C$ on $S$ is a
  smooth rational -2-curve provided $C\cdot K_S=0$. It follows from
  the lemma condition, the definition of a singular pair and \ref{l231}
   that it is sufficient to show
  that there exists a nonzero map $E\longrightarrow F$. Hence, without loss
  of generality we can  assume that $C$ is an irreducible smooth rational
  curve.

 Let us denote $r=r(F)$, $D=c_1(F)$ and
\begin{equation}\label{delta}
\Delta G=c_1^2(G)-2c_2(G),
\end{equation}
  where $G$ is any coherent sheaf on $S$.

 Summing \ref{l232} and the lemma conditions, we get
 $$r(E)=r,\qquad\qquad c_1(E)=D-C.$$
 Since the sheaves  $E$ and $F$ are exceptional, it  follows from
 \ref{l212} that
 $$\Delta F=\frac{D^2+1}{r}-r^2\qquad\mbox{and}\qquad
   \Delta E=\frac{(D-C)^2+1}{r}-r^2, $$
   that is
\begin{equation}\label{eq1}
\Delta E=\Delta F-\frac{2}{r}(D\cdot C+1)
\end{equation}
  Substituting (\ref{delta}) in this equality, we obtain
$$\frac{D\cdot C+1}{r}=c_2(E)-c_2(F)+D\cdot (D-C)+1.
$$
 Since the right part of the last equality is integer, we get

\begin{equation}\label{eq2}
D\cdot C\equiv -1(mod\; r).
\end{equation}

 Remark \ref{l234} implies that
  $F'=F|_C=x{\cal O}_C(d)\oplus y{\cal O}_C(d+1)$,
 where
 $x+y=r$ and $xd+y(d+1)=D\cdot C$. Hence by (\ref{eq2}),
 $$F'={\cal O}_C(d)\oplus (r-1){\cal O}_C(d+1).$$
 Besides,

\begin{equation}\label{eq3}
d+1=\frac{D\cdot C+1}{r}.
\end{equation}

  We see that there exists an exact sequence:

\begin{equation}\label{eq5}
0\longrightarrow G\longrightarrow F\longrightarrow {\cal O}_C(d)
\longrightarrow 0.
\end{equation}
Using this triple, let us calculate the discrete invariants of $G$.
 It is easily shown that
\begin{equation}\label{eq4}
r(G)=r(F)=r(E),\qquad c_1(G)=c_1(E),\qquad \Delta G=\Delta F-
\Delta {\cal O}_C(d).
\end{equation}
  For calculating $\Delta {\cal O}_C(d)$ let us use Riemann-Roch formula
  (\ref{l211}):
$$\chi({\cal O},{\cal O}_C(d))=r({\cal O}_C(d))+
(H\cdot c_1({\cal O}_C(d)))/2+ (\Delta {\cal O}_C(d))/2.$$

Since $r({\cal O}_C(d))=0$ and $c_1({\cal O}_C(d))\cdot H=C\cdot H=0$,
 we get
$$\chi({\cal O},{\cal O}_C(d))=\frac{1}{2} \Delta {\cal O}_C(d).$$
 On the other hand,
$$\chi({\cal O},{\cal O}_C(d))=\chi({\cal O}_{{\BBB P}^1},
{\cal O}_{{\BBB P}^1}(d))=d+1,$$
  whereby
$$\Delta {\cal O}_C(d)=2(d+1).$$

 Combining  (\ref{eq3}), (\ref{eq4}) and the last equality, we obtain
$$\Delta G=\Delta F- \frac{2}{r}(D\cdot C+1).$$
 Therefore, $\Delta G=\Delta E$ (see (\ref{eq1})).

 We see that the sheaf $G$ has the same  discrete invariants as
 the  exceptional bundle $E$. Hence
$\chi(E,G)=\chi(E,E)=1$, i.e., $$h^0(E,G)+h^2(E,G)\ge 1.$$

  Suppose, $h^2(E,G)\not= 0$; then by Serre duality, $h^0(G,E(K_S))\not= 0$.
  That is there exists a nonzero map:
\begin{equation}\label{eq6}
G\longrightarrow E(K_S).
\end{equation}
  Since the bundle $F$ is $\mu_H$-stable, it  follows from  (\ref{eq5})
 that  $G$ is also  $\mu_H$-stable and
$\mu_H(G)=\mu_H(F)=\mu_H(E).$

 Summing lemma \ref{l116} and (\ref{eq6}) , we get $G\subset E(K_S)$ .
 Besides, the first Chern class of torsion  sheaf
 $E(K_S)/G$ is equal to $rK_S$.  This contradicts the effectivity of
 $c_1(E(K_S)/G)$. Thus,
 $h^0(G,E(K_S))=h^2(E,G)=0$. Consequently the equality $\chi(E,G)=1$ yields
 that there exists a nonzero morphism $E\longrightarrow G$. Since
 $E$ and $G$ have not torsion and $G$ is a subsheaf of $F$, we obtain
  ${\rm Hom}(E,F)\not= 0$. That is, the exceptional pair
 $(E,F)$ is singular.

\vspace{2cm}

%

%% FOLLOWING LINE CANNOT BE BROKEN BEFORE 80 CHAR
%%%%%%%%%%%%%%%%%%%%%%%%%%%%%%%%%%%%%%%%%%%%%%%%%%%%%%%%%%%%%%%%%%%%%%%%%%%%%%%%%%%%%
\subsection{Structure of Rigid Sheaves.}

   In the paper \cite{OK} it was proved that any rigid bundle on Del Pezzo
surface is a direct sum of exceptional bundles.
  But there exist indecomposable
rigid bundles $E$ with ${\rm Hom}(E,E)\not= \CC$ provided  $H=-K_S$ is
nef. For example, consider a -2-curve $C$ on $S$ with $C\cdot K_S=0$.
It is easily shown that $({\cal O_S},{\cal O_S}(C))$  is an exceptional
singular pair. Denote by $E$ a nontrivial extension of
${\cal O_S}$ by ${\cal O_S}(C)$:
$$0\longrightarrow {\cal O}(C)\longrightarrow E\longrightarrow
{\cal O}\longrightarrow 0.$$
  It can be proved that $E$ is rigid and ${\rm Hom}(E,E)\cong \CC ^2$.

  In this section we prove that any rigid bundle on the surface $S$ with
  nef  anticanonical class  and  $K_S^2>0$ has similar structure.
  Unfortunately, a structure of rigid bundles on $S$ with
  $K_S^2=0$ is not known. Further, we shall assume that $K_S^2>0$.

{\sc Definition.} We say that a torsion free sheaf  $F$ has an {\it
  exceptional  filtration\/} whenever  there exists a filtration of $F$
 $$Gr(F)=(x_nE_n,x_{n-1}E_{n-1},\ldots,x_1E_1),$$
 where $\VEC E,n$ is an exceptional collection of bundles such that
 $\mu_H(E_i)\leq \mu_H(E_{i+1})$ for $i=1,2,...,n-1$.

 \vspace{3ex}
 The aim of this section is to prove the following theorem:

\TH{Theorem.}\label{l241} Let $S$ be a smooth projective surface over
  $\CC$  with the anticanonical class without base components and
  $K_S^2>0$. Then
\begin{enumerate}
\item Any torsion free rigid sheaf on $S$ is a direct sum of
   $\mu_H$-semistable rigid bundles.
\item  Any indecomposable rigid sheaf without torsion on $S$ is
   $\mu_H$-semistable.
\item  Any $\mu_H$-semistable rigid sheaf has an exceptional
   filtration. Moreover, all pairs of associated exceptional collection have
   zero or singular type.
\end{enumerate}
\ETH

\vspace{2ex}

  We shall use the vector slope
$$\bar\gamma (E)=(\mu_H(E),\mu_A(E),\frac{c_1^2(E)-2c_2(E)}{r(E)}),$$
  where $A$ is an ample divisor, $H$ is the anticanonical class of $S$, and
 $\mu_D(E)=\frac{c_1(E)\cdot D}{r(E)}\quad$ with $D=\quad H$ or $A$.

  It can easily be checked that the stability with respect to this slope is
  the Gieseker stability \ref{l133}. In particular, any
  $\bar\gamma $-semistable sheaf has the filtration with isotypic quotients
   (\ref{l137}).

  The slope $\mu_H(E)$ has Mumford-Takemoto type and satisfy the axiom
  {\sc SLOPE.4}.

\TH{Lemma.}\label{l242} Let  $F$ be $\bar \gamma $-semistable
  rigid sheaf on $S$ with  $K_S^2>0$. Suppose $F$ has a filtration with
 $\bar \gamma $-stable isomorphic one to another quotients:
  $$Gr(F)=(G_n,G_{n-1},\ldots,G_1)\qquad (G_i\cong E\quad\forall i);$$
  Then $F$ is multiple of the exceptional bundle $E$, i.e $F=nE$.
\ETH

\PR Consider the spectral sequence associated with the filtration of $F$,
  which converge to the groups ${\rm Ext}^k(F,F)$. Its $E_1$-term
  has the form:
  $$E_1^{pq}=\bigoplus\limits_{i}^{}{\rm Ext}^{p+q}(G_i,G_{p+i}).$$

  Since the quotients $G_i\cong E$ are $\bar \gamma $-stable, we see
  that they are $\mu_H$-semistable (see remark \ref{l111}).
  Hence from lemma \ref{l135} it follows that the sheaf  $E(K_S)$
  is also  $\mu_H$-semistable.
  On the other hand, the square of the canonical class of our surface
  is positive. Thus,
  $$\mu_H(E(K_S))=\mu_H(E)+K_S\cdot H<\mu_H(E).$$

  Now, using Serre duality and lemma \ref{l115} we have
 $${\rm Ext}^2(E,E)=0.$$

  Thus, $E_1$-term of the spectral sequence has the form:

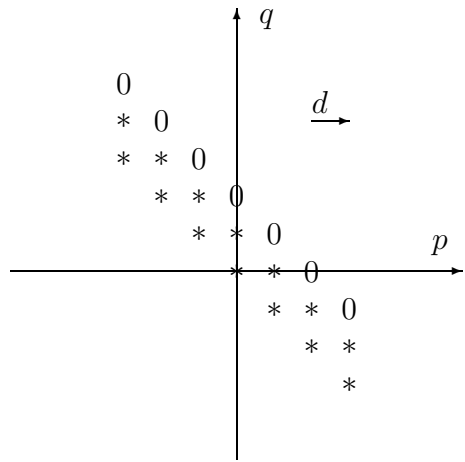
\begin{center}
\begin{picture}(70,80)(-35,-25)
%%%%%%%%%%%%%%%%%%%%%%%%%%%%%%%%%%%%%%%%%%%%%%%%%
% оформление картинки первого члена спектральной последовательности
%%%%%%%%%%%%%%%%%%%%%%%%%%%%%%%%%%%%%%%%%%%%%%%%%%%%%%%%%%%%%%%%%
  \put(3,33){$q$}
  \put(26,3){$p$}
 \put(10,20){\vector(1,0){5}}
  \put(10,21){$d$}
 \put(-30,0){\vector(1,0){60}}
 \put(0,-25){\vector(0,1){60}}
%%%%%%%%%%%%%%%%%%%%%%%%%%%%%%%%%%%%%%%%%%%%%%%%%%%%%%%%%%%%%%%%%%
%
% рисование собственно первого члена, начиная с левого верхнего элемента
% и по столбцам
%%%%%%%%%%%%%%%%%%%%%%%%%%%%%%%%%%%%%%%%%%%%%%%%%%%%%%%%%%%%%%%%%%

  \put(-15,25){\makebox(0,0){$0$}}
  \put(-15,20){\makebox(0,0){$*$}}
  \put(-15,15){\makebox(0,0){$*$}}

\put(-10,20){\makebox(0,0){$0$}}
\put(-10,15){\makebox(0,0){$*$}}
\put(-10,10){\makebox(0,0){$*$}}

  \put(-5,15){\makebox(0,0){$0$}}
  \put(-5,10){\makebox(0,0){$*$}}
  \put(-5,5){\makebox(0,0){$*$}}

  \put(0,10){\makebox(0,0){$0$}}
  \put(0,5){\makebox(0,0){$*$}}
  \put(0,0){\makebox(0,0){$*$}}

  \put(5,5){\makebox(0,0){$0$}}
  \put(5,0){\makebox(0,0){$*$}}
  \put(5,-5){\makebox(0,0){$*$}}

  \put(10,0){\makebox(0,0){$0$}}
  \put(10,-5){\makebox(0,0){$*$}}
  \put(10,-10){\makebox(0,0){$*$}}

  \put(15,-5){\makebox(0,0){$0$}}
  \put(15,-10){\makebox(0,0){$*$}}
  \put(15,-15){\makebox(0,0){$*$}}

            \end{picture}
\end{center}

  This yields that ${\rm Ext}^1(G_n,G_1)=E_1^{1-n,n}=
  E_\infty^{1-n,n}$ .
  On the other hand, $$E_\infty^{1-n,n}\subset {\rm Ext}^1(F,F)=0.$$
   But, $G_i\cong E\quad \forall i$. Consequently $E$ is a torsion free
   rigid sheaf. From lemma \ref{l221} it follows that $E$ is locally free.
  Besides, since $E$ is $\bar\gamma$-stable, we see that it is simple
  and ${\rm Ext}^2(E,E)=0$, whereby $E$ is an exceptional bundle.

  Finally, since  $G_i\cong E\quad \forall i$  and $E$ is exceptional, we
  have
 $${\rm Ext}^1(G_i,G_j)=0\quad \forall i,j.$$
 This implies the equality:
  $$F=\bigoplus\limits_{i}^{}G_i=nE.$$

\TH{Lemma.}\label{l243} Suppose $F$ is a rigid $\bar \gamma $-semistable
  sheaf on $S$; then $F$ is a direct sum of exceptional bundles.
\ETH

\PR From proposition \ref{l137} it follows that $F$ has a filtration
  with isotypic quotients:
  $$0=F_{n+1}\subset F_n\subset \cdots \subset F_2\subset F_1=F,$$
  where $G_i=F_i/F_{i+1}$ are $\bar \gamma $-semistable and they have
  filtrations with isomorphic one to another $\bar\gamma $-stable quotients.
  Besides,
  $$\forall i:\qquad \bar \gamma (G_i)=\bar \gamma (F),\qquad
  {\rm Hom}(F_{i+1},G_i)=0.$$

  Let us apply Mukai lemma (\ref{l214}) to the exact sequence:
  $$0\longrightarrow F_{i+1}\longrightarrow F_i\longrightarrow G_i
  \longrightarrow 0.$$
  It can be proved by induction on $i$ that $G_i$ and $F_{i+1}$ are
  torsion free rigid sheaves.

  Note that each $G_i$ satisfies assumption of the previous lemma.
  Therefore we have $G_i=x_iE_i$, where $E_i$ are exceptional bundles.

  Since all $G_i$ are $\bar \gamma $-semistable, we see that they are
  $\mu_H$-semistable (\ref{l111}). Moreover,
  $$\forall i:\qquad \bar \gamma (G_i)=\bar \gamma (F),\qquad
    \Longrightarrow \qquad\mu_H(G_i)=\mu_H(F).$$
  Whence, by the same argument as before, we have
  ${\rm Ext}^2(G_i,G_j)=0\quad\forall i,j$.
  Thus the $E_1$-term of the spectral sequence associated with the
  filtration of $F$ has the same form as in the proof of lemma \ref{l242}
  Therefore,  $${\rm Ext}^1(G_n,G_1)=0.$$
  But in this case, the filtration quotients of $F$ are different.
  To complete the proof we need an information about the
  groups  ${\rm Ext}^1(G_i,G_j)$ for $i<j$.

  Let us remember that  $G_i=x_iE_i$, where $E_i$ are the exceptional
  bundles. Hence,
\begin{equation}\label{eq761}
{\rm Ext}^1(G_n,G_1)\quad\Longrightarrow \quad {\rm Ext}^1(E_n,E_1)=0.
\end{equation}
   By  construction, the bundles  $E_i$ are
   $\bar \gamma $-stable and
   $\bar \gamma (E_i)=\bar \gamma (F)\quad\forall i$.
  Since $E_i$ are $\bar\gamma$-stable, it follows from lemma \ref{l134}
  that $E_n\cong E_1$ provided ${\rm Hom}(E_n,E_1)\not= 0$ or
  ${\rm Hom}(E_1,E_n)\not= 0$.

  Suppose $E_n\cong E_1$; then from (\ref{eq761}) we have
   ${\rm Ext}^1(E_1,E_n)=0$.

  Assume that  $E_n\not\cong E_1$; then we obtain
\begin{equation}\label{eq762}
{\rm Ext}^0(E_n,E_1)={\rm Ext}^0(E_1,E_n)=0.
\end{equation}
   Since $\bar \gamma (E_1)=\bar\gamma(E_n)$, we get
   $\mu_H(E_1)=\mu_H(E_n)$.

   Combining  $\mu_H$-stability of exceptional bundles on $S$,
   Serre duality and the inequality $K_S^2>0$, we obtain
   $${\rm Ext}^2(E_n,E_1)={\rm Ext}^2(E_1,E_n)=0.$$
   Combining this with (\ref{eq762}), we get
   $$\chi(E_1,E_n)=-h^1(E_1,E_n);\qquad\chi(E_n,E_1)=-h^1(E_n,E_1)=0.$$

   On the other hand, it  follows from Riemann-Roch theorem for
   exceptional sheaves
   (\ref{l212}) and the equality $\mu_H(E_1)=\mu_H(E_n)$ that
   $\chi(E_1,E_n)=\chi(E_n,E_1)$. That is,
   $$h^1(E_1,E_n)=h^1(E_n,E_1)=0.$$

   Thus we proved that ${\rm Ext}^1(E_1,E_n)=0$. This yields that
   ${\rm Ext}^1(G_1,G_n)=0$.

   Let us remark that the second term of filtration (i.e. $F_2$)
   satisfies the assumptions of our lemma. But
   $$Gr(F_2)=(G_n,G_{n-1},\ldots,G_2).$$
   By the inductive hypothesis, it can be assumed that
   $$F_2=\bigoplus\limits_{i=2}^{n}G_i.$$
   That is the sheaf $F$ is included in the exact sequence:
   $$0\longrightarrow \bigoplus\limits_{i=2}^{n}G_i
      \longrightarrow F\longrightarrow G_1\longrightarrow 0.$$
   Since ${\rm Ext}^1(G_1,G_n)=0$, we obtain  $F=\tilde F\oplus G_n$,
   where $\tilde F$ is $\bar \gamma $-semistable rigid sheaf with
   $$Gr(\tilde F)=(G_{n-1},G_{n-2},\ldots,G_1).$$

   Using the inductive hypothesis again, we have
   $\tilde F=\bigoplus\limits_{i=1}^{n-1}G_i$. That is,
   $$F=\bigoplus\limits_{i=1}^{n}G_i=
   \bigoplus\limits_{i=1}^{n}x_iE_i.$$
   This completes the proof.

\TH{Lemma.}\label{l244} Any rigid $\mu_H$-semistable sheaf  $F$
   on the surface $S$ has an exceptional filtration:
   $$Gr(F)=(x_mE_m,x_{m-1}E_{m-1},\ldots,x_1E_1)$$
   such that all exceptional pairs of collection
   $(E_1,E_2,\ldots,E_n)$ are singular or zero-pairs.
\ETH

\PR  Suppose $F$ is $\bar \gamma $-semistable; then by the previous lemma,
 $F=\bigoplus\limits_{i}^{}x_iE_i$, where
$E_i$ are $\bar \gamma $-stable exceptional bundles with equal one to
 another $\bar \gamma $-slopes. Without loss of generality it can be assumed
 that $E_i\not\cong E_j$ for
 $i\not= j$.  Using lemma \ref{l134}, we have
\begin{equation}\label{eq771}
{\rm Ext}^0(E_i,E_j)=0 \quad\forall i,j.
\end{equation}

 On the other hand the equality $\bar \gamma (E_i)=\bar \gamma (E_j)$ yields
 that  $\mu_H(E_i)=\mu_H(E_j)$.
  Now the equality
\begin{equation}\label{eq772}
{\rm Ext}^2(E_i,E_j)=0 \quad\forall i,j
\end{equation}
can be proved by the standard method on the surface $S$ with $K_S^2>0$.

  Finally, since $F$ is rigid, we have
$$0={\rm Ext}^1(F,F)={\rm Ext}^1(\bigoplus\limits_{i}^{}x_iE_i,
\bigoplus\limits_{i}^{}x_jE_j)=\bigoplus\limits_{i,j}^{}{\rm Ext}^1
(E_i,E_j),$$
  i.e., ${\rm Ext}^1(E_i,E_j)=0$.

 Combining this with equalities (\ref{eq771}) and (\ref{eq772}), we see that
 each pair of bundles in the collection $\VEC E,n$
 is an exceptional zero-pair.

  Now we suppose that $F$ is not $\bar\gamma$-semistable. Consider its
  Harder-Narasimhan filtration:
$$Gr(F)=(G_n,G_{n-1},\ldots,G_1)$$
   (see proposition \ref{l123}).

   Since $G_i$ are $\bar \gamma $-semistable and
$\bar \gamma (G_i)>\bar \gamma (G_{i-1})$ for all $i$, we get
\begin{equation}\label{eq773}
{\rm Ext}^0(G_i,G_j)=0 \quad\forall i>j.
\end{equation}

 Note that lemma \ref{l118} and $\mu_H$-semistability of the sheaf $F$
 imply $\mu_H$-semistability of the quotients $G_i$ and the equality
 $\mu_H(G_i)=\mu_H(F)$.
 Therefore, as before,
\begin{equation}\label{eq774}
{\rm Ext}^2(G_i,G_j)=0 \quad\forall i,j.
\end{equation}

  Combining (\ref{eq773}) and (\ref{eq774}), we see that the
 $E_1$-term of the spectral sequence associated with Harder-Narasimhan
 filtration of $F$ has the form:

\begin{center}
\begin{picture}(70,80)(-35,-25)
%%%%%%%%%%%%%%%%%%%%%%%%%%%%%%%%%%%%%%%%%%%%%%%%%
% оформление картинки первого члена спектральной последовательности
%%%%%%%%%%%%%%%%%%%%%%%%%%%%%%%%%%%%%%%%%%%%%%%%%%%%%%%%%%%%%%%%%
  \put(3,33){$q$}
  \put(26,3){$p$}
 \put(10,20){\vector(1,0){5}}
  \put(10,21){$d$}
 \put(-30,0){\vector(1,0){60}}
 \put(0,-25){\vector(0,1){60}}
%%%%%%%%%%%%%%%%%%%%%%%%%%%%%%%%%%%%%%%%%%%%%%%%%%%%%%%%%%%%%%%%%%
%
% рисование собственно первого члена, начиная с левого верхнего элемента
% и по столбцам
%%%%%%%%%%%%%%%%%%%%%%%%%%%%%%%%%%%%%%%%%%%%%%%%%%%%%%%%%%%%%%%%%%

  \put(-15,25){\makebox(0,0){$0$}}
  \put(-15,20){\makebox(0,0){$*$}}
  \put(-15,15){\makebox(0,0){$0$}}

\put(-10,20){\makebox(0,0){$0$}}
\put(-10,15){\makebox(0,0){$*$}}
\put(-10,10){\makebox(0,0){$0$}}

  \put(-5,15){\makebox(0,0){$0$}}
  \put(-5,10){\makebox(0,0){$*$}}
  \put(-5,5){\makebox(0,0){$0$}}

  \put(0,10){\makebox(0,0){$0$}}
  \put(0,5){\makebox(0,0){$*$}}
  \put(0,0){\makebox(0,0){$*$}}

  \put(5,5){\makebox(0,0){$0$}}
  \put(5,0){\makebox(0,0){$*$}}
  \put(5,-5){\makebox(0,0){$*$}}

  \put(10,0){\makebox(0,0){$0$}}
  \put(10,-5){\makebox(0,0){$*$}}
  \put(10,-10){\makebox(0,0){$*$}}

  \put(15,-5){\makebox(0,0){$0$}}
  \put(15,-10){\makebox(0,0){$*$}}
  \put(15,-15){\makebox(0,0){$*$}}

            \end{picture}
\end{center}

  This spectral sequence is convergent to the groups
   ${\rm Ext}^k(F,F)$ of the rigid sheaf. Hence,
\begin{equation}\label{eq775}
{\rm Ext}^1(G_i,G_j)=0 \quad\forall i\geq j.
\end{equation}

 In particular,  $G_i$ are rigid $\bar \gamma $-semistable sheaves. By the
 previous lemma, $G_i=\bigoplus\limits_{k}^{}x_{ik}E_{ik}$, where
  $E_{ik}$ are exceptional bundles. Besides, any pair
 $(E_{ik},E_{is})$ has the zero type.

   Combining (\ref{eq773}), (\ref{eq774}) and (\ref{eq775}),
   we obtain that
${\rm Ext}^d(E_{ik},E_{js})=0$ for $i>j$ and $d=0,1,2$.
 In other words, the set of all bundles $E_{ik}$ can be numerated such that
the collection $\VEC E,m$ will be exceptional.
 It remains to remark that all bundles $E_i$ have the
 $\mu_H$-slope coincided with
 $\mu_H(F)$. Thus it follows from proposition \ref{l231} that each pair of
 this collection is either singular or a zero-pair.

\vspace{2ex}

The proof plan of the
  main theorem is clear. We consider the spectral sequence
  associated with Harder-Narasimhan filtration of a rigid torsion free
  sheaf for obtaining the information about the groups
  ${\rm Ext}^1(G_i,G_j)$, where $G_i$ are quotients of this filtration.
  For this we need the following last statement.

\TH{Lemma.}\label{l245} Let  $G_1$ and $G_2$ be  $\mu_H$-semistable
  rigid sheaves on the surface $S$. Suppose
   $\mu_H(G_2)>\mu_H(G_1)$; then the equality ${\rm Ext}^1(G_2,G_1)=0$
implies ${\rm Ext}^1(G_1,G_2)=0$.
\ETH

\PR It follows from lemma \ref{l244} that each of $G_i$ has the exceptional
    filtration:
$$Gr(G_1)=(x_{n1}E_{n1},\ldots,x_{11}E_{11});\quad
  Gr(G_2)=(x_{m2}E_{m2},\ldots,x_{12}E_{12}).$$
  Moreover, $\mu_H$-slopes of $E_{ij}$ do not depend on the first index,
  i.e,
$\mu_H(E_{i1})=\mu_H(G_1)$ and $\mu_H(E_{j2})=\mu_H(G_2)$.

  Denote by $G'_i$ the restriction of the sheaves $G_i$ to an elliptic
  curve $D\!\in\!|-K_S|$. It is obvious that the sheaves $G'_i$
  have the filtrations:
$$Gr(G'_1)=(x_{n1}E'_{n1},\ldots,x_{11}E'_{11}),\qquad
  Gr(G'_2)=(x_{m2}E'_{m2},\ldots,x_{12}E'_{12}),$$
  where $E'_{ki}=E_{ki}|_D$. Further, since $E_{ki}$ are exceptional
  bundles, we see that $E'_{ki}$ are stable with respect to the standard
  slope $\mu$ on a curve (see lemma \ref{l228}). Moreover,
 $$\mu (E'_{ki})=\mu_H(E_{ki})=\mu_H(G_i)=\mu(G'_i).$$
 Now it follows from lemma \ref{l118} that $G'_i$ are $\mu$-semistable
 and $\mu(G'_2)>\mu(G'_1)$. Thus from lemma \ref{l115} it follows that
  ${\rm Hom}(G'_2,G'_1)=0$.

  Using the last equality and the long cohomology sequence associated with
  the exact triple
$$0\longrightarrow G^*_2\otimes G_1(K_S)\longrightarrow G^*_2\otimes G_1
   \longrightarrow (G^*_2\otimes G_1)|_D\longrightarrow 0,$$
we obtain
$${\rm Ext}^1(G_2,G_1(K_s))\subset {\rm Ext}^1(G_2,G_1).$$
 Now the proof follows from Serre duality.

\vspace{2ex}

{\sc Proof of theorem} \ref{l241}.
  Let  $F$ be any torsion free rigid sheaf on $S$. Consider its
  Harder-Narasimhan filtration with $\mu_H$-semistable quotients
$$Gr(F)=(G_n,G_{n-1},\ldots,G_1).$$
  It follows from the inequalities $\mu_H(G_i)>\mu_H(G_j)$  for $i>j$ that
$${\rm Hom}(G_i,G_j)=0 \;\mbox{for}\; i>j\quad\;\mbox{and}\quad\;
  {\rm Ext}^2(G_j,G_i)=0\;\mbox{for}\; i\ge j.$$

  Therefore  $E_1$-term of the spectral sequence associated with this
  filtration has the form:
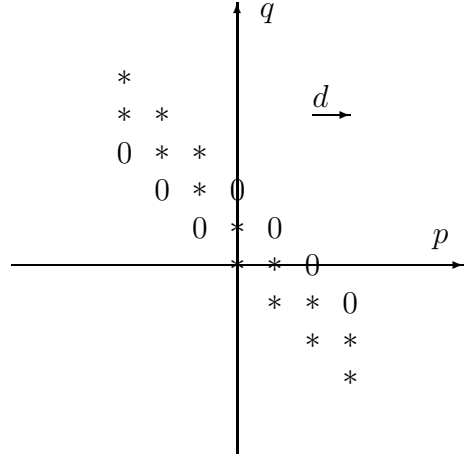
\begin{center}
\begin{picture}(70,80)(-35,-25)
%%%%%%%%%%%%%%%%%%%%%%%%%%%%%%%%%%%%%%%%%%%%%%%%%
% оформление картинки первого члена спектральной последовательности
%%%%%%%%%%%%%%%%%%%%%%%%%%%%%%%%%%%%%%%%%%%%%%%%%%%%%%%%%%%%%%%%%
  \put(3,33){$q$}
  \put(26,3){$p$}
 \put(10,20){\vector(1,0){5}}
  \put(10,21){$d$}
 \put(-30,0){\vector(1,0){60}}
 \put(0,-25){\vector(0,1){60}}
%%%%%%%%%%%%%%%%%%%%%%%%%%%%%%%%%%%%%%%%%%%%%%%%%%%%%%%%%%%%%%%%%%
%
% рисование собственно первого члена, начиная с левого верхнего элемента
% и по столбцам
%%%%%%%%%%%%%%%%%%%%%%%%%%%%%%%%%%%%%%%%%%%%%%%%%%%%%%%%%%%%%%%%%%

  \put(-15,25){\makebox(0,0){$*$}}
  \put(-15,20){\makebox(0,0){$*$}}
  \put(-15,15){\makebox(0,0){$0$}}

\put(-10,20){\makebox(0,0){$*$}}
\put(-10,15){\makebox(0,0){$*$}}
\put(-10,10){\makebox(0,0){$0$}}

  \put(-5,15){\makebox(0,0){$*$}}
  \put(-5,10){\makebox(0,0){$*$}}
  \put(-5,5){\makebox(0,0){$0$}}

  \put(0,10){\makebox(0,0){$0$}}
  \put(0,5){\makebox(0,0){$*$}}
  \put(0,0){\makebox(0,0){$*$}}

  \put(5,5){\makebox(0,0){$0$}}
  \put(5,0){\makebox(0,0){$*$}}
  \put(5,-5){\makebox(0,0){$*$}}

  \put(10,0){\makebox(0,0){$0$}}
  \put(10,-5){\makebox(0,0){$*$}}
  \put(10,-10){\makebox(0,0){$*$}}

  \put(15,-5){\makebox(0,0){$0$}}
  \put(15,-10){\makebox(0,0){$*$}}
  \put(15,-15){\makebox(0,0){$*$}}

            \end{picture}
\end{center}

   Since the sequence is convergent to the groups ${\rm Ext}^i(F,F)$ of the
   rigid sheaf, we obtain
$$0=E_\infty^{-1,2}=E_1^{-1,2}=
    \bigoplus\limits_{i}^{}{\rm Ext}^1(G_i,G_{i-1}),$$
$$0=E_\infty^{0,1}=E_1^{0,1}=
    \bigoplus\limits_{i}^{}{\rm Ext}^1(G_i,G_{i}),$$
   That is  $G_i$ are rigid  $\mu_H$-semistable sheaves and
   ${\rm Ext}^1(G_i,G_{i-1})=0$.

   By the previous lemma, the groups ${\rm Ext}^1(G_{i-1},G_{i})$
   are also trivial. In particular,
$${\rm Ext}^1(G_{1},G_{2})=0.$$

   Let $F_2$ be the first term of the filtration $Gr(F)$, i.e.,
\begin{equation}\label{filtr}
0\longrightarrow F_2\longrightarrow F\longrightarrow G_1\longrightarrow 0.
\end{equation}
   Note that $Gr(F_2)=(G_n,G_{n-1},\ldots,G_2)$ is also Harder-Narasimhan
   filtration and $\mu_H(G_2)>\mu_H(G_1)$. Taking into account corollary
   \ref{l125}, we obtain ${\rm Hom}(F_2,G_1)=0$. In addition, the sheaves
   $F_2$ and $G_1$ have not torsion. Hence we can apply to these sheaves
   lemma \ref{l222}. That is, ${\rm Ext}^2(G_1,F_2)=0$.

   Now we apply Mukai lemma to (\ref{filtr}) to  obtain
$$h^1(F,F)\ge h^1(F_2,F_2)+h^1(G_1,G_1).$$
   That is the sheaf $F_2$ is also rigid. The number of its Harder-Narasimhan
   filtration quotients is less then $n$. Hence by the inductive
   hypothesis we have $F_2=\bigoplus\limits_{i=2}^{n}G_i,$ and
$$0\longrightarrow \bigoplus\limits_{i=2}^{n}G_i
   \longrightarrow F\longrightarrow G_1\longrightarrow 0.$$
   Let us remember that ${\rm Ext}^1(G_1,G_2)=0.$ Therefore,
   $F=\tilde F\oplus G_2$, where $\tilde F$ is a torsion free rigid sheaf.
   We apply again the inductive hypothesis to the sheaf $\tilde F$ to obtain
$$F=\bigoplus\limits_{i=1}^{n}G_i.$$

   Thus any rigid sheaf without torsion on $S$ is a direct sum of
   $\mu_H$-semistable rigid ones. They are locally free by virtue of
   \ref{l221}. In particular, if $F$ is indecomposable then it is
   $\mu_H$-semistable. This concludes the proof of the first and the second
   theorem   statements. The last one is equivalent to  lemma \ref{l244}.

\vspace{2ex}

\TH{Corollary.}\label{l246} Any torsion free rigid sheaf on Del Pezzo surface
   $X$ is a direct sum of exceptional bundles
\ETH

\PR Since the anticanonical class of Del Pezzo surface is ample, we see that
    an exceptional pairs on $X$ cannot be singular (see \ref{l232}).
    On the other hand, we have proved that any indecomposable torsion free
    rigid sheaf on $S$ (in particular, on $X$) has an exceptional filtration.
    Besides, all pairs in associated exceptional collection are singular
    or zero pairs.  Thus the quotients of the exceptional filtration of
    any torsion free rigid sheaf on $X$ are its direct summands.
%

%% FOLLOWING LINE CANNOT BE BROKEN BEFORE 80 CHAR
%%%%%%%%%%%%%%%%%%%%%%%%%%%%%%%%%%%%%%%%%%%%%%%%%%%%%%%%%%%%%%%%%%%%%%%%%%%%%%%%%%%%%
\subsection{Structure of superrigid sheaves.}

  In the previous section we have proved that  any $\mu_H$-semistable
   sheaf has
  the exceptional filtration  and any torsion free rigid sheaf
  is a direct sum of $\mu_H$-semistable rigid bundles.
  Whereby for classifying rigid bundles we need a description of exceptional
  bundles and collections of ones. This description is the subject-matter
  of the next part. But for it we need the following theorem

 \TH{Theorem.}\label{l251} Let $S$ be a smooth projective surface over
    $\CC$ with the anticanonical class $H$ without base components
    and $H^2>0$. Then the following statements hold true.

 1. For any exceptional collection of bundles $\VEC E,n$ on $S$ such that
 $\forall i\quad\mu_H (E_i)\le \mu_H (E_{i+1})$ there exists a superrigid
 bundle $E$ $$({\rm Ext}^1(E,E)={\rm Ext}^2(E,E)=0)$$ such that
  $Gr(E)=(x_nE_n,x_{n-1}E_{n-1},\ldots,x_1E_1)$. We say that this bundle
  is {\it associated with the exceptional collection\/}.

 2. Any superrigid torsion free sheaf $E$ has the exceptional filtration
$$Gr(E)=(x_nE_n,x_{n-1}E_{n-1},\ldots,x_1E_1),
$$
   i.e. the collection $\VEC E,n$ is exceptional and the
   $\mu_H$-slopes of bundles $E_i$ satisfy the inequalities:
   $\mu_H(E_i)\leq \mu_H(E_{i+1})\quad\forall i$.

3. Suppose a superrigid torsion free sheaf $E$ has two exceptional
   filtrations:
$$Gr(E)=(x_nE_n,x_{n-1}E_{n-1},\ldots,x_1E_1)=
(y_mF_m,y_{m-1}F_{m-1},\ldots,y_1F_1);$$
   then $m=n$ and the exceptional collection $\VEC F,m$ can be obtained from
   $\VEC E,n$ by mutations of  neighboring zero-pairs
   $(E_i,E_{i+1})$.
\ETH

\vspace{2ex}

    Note that this theorem is obvious provided $S$ is Del Pezzo surface
    (see corollary \ref{l246}). But if $-K_S$ is nef then this theorem is
    interesting and its proof is difficult.

    Let us state and prove several lemmas.

\TH{Remark.}\label{l252} Suppose a sheaf $F$ has a filtration
    $Gr(F)=(G_n,G_{n-1},\ldots,G_1)$ such that
$Gr(G_i)=(E_{ik_i},\ldots,E_{i1})$; then there exists the filtration
 $$Gr(F)=(E_{nk_n},\ldots,E_{n1},\ldots,E_{1k_1},\ldots,E_{11}).$$
 And back to front, the neighboring quotients can be "join".
\ETH

\TH{Remark.}\label{l253} Suppose $Gr(F)=(G_n,G_{n-1},\ldots,G_1)$ is a
 filtration of a sheaf $F$ such that
  ${\rm Ext}^1(G_i,G_{i+1})={\rm Ext}^1(G_{i+1},G_i)=0$  for same $i$;
  then the  sheaf  $F$ has the filtration
   $Gr(F)=(G_n,G_{n-1},\ldots,G_i,G_{i+1},\ldots G_1)$.
\ETH

\TH{Lemma.}\label{l254} Suppose $F$ is an indecomposable rigid bundle on
  $S$ with $K_S^2>0$; then $F$ has the following filtrations:

  a) $Gr_R(F)=(Q_n,Q_{n-1},\ldots,Q_1)$ such that $\forall i\quad
  Q_i=\bigoplus y_{is}E_{is};$ $E_{is}$ are exceptional bundles,
  the collection $(E_{11},\ldots,E_{1m_1},\ldots,E_{n1},\ldots,E_{nm_n})$
  is exceptional and for each bundle   $E_{is}\quad (i=1,\ldots,n-1)$
  there is $E_{i+1,k}$ such that the pair $(E_{i,s},E_{i+1,k})$ is singular.

  b) $Gr_L(F)=(G_n,G_{n-1},\ldots,G_1)$ such that $\forall i\quad
  G_i=\bigoplus x_{is}E_{is};$ $E_{is}$ are exceptional bundles,
  the collection $(E_{11},\ldots,E_{1k_1},\ldots,E_{n1},\ldots,E_{nk_n})$
  is exceptional and for each bundle   $E_{is}\quad (i=2,\ldots,n)$
  there is $E_{i-1,l}$ such that the pair $(E_{i-1,l},E_{is})$ is singular.
\ETH

\PR Let us construct the first filtration. Similarly
  the second one is constructed.

  By the theorem about rigid bundles (\ref{l241}) the sheaf $F$ has an
  exceptional filtration
    $$Gr(F)=(x_nE_n,x_{n-1}E_{n-1},\ldots,x_1E_1).$$
  Let us partition  the exceptional collection
    $$\VEC E,n,$$
  associated with this filtration into subcollections
  $$(E_{i_{s-1}+1},E_{i_{s-1}+2},\ldots,E_{i_{s}}),$$ where $i_0=0$
   such that the  following conditions hold.

      {\sc PART.1:}
{\sl  Any pair of each subcollection has the zero type whenever this
  subcollection

  has greater than one bundle.}

      {\sc PART.2:}
{\sl  For last bundle  $E_{i_s}$ of each subcollection there is
      $E_j$ in the next

       subcollection such that the pair $(E_{i_s},E_j)$
      is singular.}

\noindent  There exists at most one singular pair in this collection because
   the bundle $F$ is indecomposable. This implies that this partition exists.

    Denote by  $Q_s$ the direct sum
    $$\bigoplus\limits_{j=i_{s-1}+1}^{i_s}x_jE_j$$
    of exceptional bundles from the  subcollection with index $s$.
    It follows from remark \ref{l252} that $F$ has filtration
    $Gr(F)=(Q_k,Q_{k-1},\ldots,Q_1),$ where $k$ is the quantity of
    all subcollections. Note that this filtration consists with
    $Gr_R$ if and only if the collection decomposed into subcollections
    satisfies the conditions {\sc PART.1 } and the following

    {\sc PART.2R:}
{\sl  for any bundle $E_{i}$ of each subcollection there is a bundle

      $E_j$ in the next subcollection such that the pair $(E_{i},E_j)$
      is singular.}

    For constructing the required collection we shall intermix the bundles
    of subcollections and move sometimes bundles from  subcollection to
    next one.

    Suppose there is a bundle $E_{\alpha}$ in the first subcollection
    such that for all bundles $E_{\beta}$ in the second one
      $(E_\alpha,E_\beta)$ are
     zero-pairs. Let us move $E_\alpha$ to the second subcollection.
    Since this moving can be realized by permutations of members in
    neighboring zero-pairs, we see that the obtained collection is
    exceptional. Besides, it satisfy the conditions
    {\sc PART.1} and {\sc PART.2.} It is clear that after a finite number of
    such movings we get the exceptional collection decomposed into
    subcollections such that for each bundle $E_\alpha$ of the first
    subcollection there is $E_\beta$ in the second one such that
    $(E_\alpha,E_\beta)$ is  singular pair. Let us mention, that one can
    do the some thing with an arbitrary pair of neighboring subcollections.
    Let this process be called the {\it displacement\/}.

    Let us do the displacement with each pair of the neighboring
    subcollections, beginning from the first one. The number of the
    subcollections cannot be changed during the process. The number of the
    bundles in the first subcollection can lessen only. Two latter
    subcollections will satisfy the condition {\sc PART.2.R}. But
    since we moved
    the bundles from the left to the right, now one can find two neighboring
    subcollections (with the numbers $s$ and $s+1$, for example) satisfying
    the following conditions. Any pair $(E_i,E_j)$ with $E_i$ belonging to
    the $s$-th subcollection and $E_j$ belonging to the $(s+1)$-th
    subcollection has the type zero. Moreover, one can warrant that
    two latter subcollection satisfy the condition {\sc PART.2R} only.

    Let us join, if it is necessary, the neighboring subcollections to
    satisfy the conditions {\sc PART.1} and {\sc PART.2}.

    Let us do the displacement with each pair of neighboring subcollection
    and join all what is possible to join, ets...

    This process cannot be repeated ad infinitum. Really, there exists
    $k_0\in\NN$ such that for any $k>k_0$ the number of the subcollections
    will not change after doing the $k$-th step --- "the displacement and
    the join". After some successive step the number of bundles in the
    first subcollection will not change, ets... Thus, the number of the
    subcollections and the number of the bundles in each subcollection
    will not change since some moment. That means , any bundle does not go
    from one subcollection to another. Hence we are done.

\vspace{2ex}

   I should like to remark that the last lemma is trivial provided each
   pair of the exceptional collection associated with a rigid bundle
   is singular. But it is not the case. On the surface $S$ there is
 an exceptional collection of bundles of equal $\mu_H$-slope consisted of
   both singular and zero-pairs. Besides, the superrigid bundle
   associated with this collection
    is indecomposable. The following example was found by Yu.~B.~
   Zuev.

   Let $X_1$ be a surface obtained from ${\BBB P}^2$ by blowing up a point
   $x_1:\quad \sigma_1:X_1\longrightarrow {\BBB P}^2$. Denote by $e$ the
   preimage of $x_1$ ($e=\sigma_1^{-1}(x_1)$). Let us choose two points
   $x_2$ and $x_3$ on $e$. Suppose $X$ is obtained from $X_1$ by blowing up
   $x_2$ and $x_3:\quad \sigma_2:X\longrightarrow X_1.$
   By definition, put
   $$e_1=\sigma_2^{-1}(e),\quad e_2=\sigma_2^{-1}(x_2),\quad
   e_3=\sigma_2^{-1}(x_3).$$
   We see that the curves $e_1,e_2$ and $e_3$ are exceptional, that is
   $e_i^2=-1$ for $i=1,2,3$. Besides, $(e_1-e_2)$ and $(e_1-e_3)$ are
   connected $-2$-curves. It can easily be checked that the collection of
   line bundles
   $$\Bigl({\cal O}_X, {\cal O}_X(e_1-e_2), {\cal O}_X(e_1-e_3)\Bigr)$$
   is exceptional. Moreover, the pairs $\Bigl({\cal O}_X,{\cal O}_X(e_1-e_2)
   \Bigr) $ and $\Bigl({\cal O}_X,{\cal O}_X(e_1-e_3)\Bigr)$
   are singular and
   $\Bigl({\cal O}_X(e_1-e_2),{\cal O}_X(e_1-e_3)\Bigr)$
   is a zero-pair.

\TH{Lemma.}\label{l255} Let $F$ be an indecomposable torsion free rigid
   sheaf on the surface $S$. Assume that
 $$ Gr(F)=(x_nE_n,x_{n-1}E_{n-1},\ldots,x_1E_1)=
          (y_mF_m,y_{m-1}F_{m-1},\ldots,y_1F_1)$$
  are two exceptional filtrations of $F$, i.e. the collections
 $$\VEC E,n \quad\mbox{and} \VEC F,m$$ are exceptional.  Suppose
 \begin{equation}\label{slope}
 \mu_H(E_i)=\mu_H(F_j)=\mu_H(F)\qquad \forall i,j;
 \end{equation}
 then $m=n$ and the collection $\VEC E,n$ can be obtained from $\VEC F,m$
 by mutations of neighboring zero-pairs.
\ETH

\PR  It follows from  proposition \ref{l231}
     that each pair of these collections has zero or singular type. Let us
     show that any such collection can be ordered by $\bar\gamma$-slope
      by virtue of  permutations of neighboring zero-pairs members only.
     In this case the lemma follows from the uniqueness of Harder-Narasimhan
     filtration (\ref{l123}).

     The possibility of such ordering is obtained by induction on the
     number of collection members from the following arguments.

     Suppose $(E,F)$ is a singular pair; then it follows from lemma \ref{l232}
     that ranks of the sheaves $E$ and $F$ are equal and $c_1(F)-c_1(E)=C$
     is an effective -2-divisor. (Recall that the $\bar \gamma $-slope is the
     vector
      $$(\mu_H,\mu_A,\frac{c_1^2-2c_2}{r}),$$
     where $\mu_A=\frac{c_1\cdot A}{ r}$, and $A$ is an ample divisor.)
     Since $A$ is ample, we get $\mu_A(E)<\mu_A(F)$. By assumption we have
     $\mu_H(E)=\mu_H(F)$. Therefore,
     $\bar \gamma (E)<\bar \gamma (F)$.

\TH{Lemma.}\label{l256} Suppose $(E,F)$ is an exceptional singular pair on
   $S$ and $G$ is a torsion free sheaf; then the following implications hold
   true.

a) ${\rm Ext}^2(G,E)=0\quad \Longrightarrow\quad  {\rm Ext}^2(G,F)=0$;

b) ${\rm Ext}^0(G,F)=0\quad \Longrightarrow\quad  {\rm Ext}^0(G,E)=0$;

c) ${\rm Ext}^0(E,G)=0\quad \Longrightarrow\quad  {\rm Ext}^0(F,G)=0$;

d) ${\rm Ext}^2(F,G)=0\quad \Longrightarrow\quad  {\rm Ext}^2(E,G)=0$;

\ETH

\PR   The lemma follows from the cohomology tables associated with the
   exact triple
 $$0\longrightarrow E\longrightarrow F\longrightarrow Q\longrightarrow 0,$$
   where $Q$ is a torsion sheaf. Besides, since $Q$ has zero rank and $G$ is
   torsion free, we see that  ${\rm Hom}(Q,G)=0$. Moreover, using Serre
   duality, we have ${\rm Ext}^2(G,Q)=0$.

\TH{Lemma.}\label{l257} Let $F$ be a
 $\mu_H$-semistable rigid bundle. Let
$$Gr(F)=(x_nE_n,x_{n-1}E_{n-1},\ldots,x_1E_1)$$ be its exceptional
 filtration and $G$ sheaf without torsion. Then

a) ${\rm Ext}^i(G,F)=0\; \forall i\quad \Longleftrightarrow \quad
     {\rm Ext}^i(G,E_k)=0\; \forall i,k$;

b) ${\rm Ext}^i(F,G)=0\; \forall i\quad \Longleftrightarrow \quad
     {\rm Ext}^i(E_k,G)=0\; \forall i,k$.

\ETH

\PR  Without loss of generality it can be assumed that $F$ is indecomposable.
     Consider its filtration
   $Gr_R(F)=(Q_m,Q_{m-1},\ldots,Q_1)$ from lemma \ref{l254}.
     By \ref{l255}, we can assume without loss of generality that
  $$Q_j=\bigoplus\limits_{i=s_j}^{s_{j+1}-1}y_iE_i\qquad
  1=s_1<s_2<\cdots<s_m<s_{m+1}=n+1.$$

  To prove the first statement of our lemma it is sufficient to check the
  following implication
  $${\rm Ext}^i(G,F)=0\; \forall i\quad \Longrightarrow \quad
     {\rm Ext}^i(G,E_k)=0\; \forall i,k.$$
(The another implication follows from \ref{l124}.)
   Let us apply the functor ${\rm Ext}^\cdot(G,\cdot)$ to each of the
sequences:
   $$0\longrightarrow F_{j+1}\longrightarrow F_j\longrightarrow
   Q_j\longrightarrow 0,$$
   where $F_j$ are terms of the filtration $Gr_R(F)$, $F_1=F$ and $F_n=Q_n$.
   We see that  for $j=2,3,\ldots,n$ ${\rm Hom}(G,F_j)=0$.
   In particular, ${\rm Hom}(G,Q_n)=0$.

{\sl Step 1.} Let us show that  ${\rm Hom}(G,E_i)=0$ for all $i$.

  It follows from the equality ${\rm Hom}(G,Q_n)=0$ that
  ${\rm Hom}(G,E_i)=0$ for $s_m\leq i\leq n$.
  By the construction of the filtration $Gr_R(F)$, for each direct
  summand  $E_\alpha$ of the bundle $Q_{n-1}$ there exists a bundle
  $E_\beta$ with $s_m\leq \beta\leq n$ such that
  $(E_\alpha,E_\beta)$ is a singular pair. The application of lemma
  \ref{l256}  to this pair yields, ${\rm Hom}(G,E_\alpha)=0$.

  In the same way, using the properties of the filtration $Gr_R(F)$ and
   \ref{l256}, we conclude the first step.

\vspace{2ex}

{\sl Step 2.} Let us check that ${\rm Ext}^2(G,E_i)=0$ for all $i$.

  Now let us intermix the bundles in collection $\VEC E,n$ to obtain
  the filtration $Gr_L(F)$. Without loss of generality we can assume that
  $$G_j=\bigoplus\limits_{i=s_j}^{s_{j+1}-1}y_iE_i\qquad
  1=s_1<s_2<\cdots<s_m<s_{m+1}=n+1.$$
  are quotients of $Gr_L(F)$.

  Applying the functor  ${\rm Ext}^\cdot(G,\cdot)$ to the exact triple
  $$0\longrightarrow F'_2\longrightarrow F\longrightarrow
  G_1\longrightarrow 0,$$
  where $F'_2$ is the first term of the filtration $Gr_L(F)$,
  we get, ${\rm Ext}^2(G,G_1)=0$. This means that ${\rm Ext}^2(G,E_i)=0$
  for any direct summand of the bundle $G_1$.

  Using the properties of the filtration $Gr_L(F)$ and lemma  \ref{l256},
  as  before, we have
      $${\rm Ext}^2(G,E_i)=0\quad \forall i.$$

  Thus it was proved that for any quotient $E_i$ of the exceptional filtration
    of the bundle $F$ the groups
  ${\rm Ext}^0(G,E_i)$ and ${\rm Ext}^2(G,E_i)$ are trivial. Hence,
  $\chi(G,E_i)\leq 0\quad\forall i$. Since Euler characteristic of sheaves is
  additive function, we have
  $$\sum\limits_{i=1}^{n}x_i\chi(G,E_i)=\chi(G,F)=0.$$   Moreover,
  all $x_i$ are natural numbers. Thus,
  $\chi(G,E_i)=0\;\forall i$.
  The first statement of the lemma was proved. Similarly the second one is
  proved.

\vspace{2ex}

{\sc Proof of theorem \ref{l251}}

1. At first assume that all pairs of the collection
$\VEC E,n$ have zero or singular type. The proof is by induction
  on the number $n$ of
  bundles in the collection. For $n=1$, there is nothing to prove.

  By induction hypothesis there exists a superrigid bundle $E'$ such that
 $$Gr(E')=(E_n,E_{n-1},\ldots,E_2).$$ Suppose the pair
  $(E_1,E_i)$ has the zero type for any $i$; then  $E=E'\oplus E_1$
  is a superrigid bundle (see \ref{l124}).

  Suppose there is an index $i$ such that $(E_1,E_i)$ is singular;
  then  ${\rm Ext}^k(E_1,E_i)=\CC$ for $k=0,1$ and
$$ {\rm Ext}^2(E_1,E_j)={\rm Ext}^k(E_j,E_1)=0\quad\forall j,k.$$
  Therefore, ${\rm Ext}^k(E',E_1)=0\;\forall k$ and
$${\rm Ext}^0(E_1,E')=V\not= 0,\quad{\rm Ext}^1(E_1,E')=W\not= 0,\quad
{\rm Ext}^2(E_1,E')=0.$$

   Consider the universal extension:
$$0\longrightarrow E'\longrightarrow E\longrightarrow W\otimes E_1
    \longrightarrow 0.$$
    By means of cohomology tables let us show that  $E$ is superrigid.
    The first table has the form:
$$\begin{array}{|c|ccccc|}
\hline
k&{\rm Ext}^k(E_1,E')&\rightarrow
 &{\rm Ext}^k(E_1,E)&\rightarrow
 &W\otimes {\rm Ext}^k(E_1,E_1)\\
\hline
0& V &&?  &&W  \\
1& W &&?  && 0 \\
2& 0 &&?  &&  0  \\
\hline
\end{array}\;. $$
Since the extension is universal, we see that the coboundary homomorphism
 $$W\rightarrow {\rm Ext}^1(E_1,E')$$ is isomorphism. Hence
$${\rm Ext}^1(E_1,E)={\rm Ext}^2(E_1,E)=0.$$
  The next tables have the form:
$$\begin{array}{|c|ccccc|}
\hline
k&{\rm Ext}^k(E',E')&\rightarrow
 &{\rm Ext}^k(E',E)&\rightarrow
 &W\otimes {\rm Ext}^k(E',E_1)\\
\hline
0&*  &&?  && 0 \\
1&0  &&?  && 0 \\
2&0  &&?  && 0 \\
\hline
\end{array}\; , $$
$$\begin{array}{|c|ccccc|}
\hline
k&W^*\otimes {\rm Ext}^k(E_1,E)&\rightarrow
 &{\rm Ext}^k(E,E)&\rightarrow
 &{\rm Ext}^k(E',E)\\
\hline
0&*  &&?  && * \\
1&0  &&?  && 0 \\
2&0  &&?  && 0 \\
\hline
\end{array}\; . $$

  Thus, $E$ is a superrigid bundle.

Now assume that $\VEC E,n$ is an arbitrary exceptional collection of bundles
 such that
 $$\mu_H(E_1)\leq\mu_H(E_2)\leq\cdots\leq\mu_H(E_n).$$
  Let us partition it into subcollections of bundles with equal $\mu_H$-
  slopes. Since all pairs in obtained subcollections are singular or
  zero-pairs, we see that there exists superrigid bundles
 $\FAM G,k$ constructed by these subcollections. Moreover,
  $\mu_H(G_i)<\mu_H(G_{i+1})$. Now let us remember that a pair
 $(E_i,E_j)$ of bundles has the type $hom$ provided
 $\mu_H(E_i)<\mu_H(E_{j})$, i.e ${\rm Ext}^k(E_i,E_j)=0$ for $k=1,2$ and
  ${\rm Ext}^k(E_j,E_i)=0$ for $k=0,1,2$.
  This yields that the bundle $E_i\oplus E_j$ is
  superrigid . Thus,    $\bigoplus G_i$ is the required bundle.

2. It follows from theorem \ref{l241} that a torsion free superrigid sheaf
   is a direct sum of
   $\mu_H$-semistable rigid bundles
   $F=\bigoplus\limits_{j=1}^{m}F_j$.
   Without loss of generality we can assume that
    $$\mu_H(F_j)<\mu_H(F_{j+1})\quad \forall i.$$
    Since $F$ is superrigid, we see that ${\rm Ext}^k(F_i,F_j)=0$
    for $k=1,2$ and for any pair  $i,j$. Besides, it follows from
    $\mu_H$-semistability of $F_i$ and the last inequality that
     ${\rm Hom}(F_j,F_i)=0$ for $j>i$.

     Taking into account theorem \ref{l241}, we obtain that each
    $F_j$ has the exceptional filtration
        $$Gr(F_j)=(x_{s_j-1}E_{s_j-1},\ldots,x_{s_{j-1}}E_{s_{j-1}}).$$
    Using the previous lemma and the proved fact (${\rm Ext}^k(F_j,F_i)=0$
    for $j>i$ and $k=0,1,2$), we see that the collection of the direct
    summands
    of all bundles $F_j$ (with the preservation of the order) is exceptional.
    This concludes the proof of the second theorem statement.

3. Since any torsion free rigid sheaf is locally free and direct summands
  are unambiguous determined, we see that it is sufficient for proving the
  third theorem statement in the case $F$ is an indecomposable superrigid
  bundle. But this case is proved in \ref{l255}.

  This completes the proof of the theorem.

%% FOLLOWING LINE CANNOT BE BROKEN BEFORE 80 CHAR
%%%%%%%%%%%%%%%%%%%%%%%%%%%%%%%%%%%%%%%%%%%%%%%%%%%%%%%%%%%%%%%%%%%%%%%%%%%%%%%%%%%
\section{Constructibility of Exceptional bundles.}

\subsection{Introduction to the Helix Theory.}

  In this section, following \cite{RH}, \cite{GH}, \cite{BH}, and \cite{AG},
  we tell about the general concepts and facts connected with exceptional
  sheaves on manifolds (see the definition of exceptional sheaves in 2.2)
  and exceptional objects in a derived category.

 The notion of exceptional bundles was introduced in the paper \cite{DP}.
 The main result of that paper is a description of Chern classes
 of samistable bundles on ${\BBB P}^2$ .  Exceptional bundles appeared
 there as some kind of boundary
 points.

 Further the exceptional bundles and the exceptional objects in a derived
 category of sheaves were researched on Rudakov's seminar in Moscow. It
  became clear that the exceptional objects (sheaves) organized as
 exceptional collections can generate the derived category of sheaves.
 Therefore, there exists a spectral sequence of Beilinson type associated
 with an exceptional collection. Let us remark that firstly
 a spectral sequence of  such type on ${\BBB P}^2$ appeared in \cite{D}. But
the general
 result independently on \cite{D} was obtained by A.~L.~Gorodentsev
 (\cite{GH}).

 The existence of the spectral sequence is the weighty reason for studying
 the exceptional sheaves. Besides, the exceptional bundles are
 interesting as bundles with a zero-dimensional moduli space.

 The next application of the exceptional bundles is a description of
 moduli spaces of semistable bundles. There exists such description for the
 case of projective plane (\cite{DP}) and of smooth 2-dimensional quadric
 (\cite{RS}).

  The helix theory is connected with number theory. Namely,
  A.~A.~Markov, in particular, studied solutions of the following
  Diophantus equation:  \begin{equation}\label{equmark} x^2+y^2+z^2=3xyz.
 \end{equation}
 (Now this equation is called the Markov equation and its solutions are
 called the Markov numbers.) It was proved that the Markov numbers
 coincide with ranks of the exceptional bundles on ${\BBB
 P}^2$ which form a foundation of a helix.

  A.~A.~Markov formulated the following hypothesis:

{\it Any triple of natural solutions of the equation (\ref{equmark})
  is uniquely determined by its maximal element.}

 This hypothesis can be formulated in terms of the exceptional bundles
 in the following way.

{\it Suppose $E$ and $F$ are exceptional bundles on ${\BBB P}^2$ with equal
 ranks; then either $E=F(n)$ or $E^*=F(n)$ for some natural $n$.}

 More detail can be found in  \cite{RP}.

 \vspace{2ex}

 Now let us pass to the helix theory. The definition of a helix and the
 first results about helices appeared in \cite{RP}. The definition of the
 helix is due to A.~L.~Gorodentsev and A.~N.~Rudakov. The word "helix"
 and the idea of considering a helix as an infinite system of bundles with
 some form of periodicity is due to W.~N.~Danilov.

 Below following \cite{RH}, we tell about axioms of helix theory.

\vspace{3ex}

 We shall consider pairs of objects of a category  {\gtc U} or elements of
 a set {\gtc U} .

 {\sc Definition.} A pair $(A,B)$ is called  {\it left admissible\/} if
 a certain pair $(L_AB,A)$ is defined. The pair $(L_AB,A)$ is called a
 {\it left mutation\/} of $(A,B)$ and the object $L_AB$ is called a
 {\it left shift\/} of $B$. Similarly, a pair
 $(A,B)$ is  {\it right admissible\/} if a certain pair
 $(B,R_BA)$ is defined. The pair  $(B,R_BA)$  in this case is called a
 {\it right mutation\/} of $(A,B)$ and the object $R_BA$ is a {\it right
 shift\/} of $A$.

 The axioms are following.

 (1L) If $(A,B)$ is left admissible then $(L_AB,A)$ is right admissible and
 $R_AL_AB=B$.

 (1R) If $(A,B)$ is right admissible then $(B,R_BA)$ is left admissible and
 $L_BR_BA=A$.

 (2L) Let $(A,B,C)$ be such that the pairs $(B,C),(A,L_BC)$ and
 $(A,B)$ are left admissible. Then the pairs $(A,C),(B',L_AC)$ are left
 admissible, where  $B'=L_AB$ and $L_AL_BC=L_{B'}L_AC$.

 (2R) Let $(A,B,C)$ be such that the pairs $(B,C),(A,B)$ and
  $(R_BA,C)$ are right admissible. Then the pairs $(A,C),(R_CA,B')$
 are right admissible, where
 $B'=R_CB$ and $R_CR_BA=R_{B'}R_CA$.

 The equalities in the axioms (2L) and  (2R) are usually called the
 {\it triangle equations\/}.

  It will be convenient to denote the object $L_AL_BC$, which appeared in
  (2L) by $L^{(2)}C$ and also to set $R^{(2)}A=R_CR_BA$. In the same way, if
  $\OVEC A,n$ is a system of objects we put
 $L^{(0)}A_s=A_s,L^{(1)}A_s=L_{A_{s-1}}A_s,\ldots,
 L^{(i)}A_s=L_{A_{s-i}}L^{(i-1)}A_s$, with the condition that the
 resulting pairs are left admissible. Analogous notation will be used for
 right mutations.

 {\sc Definition.} The collection $\{A_i|\;i\!\in\!\ZZ\}$ will be called a
 {\it helix of period $n$\/} if for all $s$ the following condition  holds

 {\sc HEL:} The pairs $(A_{s-1},A_s),(A_{s-2},L^{(1)}A_s),
 \ldots,(A_{s-n+1},L^{(n-2)}A_s)$ are left admissible and
  $L^{(n-1)}A_s=A_{s-n}$.

 Further  we shall assume that (1L), (1R), (2L) and (2R) are satisfied.
Then {\sc HEL} is equivalent to

 {\sc HEL':} The pairs $(A_{s-n},A_{s-n+1}),(R^{(1)}A_{s-n},A_{s-n+2}),
 \ldots,(R^{(n-2)}A_{s-n},A_s)$ are right admissible and
   $R^{(n-1)}A_{s-n}=A_s$.

 Each collection of the form $A_i,A_{i+1},\ldots,A_{i+n-1}$ is called a
 {\it foundation\/} of the helix $\{A_i\}$. Note that a helix is uniquely
 determined by any of its foundations.

 A collection $\{B_i|\;i\!\in\!\ZZ\}$ with
 $$\begin{array}{lcl}
 B_i=LA_{i+1}& \mbox{for}& i\equiv m-1(mod\; n),\\
 B_i=A_{i-1}& \mbox{for}& i\equiv m(mod\; n),\\
 B_i=A_{i}& \mbox{for}& i\not\equiv m, m-1(mod\; n),
 \end{array}
 $$
  is called a {\it left mutation\/} of the helix at $A_m$ and is
  denoted by $L_m$.

 A collection $\{C_i|\;i\!\in\!\ZZ\}$ with
 $$\begin{array}{lcl}
 C_i=RA_{i-1}& \mbox{for}& i\equiv m+1(mod\; n),\\
 C_i=A_{i+1}& \mbox{for}& i\equiv m(mod\; n),\\
 C_i=A_{i}& \mbox{for}& i\not\equiv m, m+1(mod\; n),
 \end{array}
 $$
  is called a {\it right mutation\/} of the helix at $A_m$ and is
  denoted by $R_m$.

  The basic fact about helices is the following.

 \TH{Theorem.}\label{l311} A right or left mutation of a helix is
 again a helix.

\ETH

\vspace{2ex}

  All applications of helices are based on this theorem.

  If we recall the triangle equations, we see that the mutations of helices
   define
  an action of the braid group on the set of all helices. One of the main
  questions of the helix theory is the question about an orbit number
  of this action.

  Let us return to exceptional sheaves on surfaces and define mutations of
  exceptional pair of sheaves. (The definition of exceptional pairs and
  their types can be found in section 2.3.)

\vspace{2ex}

 {\sc Lemma-definition.} 1. Let $(E,F)$ be an exceptional $hom$-pair
 of sheaves. Consider  the canonical homomorphisms
 $${\rm Hom}(E,F)\otimes E\stackrel{lcan}{\longrightarrow }F\qquad
 \mbox{and}\qquad E\stackrel{rcan}{\longrightarrow }{\rm Hom}(E,F)^*\otimes F
. $$
 If $lcan$ is an epimorphism then the pair $(E,F)$ is left admissible and
               $$L_EF=\ker lcan.$$
  Besides, the sheaf $L_EF$ is exceptional and the pair
 $(L_EF,E)$ is also exceptional.

 The pair $(E,F)$ is right admissible provided  $rcan$ is a monomorphism.
 In this case, $$R_FE={\rm coker}rcan.$$
 In addition, the sheaf $R_FE$  and the pair $(F,R_FE)$ are also exceptional.

 In both these cases the mutation of the pair $(E,F)$ is called
  {\it regular\/}.

 Suppose $lcan$ is a monomorphism; then the left shift of $F$ is defined as
 $L_EF={\rm coker}lcan$. (The pair $(L_EF,E)$ is exceptional as well.)

 The right shift of $E$ is defined as $R_FE=\ker rcan$ whenever
  $rcan$ is an epimorphism. (The pair $(F,R_FE)$ is exceptional in this case
  also.)

2. The $ext$-pair $(E,F)$ is both left and right admissible. The following
  universal extensions define the mutations of the $ext$-pair.
 $$0\longrightarrow F\longrightarrow L_EF\longrightarrow {\rm Ext}^1(E,F)
 \otimes E\longrightarrow 0,$$
 $$0\longrightarrow {\rm Ext}^1(E,F)^*\otimes F\longrightarrow R_FE
 \longrightarrow E\longrightarrow 0.$$
  In this case, as above, $L_EF$ and $R_FE$ are exceptional
 and $(L_EF,E)$, $(F,R_FE)$ are $hom$-pairs.

3. Both a left and a right mutation of a zero-pair is permutation of pair
  terms.

\vspace{2ex}

  It follows from this lemma that there are cases when left or right mutation
  of a $hom$-pair is not defined. Moreover, there are not mutations of a
  singular pair of sheaves.

 For overcoming of these limitations, following (\cite{AG}), let us pass to the
 bounded derived category of sheaves on the surface $S$ $\quad(D^b(S))$.
 Exceptional objects and collections in this category are defined in the
 same way as in the basic category of sheaves.

\vspace{2ex}

 {\sc Lemma-Definition.} Let $(E,F)$ be an exceptional pair in $D^b(S)$.
 Objects
 $L_EF$ and $R_FE$, which complete the canonical morphisms
 $$L_EF\longrightarrow R^{\cdot}{\rm Hom}(E,F)\otimes E\longrightarrow F
 \qquad\mbox{and}\qquad E\longrightarrow
 R^{\cdot}{\rm Hom}(E,F)^*\otimes F\longrightarrow R_FE$$
 up to the distinguished triangles, just as the pairs
  $(L_EF,E),\; (F,R_FE),$ are exceptional.

 \vspace{2ex}

 The category of sheaves is imbedded into  $D^b(S)$ by morphism $\delta$.
 Any exceptional sheaf remains exceptional under this imbedding. Mutations
 in basic and in derived category are connected in the following way. If
 an exceptional pair of sheaves $(E,F)$ is left admissible then the left
 shift of  $\delta (F)$ in the derived category is quasiisomorphic to
$\delta (L_EF)$.  That is, it is a complex with a unique nonzero
 cohomology coincided with $L_EF$, and vice versa. The similar
 statement holds true in the case of a right mutation. Thus we can reckon
 that any exceptional pair of sheaves is both left and right admissible.

 \TH{Theorem.}(Gorodentsev-Orlov.) Any exceptional object of $D^b(S)$
  is quasiisomorphic to an exceptional sheaf provided  $S$ is Del Pezzo
  surface. That is, all mutations of exceptional pairs of sheaves belong
  to the basic category.
\ETH

\TH{Theorem.}(Rudakov-Gorodentsev.) 1. The above defined mutations of
sheaves and exceptional objects of derived category satisfy the axioms
(1L),(2L),(1R) and (2R).

2. An exceptional collection remains exceptional whenever some its
pair of neighboring sheaves
   is replaced by a mutation of this pair. This procedure is called a
  {\it mutation of the collection\/}.
\ETH

{\sc Definition.} Let $\sigma =\VEC E,n$ be an exceptional collection of
 sheaves or objects of $D^b(S)$. It is  {\it full\/} provided $D^b(S)$
 is generated by $\sigma$, i.e. the set of all objects of $D^b(S)$
 can be obtained from the members of $\sigma$ by means of direct summing,
 tensoring and forming cones of all possible morphisms.

\vspace{2ex}

For example, the following collection of line bundles on ${\BBB P}^2$
$$({\cal O}_{{\BBB P}^2},{\cal O}_{{\BBB P}^2}(1),
{\cal O}_{{\BBB P}^2}(2))$$
is a full exceptional collection.

\TH{Theorem.}\label{l314} (Bondal.) Let $\sigma =\VEC E,n$ be an exceptional
 collection of sheaves or objects of derived category on a manifold $X$.
 Then the following statements are true.

1. If $\sigma$ is full then its left and right mutations are full
  collections.

2. The collection of the form
$$\sigma_\infty=\{E_i|i\!\in\!\ZZ,\; E_{i+sn}=E_i(-sK_X)\}$$
 is a helix of period $n$ if and only if  $\sigma$ is full.
\ETH

\vspace{2ex}

 We see that full collections are closely connected with helices.

\vspace{2ex}

  For writing the spectral sequence mentioned at the beginning of this
  section define dual collections.

  Let  $\sigma=\VEC E,k$ be an exceptional collection. A collection
  $({}^{\vee}E_{-k},\ldots,{}^{\vee}E_{-1},{}^{\vee}E_{0}),$
  where
$${}^{\vee}E_0=E_0,{}^{\vee}E_{-1}=LE_1,{}^{\vee}E_{-2}=L^{(2)}E_2,
\ldots,{}^{\vee}E_{-k}=L^{(k)}E_k$$ is called {\it left dual\/} to
$\sigma$.  A collection $(E^{\vee}_k,E^{\vee}_{k-1},\ldots,E^{\vee}_0),$
 where
$$E_0^{\vee}=R^{(k)}E_0,E_1^{\vee}=R^{(k-1)}E_1,E_2^{\vee}=R^{(k-2)}E_2,
\ldots,E_k^{\vee}=E_k$$ is called {\it right dual\/} to $\sigma$.

  In this notation the following theorem holds true.

\TH{Theorem.}(Gorodentsev.) Let  $Q$ be an exceptional object belonging to
 the  subcategory generated by an exceptional collection $\OVEC E,k$.
 Then there exists a spectral sequence $$E^{p,q}\quad\Longrightarrow
  H^{p+q}(Q),$$ in which the $E_1$ term has the form
$$E_1^{p,q}=\bigoplus\limits_{r+s=q}^{}{\rm Hom}^r_{D^b(S)}(E_{k-p},Q)
\otimes H^s({}^{\vee}E_{k-p}).$$
 In this case we say that the spectral sequence is associated with the
left dual collection.

 Similarly, a spectral sequence associated with right dual collection can be
 written.
\ETH

\TH{Corollary.}\label{l316} Let $\OVEC E,k$ be an exceptional collection of
 sheaves on the surface $S$. Suppose the left dual collection belongs to
 the basic category, i.e. each member of the left dual collection is a
 sheaf; then for any sheaf $Q$ belonging to the category generated by this
 collection there exists a spectral sequence $E^{p,q}$, with the
  $E_1$-term of the form
  $$E_1^{p,q}={\rm Ext}^{q-\Delta_p}(E_{-p},Q)\otimes {}^{\vee}E_{-p},$$
  where $\Delta_p$ is the number of nonregular mutations needed for
  constructing the sheaf ${}^{\vee}E_{-p}$.
Besides, there exists a spectral sequence $E^{p,q}$ with
   $E_1$-term
  $$E_1^{p,q}={\rm Ext}^{k-q-\Delta_p}(Q,E_{-p})^*\otimes E_{-p}^{\vee},$$
  where $\Delta_p$ is the number of nonregular mutations needed for
  constructing the sheaf $E_{-p}^{\vee}$.

  Both these sequences  converge to  $Q$ on the principal diagonal,
  i.e.
  $E^{p,q}_\infty =0$ for $p+q\not= 0$ and
  $$Gr(Q)=(E_\infty^{0,0},E_\infty^{-1,1},\ldots,E_\infty^{-n,n}) .$$
\ETH

\vspace{2ex}

  The helix theory has the following problems.

  1. Are there full exceptional collections on a given manifold?

  2. How many orbits has the braid group under the action on the set of all
  helices?
  We say that all helices (full exceptional collections) are
   {\it constructible\/} provided the orbits is unique.

  3. Does an orbitrary exceptional collection belong to a foundation of a
  helix? In other words, is there full exceptional collection containing a
  given exceptional collection? We say that the exceptional sheaves are
  {\it constructible\/} whenever the solutions of the second and the last
  problems are positive.

  4. We can consider an action of braid group on the set of exceptional
  collections generated one and the same derived subcategory of $D^b(X)$.
  How many orbits has this action?

  5. Description of stable subgroups of braid group action.

  Full collections were found on  ${\BBB P}^n$, Del Pezzo surfaces,
   $G(2,4)$. Besides, the following theorem proved by D.~Orlov in \cite{O}.

  \TH{Theorem.}\label{l317}(Orlov.) 1. Let ${\BBB P}(E)\rightarrow M$
  be a projectivisation of a vector bundle on a manifold $X$. Suppose there
  is a full exceptional collection on $X$; then there exists such
  collection on ${\BBB P}(E)$.

  2. Let $\tilde X$ be obtained from $X$ by blowing up a smooth regular
  submanifold $Y$. Suppose there are full exceptional collections on $X$
  and on $Y$; then there is a full exceptional collection on $\tilde X$ as
  well.
\ETH

\vspace{2ex}

  In the papers \cite{RP},\cite{RQ},\cite{OK} it was proved that all
  exceptional
  sheaves and all helices on Del Pezzo surfaces are constructible. The
  constructibility of helices on ruled surfaces with the rational base
  and on ${\BBB P}^3$ was  proved in (\cite{NOG}).

\vspace{2ex}

  In the last part of our paper we shall prove the following theorem.

\TH{Theorem.}\label{l318}1. Let $ \sigma$ be an exceptional collection of
  bundles on a smooth projective surfaces $S$ with anticanonical class
  without base components and $K_S^2>0$. Suppose rank of each bundle of this
  collection is greater then 1; then there is a full exceptional collection
  $\tau$ such that $\sigma $ is a subcollection of $\tau$. Moreover, $\tau$
  can be obtained by mutations from the basic full collection. In other
  words, all helices on $S$ are constructible.

2. The condition about ranks can be rejected provided  $K_S^2>1$.
\ETH

%% FOLLOWING LINE CANNOT BE BROKEN BEFORE 80 CHAR
%%%%%%%%%%%%%%%%%%%%%%%%%%%%%%%%%%%%%%%%%%%%%%%%%%%%%%%%%%%%%%%%%%%%%%%%%%%%%%%%%
\subsection{Restriction of Superrigid Bundles to an Exceptional Curve.}

   Let us remember that we deal with the surface $S$ with anticanonical
   class $H=-K_S$ without base components. This means that $H$ is nef.

   In the beginning of section 2.4 we limited class of considered
   surfaces by the condition: $K_S^2>0$. Using the description of surfaces
   with numerically effective anticanonical class from section 2.1, we see
   that the surfaces satisfying such condition are the following:  ${\BBB
   P}^2, {\BBB P}^1\times {\BBB P}^1$  or surfaces obtained from ${\BBB
  P}^2$ by blowing up at most 8 points. Among them only the quadric has
  not blowing down onto ${\BBB P}^2$.

   The helix theory on the quadric is known (\cite{RQ}, \cite{ZUZ}).
   Further we can assume that the surface $S$ satisfy the following
   conditions.

1. $K_S^2>0.$

2. There is a blowing down of  $S$ onto ${\BBB P}^2$.

   Suppose $S$ is obtained from ${\BBB P}^2$ in the following way
$$S\stackrel{\sigma_d}{\longrightarrow}S_{d-1}
\stackrel{\sigma_{d-1}}{\longrightarrow}\cdots
\stackrel{\sigma_{2}}{\longrightarrow}S_1
\stackrel{\sigma_{1}}{\longrightarrow}S_0={\BBB P}^2,$$
  where $\sigma_i$ is blow up a point $p_{i-1}\in S_{i-1}$ and $d\leq 8$.
  By definition, put $e_i=\sigma_i^{-1}(p_{i-1}$.

It is clear that $e_i$ are exceptional $-1$-curves for all $i$ and $e_d$
is irreducible. We see that $e_d$ is a smooth rational curve.

   It is known, that the divisors $h,e_1,\ldots,e_d$ generate the group
   ${\rm Pic}(S)$ (here $h$ is the preimage of a line on
   ${\BBB P}^2$). Besides,

   $$he_i=e_ie_j=0\quad\mbox{for}\quad i\not=j\quad\mbox{and}
   \quad e_i^2=-1.$$
   $$K_S=-3h+\sum\limits_{i=1}^{d}e_i.$$

\TH{Remark.}\label{l321} The divisor $h$ is numerically effective.
\ETH

\PR In reality, a line on ${\BBB P}^2$ has not base components. Hence, its
  preimage is base set free as well. Therefore, a cup product $h$ with any
  curve on $S$ is nonnegative.

\TH{Lemma.}\label{l322} Let $E$ and $F$ be exceptional bundles on $S$
  with equal $\mu_H$-slopes and let $e=e_d$ be an irreducible exceptional
  curve. Suppose $E=F$ or $c_1(E)-c_1(F)=C$ is -2-curve; then

  either ${\rm Ext}^2(E,F(-e))=0$

  or $K_S^2=1$ and $(E,F)$ is exceptional pair of the form
  $\Bigl({\cal O}_S(D),{\cal O}_S(D+e+K_S)\Bigr)$, where $D$ is some divisor
  of   ${\rm Pic}(S)$.
\ETH

\PR By Serre duality theorem,
$${\rm Ext}^2(E,F(-e))^*\cong{\rm Hom}(F,E(e+K_S)).$$
  Suppose $K_S^2>1$; then
$$\mu_H(E(e+K_S))=\mu_H(E)-K_S\cdot e-K_S^2<\mu_H(E)=\mu_H(F).$$
  Hence the equality ${\rm Hom}(F,E(e+K_S))=0$ follows from
$\mu_H$-stability of exceptional bundles on $S$ and \ref{l115}

  Now suppose $K_S^2=1$ and ${\rm Hom}(F,E(e+K_S))\not= 0$. It follows
  from the equality $\mu_H(F)=\mu_H(E(e+K_S))$  and \ref{l116} that
  there exists the exact triple:
\begin{equation}\label{equ91}
0\longrightarrow F\longrightarrow E(e+K_S)\longrightarrow Q
\longrightarrow 0,
\end{equation}
  where $Q$ is a torsion sheaf. Denote by $r$ the rank of bundles $E$ and
  $F$. (Let us remember that $r(E)=r(F)$). We get
  $$c_1(Q)=c_1(E)-c_1(F)+r(e+K_S)=C+r(e+K_S).$$
 Recall that the first Chern class of a torsion sheaf must be "nonnegative",
 i.e. either effective or trivial. Assume that $E=F$; then
  $c_1(Q)=r(e+K_S)$. It is not possible, since this divisor is ineffective.

 Assume that $E\not= F$. Then, by the
 lemma conditions, $C=c_1(E)-c_1(F)$ is
 -2-divisor such that  $C\cdot K_S=0$ (recall that
$\mu_H(E)=\mu_H(F)$ and $r(E)=r(F)$).  Such divisors were described by
 Yu.~I.~Manin in \cite{MAN}. Using his results, we can state that if
 $C=ah-\sum\limits_{i=1}^{d}b_ie_i$ then $|a|\leq 3$. Moreover,
 $C=3h-e_j-\sum\limits_{i=1}^{d}e_i$ whenever $a=3$.

 We assume that sequence (\ref{equ91}) exists. In this case the
 divisor $C+r(e+K_S)$ is effective. Whereby, the cup product $h\cdot
 (C+r(e+K_S))$ is nonnegative (\ref{l321}). Thus,
 $C=3h-e_j-\sum\limits_{i=1}^{d}e_i$ and $r=1$.

 We have $C+r(e+K_S)=2e_d-2e_j$
(recall that $e=e_d$). The curve $e_d$ is irreducible. Hence, $2e_d-2e_j$
  is ineffective when $j\not= d$. Therefore, $e_d=e_j$ and $C=-K_S-e$.
 Thus the pair $(E,F)$ is equal  to
   $$\Bigl({\cal O}_S(D),{\cal O}_S(D+e+K_S)\Bigr).$$
 This concludes the proof.

 \TH{Corollary of the proof.}\label{l323} Suppose $C$ is -2-divisor with
 $C\cdot K_S=0$ and $e=e_d$; then the divisor
 $C+e+K_S$ is nonpositive.
\ETH

\TH{Lemma.}\label{l324} Let $(E,F)$ be an exceptional pair of bundles on $S$
  with $\mu_H(E)<\mu_H(F)$ and $e=e_d$ the irreducible exceptional curve.
  Then
        $${\rm Ext}^2(E,F(-e))=0.$$
\ETH

\PR By Serre duality theorem we have
$${\rm Ext}^2(E,F(-e))^*\cong {\rm Hom}(F,E(e+K_S)).$$
 But,
 $$\mu_H(E(e+K_S))=\mu_H(E)+1-K_S^2\leq \mu_H(E)<\mu_H(F)$$
 and the proof follows from  $\mu_H$-stability of exceptional bundles
 on $S$ and lemma \ref{l115}.

\TH{Lemma.}\label{l325} Suppose $E$ and $F$ are rigid $\mu_H$-semistable
  bundles on $S$. Assume that they have exceptional filtrations
$$Gr(E)=(x_nE_n,x_{n-1}E_{n-1},\ldots,x_1E_1),\qquad
  Gr(F)=(y_mF_m,y_{m-1}F_{m-1},\ldots,y_1F_1).$$
  In addition we assume that the following conditions hold.

1. ${\rm Ext}^k(F,E)=0 \quad \forall k=0,1,2$.

2. $\mu_H(E)<\mu_H(F)<\mu_H(E)+K_S^2$.

3. Provided $K_S^2=1$ the exceptional collections
 $\VEC E,n$, $\VEC F,m$ have not pairs of the form
 $\Bigl({\cal O}_S(D),{\cal O}_S(D+e+K_S)\Bigr)$, where
  $D\!\in\!{\rm Pic}(S)$, and $e=e_d$ is irreducible exceptional curve.

  Then the restrictions of $E$ and $F$ to $e$ have the form
 $$E'=E|_e=\alpha {\cal O}_e(s-1)\oplus \beta {\cal O}_e(s),\quad
    F'=F|_e=\gamma {\cal O}_e(s-1)\oplus \delta {\cal O}_e(s)\oplus
    \epsilon {\cal O}_e(s+1),$$
 where $\alpha ,\beta, \gamma ,\delta,\epsilon$ are nonnegative integer
 numbers with  $\gamma \epsilon =0$.
\ETH

\PR It follows from the lemma conditions and \ref{l241} that all pairs
$(E_i,E_j)$ and $(F_i,F_j)$ for $i<j$ are exceptional singular or zero-pairs.
 By proposition \ref{l232} we have
$$\mu_H(E_i)=\mu_H(E_j),\; \mu_H(F_i)=\mu_H(F_j)$$
 and remainders of the first Chern classes
$$c_1(E_j)-c_1(E_i)\mbox{ and } c_1(F_j)-c_1(F_i)$$
  are -2-divisors. Since among these pairs there are not the pairs of the
  kind
$$\Bigl({\cal O}_S(D),{\cal O}_S(D+e+K_S)\Bigr)$$ (in the case
  $K_S^2=1$), it follows from lemma  \ref{l322} that
$${\rm Ext}^2(E_i,E_j(-e))={\rm Ext}^2(E_i,E_i(-e))=
{\rm Ext}^2(F_i,F_i(-e))={\rm Ext}^2(F_i,F_j(-e))=0$$
 for any pair $i,j$. Thus the equalities
$${\rm Ext}^2(E,E(-e))={\rm Ext}^2(F,F(-e))=0$$
  follow from \ref{l322}

 Since the bundles $E$ and $F$ are rigid, using the exact
 triples $$0\longrightarrow E^*\otimes E(-e)\longrightarrow E^*\otimes E
\longrightarrow (E^*\otimes E)|_e\longrightarrow 0,$$
$$0\longrightarrow F^*\otimes F(-e)\longrightarrow F^*\otimes F
\longrightarrow (F^*\otimes F)|_e\longrightarrow 0$$
  we get,
      ${\rm Ext}^1(E',E')={\rm Ext}^1(F',F')=0$.

 By Grothendieck theorem \cite{OKON}, any bundle on a projective line (in
 particular, $E'$ and $F'$ on $e$) is a direct sum of line bundles. From
 rigidity of $E'$ and $F'$ and Bott formula,
  which calculates cohomologies of
 line bundles on the projective line (\cite{BOTT}) we obtain
 $$E'=E|_e=\alpha {\cal O}_e(s-1)\oplus \beta {\cal O}_e(s),\quad
    F'=F|_e=\gamma {\cal O}_e(s'-1)\oplus \delta {\cal O}_e(s').$$

  Using the first and the second conditions of the lemma, let us show that
                   $$s\leq s'\leq s+1.$$
  Note that from condition 1 and proposition \ref{l231} it follows that each
  pair $(E_i,F_j)$ is exceptional. Besides,
         $$\mu_H(E_i)<\mu_H(F_j)<\mu_H(E_i)+K_S^2.$$
        Applying lemma \ref{l324} to the pairs $(E_i,F_j)$, we get
   ${\rm Ext}^2(E_i,F_j(-e))=0$. This means that
         ${\rm Ext}^2(E,F(-e))=0$.

 By virtue of the inequalities on $\mu_H$-slopes and \ref{l231}, the pairs
  $(E_i,F_j)$ have the type $hom$. In particular,
  ${\rm Ext}^1(E_i,F_j)=0$. Hence we have ${\rm Ext}^1(E,F)=0$. Now it
  follows from a long exact cohomology sequence associated with the
  restriction sequence to the exceptional curve
$$0\longrightarrow E^*\otimes F(-e)\longrightarrow E^*\otimes F
   \longrightarrow (E^*\otimes F)|_e\longrightarrow 0$$
 that ${\rm Ext}^1(E',F')=0$.

  By Serre duality and the first lemma condition we get
 $${\rm Ext}^k(E,F(K_S))=0\qquad\mbox{for}\qquad k=0,1,2.$$
  The second lemma condition yields the inequality
$$\mu_H(F(K_S))<\mu_H(E)<\mu_H(F(K_S))+K_S^2.$$
  Repeating the reasoning for rigid bundles $F(K_S)$ and $E$,
  we see that ${\rm Ext}^1(F(K_S)|_e,E')=0$.

  Note that $F(K_S)|_e=F'(-1)$, i.e.
${\rm Ext}^1(F'(-1),E')=0$. Now, the inequality
        $(s\le s'\le s+1)$
 follows from Bott formula. This completes the proof.

\TH{Corollary.}\label{l326} Assume that an ordered collection of
 $\mu_H$-semistable rigid
  bundles $$\VEC E,m$$ satisfies the following conditions.

1. ${\rm Ext}^k(E_i,E_j)=0$ for $j>i,\; k=0,1,2$.

2. $\mu_H(E_1)<\mu_H(E_2)<\cdots <\mu_H(E_m)<\mu_H(E_1)+K_S^2$.

3. Provided $K_S^2=1$ the exceptional collections corresponding to
   the exceptional filtrations of all  $E_i$ have not pairs of the form
 $\Bigl({\cal O}_S(D),{\cal O}_S(D+e+K_S)\Bigr)$, where $e=e_d$ is the
  irreducible exceptional curve.

 Then there is a number $i$ such that
 $$\Bigl(E_i\oplus \cdots\oplus E_m\oplus E_1(-K_S)
 \oplus \cdots\oplus E_{i-1}(-K_S)
   \Bigr)|_e=\alpha {\cal O}_e(s)\oplus \beta {\cal O}_e(s+1).$$
\ETH

\PR
  We shall say that an
  ordered pair of rigid $\mu_H$-semistable bundles $(E,F)$
  on $S$ {\it has a zero type of  decomposition\/} whenever
 $$(E\oplus F)|_e=\alpha {\cal O}_e(s)\oplus \beta {\cal O}_e(s+1)
.$$
 It {\it has a first type of  decomposition \/} whenever

 $$E|_e=\alpha {\cal O}_e(s)\oplus \beta {\cal O}_e(s+1),\qquad
   F|_e=\gamma {\cal O}_e(s+1)\oplus \delta {\cal O}_e(s+2),$$
  with $\alpha\cdot\delta\not= 0$.

 From the previous lemma it follows that each pair from our collection has
 either the zero or the first type of  decomposition.

  We see that the
  lemma statement holds true provided the pair $(E_1,E_i)$ has the
  zero type of  decomposition for all $i$.

  In the converse case, denote by $i$ the minimal number such that the pair
   $(E_1,E_i)$ has the first type of  decomposition. Note that,
   in this case, $\forall j\!<\!i\!\le k$ the pair $ (E_j,E_k)$ has the
   first type
   of   decomposition and the pair $(E_s,E_l)$  has the zero type of
   decomposition whenever $i\le s<l$ or $s<l\le i$.
    Besides, if a pair $(E,F)$ has the first type
   of decomposition then the pair $(F,E(-K_S))$ has the zero type of
   decomposition.

   Thus, each pair of the collection Следовательно, каждая пара в наборе
$$\Bigl(E_i\oplus \cdots\oplus E_m\oplus E_1(-K_S)\oplus\cdots
\oplus E_{i-1}(-K_S)  \Bigr)$$
has the zero type of a   decomposition. This completes the proof.
%

%% FOLLOWING LINE CANNOT BE BROKEN BEFORE 80 CHAR
%%%%%%%%%%%%%%%%%%%%%%%%%%%%%%%%%%%%%%%%%%%%%%%%%%%%%%%%%%%%%%%%%%%%%%%%%%%%%%%%%
\subsection{Equivalence of Collections and the Key Exact Sequence.}

{\sc Definition.} We shall say that an exceptional collection
   $\sigma =\VEC E,k$ (of sheaves or of objects in $D^b(S)$) on $S$
   is {\it constructible\/} whenever there is a full exceptional
   collection $(E_1,\ldots,E_k,E_{k+1},\ldots,E_n)$ containing $\sigma $
    such that it is obtained from the basic collection $$\sigma_0=({\cal
    O}_S,{\cal O}_S(h),{\cal O}_S(2h), {\cal O}_{e_1}(-1),\ldots,{\cal
        O}_{e_d}(-1))$$ by mutations. Here $h$ is the preimage of a line
   on ${\BBB P}^2$ and $e_i$ are the blow up divisors.  (It follows from
    \cite{O} that the basic collection is exceptional and full.)

    We say that an exceptional collection $\sigma$ is {\it equivalent\/} to
   an exceptional collection $\tau$ whenever the following condition
   holds true.  The collection $\sigma $ is constructible if and only if
   $\tau$ is constructible.

\TH{Lemma.}\label{l331} a) Suppose an exceptional collection $\sigma$ is
   obtained from an exceptional collection $\tau$ by mutations; then these
   collections are equivalent.

    b) An exceptional collection $\VEC E,k$ is equivalent to the following
    collections:
    $$(E_k(K_S),E_1,\ldots,E_{k-1})\qquad\mbox{ and }\qquad
      (E_2,\ldots,E_k,E_1(-K_S)).$$
\ETH

\PR a) Assume that an exceptional collection $\sigma=\VEC E,k$ is obtained
  from  $\tau=\VEC F,k$ by mutations. Since all mutations of collections are
  invertible (see the axioms of the helix theory), we can assume that
  $\tau$ is also obtained from $\sigma $ by mutations. Suppose $\sigma $ is
  constructible, i.e. there is a full exceptional collection
 $\sigma'=(E_1,\ldots,E_k,E_{k+1},\ldots,E_n)$  obtained from the basic
 collection by mutations. Then the exceptional collection
    $\tau'=(F_1,\ldots,F_k,E_{k+1},\ldots,E_n)$  is also obtained from the
    basic collection by mutations. Therefore, $\tau'$ is full (\ref{l314}).
    Besides, the basic collection and $\tau'$ belong to one and the same
    orbit of the braid group action. Thus $\tau$ is also constructible.

    b) For proving the second lemma statement it is sufficient to check that
    the collections $\sigma=\VEC E,k$ and $(E_2,\ldots,E_k,E_1(-K_S))$ are
    equivalent. Suppose $\sigma $ is constructible and
    $\sigma_1=(E_1,\ldots,E_k,E_{k+1},\ldots,E_n)$ is a full exceptional
    collection corresponding to it. By theorem \ref{l314}, if follows from
    full of the collection $\sigma_1$ that it is a foundation of a
    helix and $E_1(-K_S)=R^{n-1}E_1$. That is, the collection
    $$\sigma_2=(E_2,\ldots,E_k,E_{k+1},\ldots,E_n,E_1(-K_S))$$
    is equivalent to  $\sigma_1$. Now we shift each of the sheaves
    $E_n,E_{n-1},\ldots,E_{k+1}$ to the right over $E_1(-K_S)$ to obtain the
    full collection
    $$\tau_1=(E_2,\ldots,E_k,E_1(-K_S),RE_{k+1},RE_{k+2},\ldots,RE_n),$$
    equivalent to $\sigma_2$. Thus $\tau$ is constructible
    as well.

    Since all operations are invertible, we see that the collection $\sigma$
    is equivalent to $\tau$.

    \vspace{2ex}

    {\sc Notation.} Let $\sigma=\VEC E,k$ be an exceptional collection of
    bundles. By definition, put
  $$\mu_-(\sigma)=min\{\mu_H(E_i)\},\qquad\mu_+(\sigma)=max\{\mu_H(E_i)\}.$$

\TH{Lemma.}\label{l332} For any exceptional collection of bundles
    $\sigma=\VEC E,k$ there exists an exceptional collection of bundles
    $\tau=\VEC F,k$ equivalent to $\sigma$ such that
    $$\mu_-(\sigma)\le\mu_-(\tau)=\mu_H(F_1)\le\ldots\le
      \mu_H(F_k)=\mu_+(\tau)\le\mu_+(\sigma).$$
\ETH

    Further we shall say that the exceptional collection of bundles
    $\VEC F,n$ is {\it $hom$-collection\/} whenever
 $$\mu_H(F_1)\le  \mu_H(F_2)\le\cdots\le\mu_H(F_n).$$

\PR Let  $s$ be the minimal number such that $\mu_H(E_s)>\mu_H(E_{s+1})$.
        It follows from proposition  \ref{l231} that the exceptional pair
    $(E_s,E_{s+1})$ has the type $ext$. Consider the left mutation of this
    pair
    $$0\longrightarrow E_{s+1}\longrightarrow L_{E_s}E_{s+1}
    \longrightarrow E_s\otimes {\rm Ext}^1(E_s,E_{s+1})\longrightarrow 0.$$
    Since the sheaves $E_s$ and $E_{s+1}$ are locally free, we see that
    $L_{E_s}E_{s+1}$ is locally free as well. In view of $\mu_H$-stability
    of exceptional bundles, we have
    $$\mu_H(E_{s+1})<\mu_H(L_{E_s}E_{s+1})<\mu_H(E_s).$$

    Now suppose that $\mu_H(L_{E_s}E_{s+1})<\mu_H(E_{s-1})$; then we do
    the left mutation of this $ext$-pair, et cetera...

    It is clear that after a finite number of the mutations  we shall obtain
    an exceptional collection $\sigma'$ equivalent to the original one
    such that
    $$\mu_-(\sigma)\le\mu_-(\sigma')<\mu_+(\sigma')\le\mu_+(\sigma).$$
     Moreover, if we denote by $s'$ the minimal number such that
     $$\mu_H(E'_{s'})>\mu_H(E'_{s'+1})$$ then $s'>s$. This implies the lemma
     statement.

\TH{Lemma.}\label{l333} For any exceptional collection of bundles
    $\sigma=\VEC E,k$ there exists an exceptional $hom$-collection of
    bundles    $\tau=\VEC F,k$ equivalent to $\sigma$ such that
                  $$\mu_+(\tau)-\mu_-(\tau)<K_S^2.$$
\ETH

\PR By definition, put $\Delta\mu(\sigma)=\mu_+(\sigma)-\mu_-(\sigma).$
     Assume that $\Delta\mu(\sigma)>K_S^2$. By lemma \ref{l332}, without loss
     of generality it can be assumed that $\sigma$ is a $hom$-collection.
     We have
     $$\mu_-(\sigma)=\mu_H(E_1),\qquad\mu_+(\sigma)=\mu_H(E_n).$$
     Since $\Delta\mu(\sigma)>K_S^2$ and
     $\mu_H(E_1(-K_S))=\mu_H(E_1)+K_S^2$, we see that there is a number $s$
     such that
     $$\mu_H(E_{s-1})\le \mu_H(E_{1}(-K_S))<\mu_H(E_{s}).$$

     The collection
      $\sigma_1=(E_s,\ldots,E_n,E_1(-K_S),\ldots,E_{s-1}(-K_S))$
     is equivalent to $\sigma$ and it has the following $\mu_h$-slope
     limits:
     $$\mu_-(\sigma_1)=\mu_H(E_1(-K_S))=\mu_H(E_{1})+K_S^2,$$
     $$\mu_+(\sigma_1)=max\{\mu_H(E_{s-1}(-K_S)),\mu_H(E_{n})\}.$$

     Suppose $\mu_+(\sigma_1)=\mu_H(E_{s-1}(-K_S))$; then
      $s>1$ and
     $$\Delta\mu(\sigma_1)=\mu_H(E_{s-1})-\mu_H(E_{1})\leq K_S^2.$$
     Hence, ordering the collection $\sigma_1$ as in \ref{l332}, we obtain
     the $hom$-collection $\sigma_2$ equivalent to the original one such that
     $\Delta\mu(\sigma_2)\leq K_S^2.$

     Suppose $\mu_+(\sigma_1)=\mu_H(E_{n})=\mu_+(\sigma);$ then ordering the
     collection $\sigma_1$ by $\mu_H$-slopes, we obtain the $hom$-collection
     $\sigma_3$ equivalent to the original one such that
     $\Delta\mu(\sigma_3)\leq \Delta\mu(\sigma)- K_S^2.$

     Doing this operation for several times we all the same obtain a
     $hom$-collection equivalent to $\sigma$ with  $\Delta\mu\leq K_S^2.$

     Now let us assume that $\Delta\mu(\sigma)= K_S^2$. Denote by $s$ the
        minimal number such that  $\mu_H(E_{s})<\mu_H(E_{s+1})$.
         Consider the equivalent collection
     $$\tau=(E_{s+1},\ldots,E_n,E_1(-K_S),\ldots,E_s(-K_S)).$$
     By the choice of $s$ we have
     $$\mu_+(\tau)=\mu_H(E_{n})=\mu_H(E_{1}(-K_s))=\ldots=\mu_H(E_{s}(-K_S))=
     \mu_+(\sigma),$$
     $$\mbox{and}\qquad \mu_-(\tau)=\mu_H(E_{s+1})>\mu_-(\sigma).$$
     In other words, $\tau$ is the $hom$-collection with
     $\Delta\mu(\tau)<K_S^2$.
     This completes the proof.

\TH{Lemma.}\label{l334} Let $\sigma=\VEC E,k$  be an exceptional collection
   of bundles on the surface $S$ with $K_S^2\ge 1$. In addition assume that
   in the case $K_S^2=1$ this collection has not pairs of the form
     $({\cal O}_S(D),{\cal O}_S(D+e+K_S))$, where $D$ is some divisor and
    $e=e_d$ is the exceptional rational curve.

    Then there is an exceptional $hom$-collection $\VEC F,k$ equivalent to
   $\sigma$ such that the associated with it superrigid bundle $F$
     $\Bigl(Gr(F)=(x_nF_n,\ldots,x_1F_1)\Bigr)$
   is included in the exact sequence
      \begin{equation}\label{orl}
      0\rightarrow G\rightarrow F\rightarrow
    {\rm Hom}(F,{\cal O}_e(-1))^*\otimes {\cal O}_e(-1)\rightarrow 0,
    \end{equation}
    where $G$ is a superrigid bundle with
    $${\rm Ext}^k(G,{\cal O}_e(-1))=0\qquad \forall k.$$
     \ETH

     \PR By lemma \ref{l333} there is a $hom$-collection
     $\tau=\VEC {E'},k$ equivalent to the original one such that
     $\mu_+(\tau)-\mu_-(\tau)<K_S^2.$

     We shall separate this collection into groups of bundles with equal
     $\mu_H$-slopes. We shall construct from these groups superrigid
   $\mu_H$-semistable bundles $\FAM {\cal E},m$
     (see theorem \ref{l251}).

      We have
      ${\rm Ext}^k({\cal E}_j,{\cal E}_i)=0$ for $j>i,\quad k=0,1,2;$
        and
     $$\mu_H({\cal E}_{1})<\mu_H({\cal E}_{2})<\cdots<
       \mu_H({\cal E}_{m})<\mu_H({\cal E}_{1})+ K_S^2.$$
       In view of corollary \ref{l326} there is a number $i$ such that
$$({\cal E}_{i}\oplus\cdots\oplus{\cal E}_{m}\oplus{\cal E}_{1}(-K_S)\oplus
\cdots{\cal E}_{i-1}(-K_S))|_e=\alpha{\cal O}_e(d-1)\oplus\beta{\cal O}_e(d).
$$
    Hence there is a number  $j$ such that the superrigid bundle
    $\bar F$ associated    with the $hom$-collection
$$\tau'=(E'_j,\ldots,E'_k,E'_1(-K_S),\ldots,E'_{j-1}(-K_S))$$
     satisfy the condition
     $$\bar F|_e=\alpha{\cal O}_e(d-1)\oplus\beta{\cal O}_e(d).$$

     Thus the superrigid bundle $F$ constructed from the exceptional
     $hom$-collection
$$\tau''=\tau'(dK_S)=(E'_j(dK_S),\ldots,E'_k(dK_S),E'_1((d-1)K_S),
\ldots,E'_{j-1}((d-1)K_S))$$
   restricts to the curve $e$ in the following way:
$$F|_e=\alpha{\cal O}_e(-1)\oplus\beta{\cal O}_e.$$
     It is easily shown that the collection $\tau''$ is equivalent to the
        original one $\sigma$.

    The following equalities can be obtained by direct calculations.
    $$
    \mbox{$ h^i(F,{\cal O}_e(-1))=$}
      \left\{
       \begin{array}{ccc}
       \alpha& \mbox{ for }& i=0\\
       0     & \mbox{ for }& i>0
       \end{array}
      \right.\; ,
    $$
    $$
    \mbox{$ h^i({\cal O}_e(-1),F)=$}
      \left\{
       \begin{array}{ccc}
       \beta& \mbox{ for }& i=1\\
       0     & \mbox{ for }& i\not= 1
       \end{array}
      \right.\; .
    $$
   Consider the canonical map:
   $$F\stackrel{rcan}{\longrightarrow }{\rm Hom}(F,{\cal O}_e(-1))^*
      \otimes {\cal O}_e(-1).$$
   Since the restriction of this map to the curve $e$ is an epimorphism,
   we see that exact sequence (\ref{orl}) holds true.

   The sheaf $G$ from this sequence, as a subsheaf of a bundle, has not
   torsion. For calculating  its cohomologies let us consider cohomology
   tables corresponding to exact sequence (\ref{orl}). For this, denote
   by ${\cal L}$ the torsion sheaf ${\cal O}_e(-1)$.

   $$\begin{array}{|c|ccccc|}
   \hline
   k&{\rm Hom}(F,{\cal L})\otimes
   {\rm Ext}^k({\cal L},{\cal L})&\rightarrow
    &{\rm Ext}^k(F,{\cal L})&\rightarrow
    &{\rm Ext}^k(G,{\cal L})\\
   \hline
   0& {\rm Hom}(F,{\cal L})\otimes \CC &&{\rm Hom}(F,{\cal L})
    &&?  \\
   1& 0 &&0  && ? \\
   2& 0 &&0  && ?  \\
   \hline
   \end{array} $$
   $$\begin{array}{|c|ccccc|}
   \hline
   k&{\rm Hom}(F,{\cal L})\otimes
   {\rm Ext}^k({\cal L},F)&\rightarrow
    &{\rm Ext}^k(F,F)&\rightarrow
    &{\rm Ext}^k(G,F)\\
   \hline
   0& 0 &&*  &&?  \\
   1& * &&0  && ? \\
   2& 0 &&0  && ?  \\
   \hline
   \end{array} $$
   $$\begin{array}{|c|ccccc|}
   \hline
   k&
   {\rm Ext}^k(G,G)&\rightarrow
    &{\rm Ext}^k(G,F)&\rightarrow
    &{\rm Hom}(F,{\cal L})^*\otimes{\rm Ext}^k(G,{\cal L})\\
   \hline
   0& ? &&*  &&0  \\
   1& ? &&0  && 0 \\
   2& ? &&0  && 0  \\
   \hline
   \end{array} $$
 This concludes the lemma proof.

\vspace{2ex}

  I want to remark that the idea of the construction of exact sequence
   (\ref{orl}) on Del Pezzo surface with an exceptional bundle as $F$
   pertains to D.~O.~Orlov (\cite{OK}).

%% FOLLOWING LINE CANNOT BE BROKEN BEFORE 80 CHAR
%%%%%%%%%%%%%%%%%%%%%%%%%%%%%%%%%%%%%%%%%%%%%%%%%%%%%%%%%%%%%%%%%%%%%%%%%%%%%%%%%%
\subsection{Category Generated by a Pair.}

  In the previous section we constructed from an exceptional collection
  $\sigma$ of bundles the $hom$-collection  $\tau=\VEC F,k$ equivalent to
  $\sigma$ such that the superrigid bundle associated with  $\tau$ is
  included in exact sequence (\ref{orl}). In the next section using double
  induction, we shall show that from this sequence the constructibility of
  the collection $\tau$ follows.

  Here we check a base of one of the inductions. Namely we prove the
  following  proposition.

\TH{Proposition.}\label{l341} Suppose a superrigid sheaf $F$ on the surface
  $S$ is included in the exact sequence
 $$0\longrightarrow y_1G_1\longrightarrow F\longrightarrow y_0G_0
 \longrightarrow 0,$$
  where $(G_0,G_1)$ is an exceptional pair, $y_i$ are positive integer and
  $G_1$ is a bundle; then
\begin{enumerate}
  \item
   Assume that $G_0$ is locally free; then $F$ has a unique (to within
   permutations of quotients) exceptional filtration
   $$Gr(F)=(x_1F_1,x_0F_0)\qquad\quad(x_i\geq 0).$$
   Moreover,
     \begin{enumerate}
    \item
    the pair $(F_0,F_1)$ is obtained from the pair $(G_0,G_1)$ by mutations,
   \item
   $r(F_0)+r(F_1)\geq r(G_0)+r(G_1)
    \qquad\mbox{ whenever } x_0\cdot x_1\not= 0$,
   \item
   $r(F_0)\geq r(G_0)+r(G_1) \qquad\mbox{ provided } x_1= 0$,
   \item
   the equality of the sums of ranks is achieved if and only if $F_i=G_i$
   for $i=1,2$.
  \end{enumerate}
  \item
  Assume that   $G_0={\cal O}_e(-1)$
   for the exceptional rational curve $e=e_d$; then
   $$F=x_1F_1\oplus x_0F_0\qquad\quad(x_i\geq 0).$$
  Moreover,
   \begin{enumerate}
    \item
   the exceptional pair $(F_0,F_1)$ is obtained from the pair  $(G_0,G_1)$
   by mutations,
   \item
   $r(F_0)+r(F_1)\geq r(G_0)+r(G_1), \qquad
   \mbox{ whenever } x_0\cdot x_1\not= 0$,
   \item
   $r(F_0)\geq r(G_0)+r(G_1), \qquad\mbox{ provided } x_1= 0$,
   \item
    $F_0$  is locally free.
  \end{enumerate}
  \end{enumerate}
\ETH

\vspace{2ex}

  To prove this proposition, we need several lemmas.

\TH{Lemma.}\label{l342}
  let $A$ and $B$ be sheaves on a manifold $X$ and let
  $\varphi : V\otimes A\longrightarrow W\otimes B$
  be a morphism of  sheaves. Then
    \begin{enumerate}
 \item
   The canonical map $lcan:{\rm Hom}(A,B)\otimes A\longrightarrow B$ is an
   epimorphism provided $\varphi $ is also an epimorphism.
\item The canonical map $rcan:A\longrightarrow {\rm Hom}(A,B)^*\otimes B$ is
   a monomorphism provided $\varphi $ is a monomorphism as well.
    \end{enumerate}
\ETH

\PR In view of symmetry of statements it is sufficient to check the first
of them. At first consider the case of the one-dimensional space $W$, i.e.
 $$\varphi : V\otimes A\longrightarrow B\longrightarrow 0.$$
   Recall that the canonical map $lcan$ is determined by the element of
   ${\rm Hom}(A,B)^*\otimes {\rm Hom}(A,B)$ corresponding to the identical
   morphism ${\rm Hom}(A,B)\longrightarrow {\rm Hom}(A,B)$. Denote by
   $lcan$ this element as well. Let us define a line map
    $\psi:V\longrightarrow {\rm Hom}(A,B)$ so that
    $$\psi^*\otimes id_{{\rm Hom}(A,B)}:lcan\mapsto \varphi .$$
   Thus the following commutative diagram is arisen:
$$
\begin{array}{ccccc}
   {\rm Hom}(A,B)\otimes A& \stackrel{lcan}{\longrightarrow }&
   B&&\\
   \llap{$\scriptstyle \psi\otimes id_A $}\uparrow&&
   \llap{$\scriptstyle id_B $}\uparrow&&\\
   V\otimes A&\stackrel{\varphi }{\longrightarrow }&B&\longrightarrow &0
\end{array},
$$
  From this diagram it follows that $lcan$ is an epimorphism.

  Now suppose,
 $$\varphi : V\otimes A\longrightarrow W\otimes B\longrightarrow 0.$$
 Then
$${\rm Hom}(A,W\otimes B)\otimes A\stackrel{lcan}{\longrightarrow }
  W\otimes B\longrightarrow 0
.$$
  We see that there is a commutative diagram
$$
\begin{array}{ccccc}
   0&&0&&\\
   \uparrow&& \uparrow&& \\
   {\rm Hom}(A,B)\otimes A& \stackrel{lcan}{\longrightarrow }&
   B&\longrightarrow &0\\
   \llap{$\scriptstyle \pi_*\otimes id_A $}\uparrow&&
   \llap{$\scriptstyle \pi $}\uparrow&&\\
   {\rm Hom}(A,W\otimes B)\otimes A&\stackrel{lcan }{\longrightarrow }&
   W\otimes B&\longrightarrow &0
\end{array},
$$
  where $\pi$ is a projection $W\otimes B\longrightarrow B\longrightarrow 0$.

\TH{Lemma.}\label{l343} Let $F$ be a rigid sheaf and $(A,B)$ an exceptional
  $hom$-pair of sheaves on the surface $S$. Then the following statements
  hold.
\begin{enumerate}
\item
  If the sequence
   \begin{equation}\label{lpri}
   0\longrightarrow F\longrightarrow xA\longrightarrow yB
   \longrightarrow 0
   \end{equation}
   is exact for positive integer $\ x$ and $y$  then
\begin{enumerate}
\item
    the left mutation of pair $(A,B)$ belongs to the basic category and
    it is  regular;
\item
    either $F=wA\oplus zL_AB$ or there exists an exact sequence
$$0\longrightarrow F\longrightarrow zL_AB\longrightarrow wA\longrightarrow 0
$$
   for some nonnegative integer  $\ z$ and $w$.
\end{enumerate}
\item   If the sequence
   \begin{equation}\label{rpri}
   0\longrightarrow xA\longrightarrow yB\longrightarrow F
   \longrightarrow 0
   \end{equation}
   is exact for positive integer $x$ and $y$ then
\begin{enumerate}
\item
        the right mutation of the pair $(A,B)$ belongs to the
    basic category and it is regular;
\item
    either $F=wB\oplus zR_BA$ or there exists an exact sequence
$$0\longrightarrow wB\longrightarrow zR_BA\longrightarrow F\longrightarrow 0
$$
   for some nonnegative integer $z$ and $w$.
\end{enumerate}
\end{enumerate}
\ETH

\PR In view of duality of the lemma statements it is sufficient to prove
  the first of them.

  The regularity of the left mutation of the  $hom$-pair
  $(A,B)$ follows from its definition, sequence  (\ref{lpri}) and lemma
    \ref{l342}. Note that in this case the pair $(L_AB,A)$ has also the type
    $hom$.

   Sequence (\ref{lpri}) yields that the sheaf $F$ belongs to category
   generated by the pair $(A,B)$. Therefore there exists a spectral sequence
   $E^{p,q}$  (\ref{l316}) convergent to $F$ on the principal diagonal. Its
   $E_1$ term has the form:
$$
\begin{array}{ccc}
E_1^{-1,1}={\rm Ext}^1(B,F)\otimes L_AB&\stackrel{d}{\longrightarrow }&
E_1^{0,1}={\rm Ext}^1(A,F)\otimes A \\
E_1^{-1,0}={\rm Ext}^0(B,F)\otimes L_AB&\stackrel{d}{\longrightarrow }&
E_1^{0,0}={\rm Ext}^0(A,F)\otimes A
\end{array}.
$$
  It follows from exact sequence (\ref{lpri}) and the fact that the pair
   $(A,B)$ is exceptional that the group ${\rm Ext}^0(B,F)$ is trivial.
   Hence the spectral sequence degenerates into two exact triples:
   $$0\longrightarrow C\longrightarrow {\rm Ext}^1(B,F)\otimes L_AB
      \longrightarrow {\rm Ext}^1(A,F)\otimes A\longrightarrow 0,$$
   $$0\longrightarrow {\rm Hom}(A,F)\otimes A\longrightarrow F
   \longrightarrow C\longrightarrow 0.$$

   Assume that ${\rm Hom}(A,F)\not= 0$. Consider the cohomology table
   corresponding to the first of these triples:
       $$\begin{array}{|c|ccccc|}
   \hline
   k&V\otimes {\rm Ext}^k(A,A)&\rightarrow
    &W\otimes {\rm Ext}^k(L_AB,A)&\rightarrow
    &{\rm Ext}^k(C,A)\\
   \hline
   0& * &&*  && ? \\
   1& 0 &&0  &&?  \\
   2& 0 &&0  && ?  \\
   \hline
   \end{array} $$
   where  ${\rm Ext}^1(A,F)^*=V$,
   ${\rm Ext}^1(B,F)^*=W$. The first and the second columns are filled with
   using the properties of the pair  $(L_AB,A)$.

  From the table the equality ${\rm Ext}^1(C,A)=0$ follows. This means that
   $$F=C\oplus {\rm Hom}(A,F)\otimes A.$$
  Since $F$ is rigid, we get ${\rm Ext}^1(A,C)=0$ and  ${\rm Ext}^1(A,A)=0$.
   Hence,   ${\rm Ext}^1(A,F)=0$.  Thus,
   $$C={\rm Ext}^1(B,F)\otimes L_AB$$
   and $F$ is a direct sum of multiplicities of the sheaves $A$ and $L_AB$.

   Assume that ${\rm Hom}(A,F)=0$; then the spectral sequence degenerates
   into the exact triple
  $$0\longrightarrow F\longrightarrow {\rm Ext}^1(B,F)\otimes L_AB
  \longrightarrow {\rm Ext}^1(A,F)\otimes A\longrightarrow 0.$$
  This concludes the lemma proof.

\TH{Lemma.}\label{l344} Let $(E_0,E_1)$ be an exceptional $ext$-pair of
 sheaves on a manifold $X$ with $\chi(E_0,E_1)<-1$. In addition, assume that
 for each natural $n$ the following sheaves are determined:
  $$E_{n+1}=R_{E_n}E_{n-1},\qquad\qquad E_{-(n+1)}=L_{E_{-n}}E_{1-n}.
  $$
  Suppose that for a given sheaf $F$ and for any positive integer $n$ there
  are natural numbers $x_n,y_n,z_n,w_n$ such that the following exact
  sequences
  $$0\longrightarrow F\longrightarrow x_nE_{-(n+1)}\longrightarrow
     y_nE_{-n}\longrightarrow 0,
  $$
  $$0\longrightarrow z_nE_n\longrightarrow w_nE_{n+1}\longrightarrow F
  \longrightarrow 0
  $$
  take place; then the Euler characteristic $\chi(F,F)$ is nonpositive.
\ETH

\PR Denote by $e_n$ the images of $E_n$ in $K_0(X)$. The modulus
  $K_0(X)$ inherits the bilinear form $\chi(\cdot,\cdot)$. Denote it by
 $(\cdot,\cdot)$.

  Using the lemma conditions we have
  $$(e_0,e_0)=(e_1,e_1)=1,\qquad(e_1,e_0)=0,\qquad(e_0,e_1)=-h<-1
  .$$
  By the definition of mutations of an $ext$-pair, we get
  $$ e_{-1}=e_1+he_0,\qquad e_2=he_1+e_0
  .$$

  It follows from the exact sequences and the lemma assumptions
  that all pairs $(E_n,E_{n+1})$  for $n\!\in\!\ZZ ^*$ have the type $hom$
  and both mutations of these pairs (except for the left mutation of
 $(E_1,E_2)$ and the right one of $(E_{-1},E_0)$) are regular
  (\ref{l342}).

  The following formulae are easily obtained from the definition of
  mutations of $ext$ and $hom$-pairs.
  $$e_{-n}=he_{1-n}-e_{2-n}\qquad(n>1)
  ,$$
  $$e_{n}=he_{n-1}-e_{n-2}\qquad(n>2)
  ,$$
  $$\forall n\!\in\!\ZZ: \qquad(e_n,e_n)=1,\quad(e_{n+1},e_n)=0
$$
  $$\mbox{and for  }\quad n\not= 0\qquad(e_n,e_{n+1})=h
  .$$

  Denote by $x_n$ and $x_{n-1}$ coordinates of vector  $e_n\quad(n>0)$
  with respect to the basis  $\{e_1,e_0\}$: $\quad e_n=x_ne_1+x_{n-1}e_0$.
The recurrence relations  $$x_0=0,\quad x_1=1,\qquad x_{n+1}=hx_n-x_{n-1}
$$  are proved by induction on $n$.

  Note that the vectors   $e_{-n}\quad(n>0)$ are expressed by means of the
  same numbers:
  $$   e_{-n}=x_{n-1}e_1+x_ne_0
  .$$

  Let $V$ be a 2-dimensional vector space over ${\BBB Q}$ generated by
  $e_0,e_1$. Let us choose an affine map $U$ in ${\BBB P}(V)$ containing
  the image of $e_0$ as origin.
  $$xe_1+ye_0\quad\leadsto \frac{x}{y}e_1+e_0
  .$$
  We preserve the notations for the images of $e_n$ on $U$. Let us calculate
  coordinates $l_+$ and $l_-$ of limit points
  $e_{+\infty}=\lim_{n \to \infty}e_n,\quad
    e_{-\infty}=\lim_{n \to \infty}e_{-n}$ on $U$.
  $$  l_{+}=\lim_{n \to \infty}\frac{x_n}{x_{n-1}}=
  h- \lim_{n \to \infty}\frac{x_{n-2}}{x_{n-1}}
  =h-l_-=h-1/l_+
  .$$
  Hence, $l_+$ and $l_-$  are roots of the equation $l^2-hl+1=0$. That is,
  $$l_{\pm}=\frac{h\pm \sqrt{h^2-4}}{2}
  $$ (by assumption, $h\geq 2$).
  Taking into account the exact triples from the lemma condition, we see
  that the point $f$ on $U$ corresponding the sheaf $F$ has a coordinate
  $x\in [l_-,l_+]$.

  On the other hand, a sing of $\chi(F,F)$ is determined by a sign of
  $(e_0+xe_1)^2=x^2-hx+1$.
  Now the lemma proof follows from the inequality
  $$x^2-hx+1\le 0\qquad\mbox{for }\quad x\in [l_-,l_+]
  .$$

\TH{Corollary of the proof.}\label{l345} Under the conditions of the
  previous lemma, we have $r(E_n)\ge r(E_0)+r(E_1)$ (for
   $n\not= 0$ and $n\not= 1$ ). Moreover, $r(E_n)> r(E_0)+r(E_1)$  for
  $n\not= 0$ and $n\not= 1$ whenever both $E_0$ and $E_1$  has a positive
  rank.
  \ETH

  \PR In reality, we see that the image of the sheaf $E_n$ in $K_0(X)$
  has the form: $e_n=ae_0+be_1$ for some natural $a$ and $b$. Thus our
  statement follows from the additivity of rank function.

\vspace{3ex}

 {\sc Proof of proposition \ref{l341}} Suppose $G_i$ are locally free.
 If the pair $(G_0,G_1)$ has zero or $hom$ type then  $h^1(G_0,G_1)=0$
  and $F=y_0G_0\oplus y_1G_1$.

 If the pair $(G_0,G_1)$ is singular then $\mu_H(G_0)=\mu_H(G_1)$
 and the proof follows from the uniqueness of an exceptional filtration
 (\ref{l251}).

 Now, suppose  $G_0$ is a torsion sheaf and $G_1$ is a bundle; then the pair
 $(G_0,G_1)$ is necessarily an $ext$-pair.

Thus, let $(G_0,G_1)$ be an $ext$-pair. Following traditions, put

 $$G_{n+1}=R_{G_n}G_{n-1}\qquad\mbox{and}\qquad G_{-n}=L_{G_{1-n}}G_{2-n}
 .$$

 {\sc Step 1.} {\sl One of the following possibilities holds true}:
 $$F=x_1G_1\oplus x_2G_2,
 $$
 $$0\longrightarrow x_1G_1\longrightarrow x_2G_2\longrightarrow F
 \longrightarrow 0.
 $$
 Consider the spectral sequence convergent to $F$, which is constructed
 by the right dual collection $(G_1^{\vee},G_0^{\vee})$ (\ref{l316}).
 (Recall that  $G_1^{\vee}=G_1$ and $G_0^{\vee}=R_{G_1}G_0=G_2$.)
 Since the right mutation of the pair $(G_0,G_1)$ is nonregular, we get
 $$\Delta_0=1\quad\mbox{and}
 \quad E_1^{0,q}={\rm Ext}^{-q}(F,G_0)^*\otimes G_2.
 $$
 In addition, we do not mutations to obtain the sheaf $G_1^{\vee}$. Hence,
 $\Delta_1=0$ and
 $$E_1^{-1,q}={\rm Ext}^{1-q}(F,G_1)^*\otimes G_1
. $$
 Thus we see that $E_1$ term of the spectral sequence has the form
 $$
 \begin{array}{ccc}
 E_1^{-1,1}={\rm Ext}^0(F,G_1)^*\otimes G_1&\stackrel{d}{\longrightarrow }&
 0\\
 E_1^{-1,0}={\rm Ext}^1(F,G_1)^*\otimes G_1&\stackrel{d}{\longrightarrow }&
 E_1^{0,0}={\rm Ext}^0(F,G_0)^*\otimes G_2\\
 E_1^{-1,-1}={\rm Ext}^2(F,G_1)^*\otimes G_1&\stackrel{d}{\longrightarrow }&
 E_1^{0,-1}={\rm Ext}^1(F,G_0)^*\otimes G_2\\
 0&\stackrel{d}{\longrightarrow }&
 E_1^{0,-2}={\rm Ext}^2(F,G_0)^*\otimes G_2
 \end{array}
  .$$

 Using cohomology tables corresponding to the exact sequence from the
 proposition condition:
 $$0\longrightarrow y_1G_1\longrightarrow F\longrightarrow y_0G_0
 \longrightarrow 0
, $$
 let us calculate the groups ${\rm Ext}^k(F,G_i)$.
 $$
 \begin{array}{|c|ccccc|}
 \hline
 k&y_0{\rm Ext}^k(G_0,G_1)&\rightarrow&
   {\rm Ext}^k(F,G_1)&\rightarrow&
   y_1{\rm Ext}^k(G_1,G_1)\\
 \hline
 &0 &&? && *\\
 &* &&? && 0\\
 &0 &&? && 0\\
 \hline
 \end{array}
 $$
 $$
 \begin{array}{|c|ccccc|}
 \hline
 k&y_0{\rm Ext}^k(G_0,G_0)&\rightarrow&
   {\rm Ext}^k(F,G_0)&\rightarrow&
   y_1{\rm Ext}^k(G_1,G_0)\\
 \hline
 &* &&? && 0\\
 &0 &&? && 0\\
 &0 &&? && 0\\
 \hline
 \end{array}
 $$
  Whereby, the spectral sequence degenerates into two exact triples:
 $$0\longrightarrow {\rm Ext}^1(F,G_1)^*\otimes G_1
 \longrightarrow {\rm Ext}^0(F,G_0)^*\otimes G_2
 \longrightarrow C\longrightarrow 0
 ,$$
 $$ 0\longrightarrow C\longrightarrow F\longrightarrow
 {\rm Ext}^0(F,G_1)^*\otimes G_1\longrightarrow 0
 .$$
 Now, as in the proof of lemma \ref{l343}, using the first of these triples,
 it is easily shown that ${\rm Ext}^1(G_1,C)=0$. Therefore, if
  ${\rm Ext}^0(F,G_1)\not= 0$ then $F$ is a direct sum. In the converse
  case, the sheaf $F$ is included in the exact sequence.

 \vspace{3ex}

 {\sc Step 2.} {\sl One of the following possibilities holds true:}
 $$ F=x_{-1}G_{-1}\oplus x_0G_0
 ,$$
 $$0\longrightarrow F\longrightarrow x_{-1}G_{-1}\longrightarrow
 x_0G_0\longrightarrow 0
 .$$

 This step is checked in the same way as the first one with using the
 spectral sequence associated with the left dual collection $(G_{-1},G_0)$.

 \vspace{3ex}

 {\sc Step 3.} {\sl The sheaf $F$ is decomposed into the direct sum }:
 $$ F=x_{n-1}G_{n-1}\oplus x_nG_n
 $$
 for some $n\!\in\!\ZZ$ and nonnegative integer $x_{n-1},x_n$.
(That is $F_0=G_{n-1}$ and $F_1=G_n$  in the formulation of proposition.)

 Using the first two steps and lemma  \ref{l343}, it can be stated that for
 any $n>0$ the following exact triples
 $$0\longrightarrow x_nG_n\longrightarrow x_{n+1}G_{n+1}\longrightarrow F
 \longrightarrow 0
 ,$$
 $$0\longrightarrow F\longrightarrow x_{-n}G_{-n}\longrightarrow
 x_{1-n}G_{1-n}\longrightarrow 0
 $$
 hold unless $F=x_{n-1}G_{n-1}\oplus x_nG_n$.

 Let us show that these triples contradict the proposition conditions.

 Suppose $h^1(G_0,G_1)>1$; then it follows from these sequences and lemma
 \ref{l344} that $\chi(F,F)\le 0$. This contradicts the fact that $F$ is
 rigid.

 Suppose $h^1(G_0,G_1)=1$; then the series of the exceptional sheaves $G_n$
 is formed by $G_0,G_1,G_2$. In reality, in this case both the right and the
 left mutation of the $ext$-pair  $(G_0,G_1)$  is described by the sequence
 $$0\longrightarrow G_1\longrightarrow G_2\longrightarrow G_0
 \longrightarrow 0
 .$$
 Whence, $L_{G_0}G_1=G_2$ and $R_{G_2}G_1=G_0$. Hence there are exact
 triples:
 $$0\longrightarrow y_1G_1\longrightarrow F\longrightarrow
 y_0G_0\longrightarrow 0
 ,$$
 $$0\longrightarrow x_2G_2\longrightarrow x_0G_0\longrightarrow F
 \longrightarrow 0
 .$$
 Since $G_0$ is indecomposable, it follows from the second sequence that
 $h^1(F,G_2)\not= 0$. We apply the functor ${\rm Ext}^\cdot(\cdot,G_2)$ to
 the first triple to obtain
 $$y_0{\rm Ext}^1(G_0,G_2)\longrightarrow
 {\rm Ext}^1(F,G_2)\longrightarrow y_1{\rm Ext}^1(G_1,G_2)
 .$$
 Since the pair  $(G_2,G_0)$ is exceptional, we get ${\rm Ext}^1(G_0,G_2)=0$.
 Besides, since $(G_1,G_2)$ is a $hom$-pair, we obtain
 ${\rm Ext}^1(G_1,G_2)=0$. Thus,  $h^1(F,G_2)=0$. This contradiction proves
 the 3-th step.

 \vspace{3ex}

 {\sc Step 4.} {\sl Suppose $G_1$ is a bundle and $G_0={\cal O}_e(-1)$; then
 $F$ is locally free or $F_0$ is a bundle and $F_1={\cal O}_e(-1)$.}

  By the proposition assumption, the sheaf $F$ is included in the exact
  triple:
  $$0\longrightarrow y_1G_1\longrightarrow F\longrightarrow
  y_0{\cal O}_e(-1)\longrightarrow 0
.  $$
  Since $F$ is rigid, we see that $F$ is locally free
  whenever $F$ has not torsion (\ref{l221}).
  Therefore its direct summands are locally free as well.

  Assume that $F$ has a torsion $TF$. Since   $G_1$ is locally free, we
  obtain  the following commutative diagram:
$$\begin{array}{ccccccccc}
&&0&&0&&0&&\\
&&\uparrow&&\uparrow&&\uparrow&&\\
0&\longrightarrow &y_1G_1&\longrightarrow &F'&\longrightarrow& Q&
\longrightarrow&0\\
&&\uparrow&&\uparrow&&\uparrow&&\\
0&\longrightarrow &y_1G_1&\longrightarrow &F&\longrightarrow&
y_0{\cal O}_e(-1)&\longrightarrow&0\\
&&\uparrow&&\uparrow&&\llap{$\scriptstyle\varphi $}\uparrow&&\\
&&0&\longrightarrow &TF&\longrightarrow&TF&\longrightarrow&0\\
&&&&\uparrow&&\uparrow&&\\
&&&&0&&0&&\\
\end{array},$$ where $F'$ is torsion free.
 Since $TF$ as a subsheaf of $y_0{\cal O}_e(-1)$ and the curve $e$
 is isomorphic to the projective line, we get
$$TF\cong \bigoplus\limits_{i}^{}z_i{\cal O}_e(s_i).$$
 Hence,
$$Q\cong \left[\bigoplus\limits_{j}^{}w_j{\cal O}_e(d_j)\right]\oplus T^0,$$
 where $T^0$ is a torsion sheaf with a zero-dimensional support.

 Consider the upper row of the diagram. Assume that
 $T^0\not= 0,$. Since $G_1$ is locally free and the support of $T^0$ is
 zero-dimensional, we get ${\rm Ext}^1(T^0,y_1G_1)=0$.
 Hence $T^0$ is a direct summand of $F'$. But this contradicts the fact
 that $F'$ has not torsion. For the same reason, $${\rm Ext}^1({\cal
 O}_e(d_j),G_1)\not= 0\qquad \forall j.$$ Let us show that this yields the
 inequality $d_j\leq -1$.

 In reality, by assumption, $({\cal O}_e(-1),G_1)$ is an exceptional pair.
 Therefore, it easily follows from a calculation of cohomologies that
 $(G_1)|_e=r(G_1){\cal O}_e$. Thus,
 $${\rm Ext}^1({\cal O}_e(d_j),G_1)^*\cong
  {\rm Ext}^1(G_1,{\cal O}_e(d_j)\otimes K_S)=
  r(G_1){\rm Ext}^1({\cal O}_e,{\cal O}_e(d_j-1))\not= 0$$
 and the inequalities $d_j\leq -1$ hold for all $j$.

  On the other hand, $Q$ is a quotient of $y_0{\cal O}_e(-1)$. Hence,
  $d_j\geq -1\quad \forall j.$  From these inequalities it follows that
   $d_j=-1\quad \forall j$ and $Q=w{\cal O}_e(-1)$.

   We see that $TF=z{\cal O}_e(-1)$.

   Now, note that from the exact sequence
$$ 0\longrightarrow y_1G_1\longrightarrow F'\longrightarrow w{\cal O}_e(-1)
\longrightarrow 0$$
 it follows that ${\rm Ext}^1(F',{\cal O}_e(-1))=0$, i.e.
 $F=F'\oplus z{\cal O}_e(-1)$.

  By the previous step, $F=x_{n-1}G_{n-1}\oplus x_nG_n$.
 Therefore,  $F=x_{-1}G_{-1}\oplus x_0{\cal O}_e(-1)$ or
  $F=x_{-1}G_{-1}\oplus x_0{\cal O}_e(-1)\oplus x_{1}G_{1}. $
  Since $({\cal O}_e(-1),G_{1})$ is the $ext$-pair and $F$ is
   superrigid sheaf, we see that the last inequality is impossible.
  On the other hand, $F'$ is locally free, as rigid sheaf without torsion
(  \ref{l221}). Thus the sheaf $x_0F_0=x_{-1}G_{-1}=F'$ is locally free as
 well.

 \vspace{3ex}

 {\sc Step 5.}        {\sl $r(F_0)+r(F_1)>r(G_0)+r(G_1)$
 for $x_0\cdot x_1\not= 0$, and $r(F_0)\geq r(G_0)+r(G_1)$
 for $x_1=0$\/}.

 Since  $F_0$ and $F_1$ are direct summands of a superrigid sheaf, we
 obtain that the pair $(F_0,F_1)$ is exceptional and it has the type $hom$.
 Therefore it does not coincides with the pair $(G_0,G_1)$.
  In view of this the first inequality follows from corollary \ref{l345}.

 Suppose $F=x_0F_0$; then $F_0\not= G_0,G_1$. By the same argument,
 $r(F_0)>r(G_0)+r(G_1)$. The equality of ranks is possible here only if
 $F_0=G_1$ and $G_0={\cal O}_e(-1)$. This completes the proof.

%% FOLLOWING LINE CANNOT BE BROKEN BEFORE 80 CHAR
%%%%%%%%%%%%%%%%%%%%%%%%%%%%%%%%%%%%%%%%%%%%%%%%%%%%%%%%%%%%%%%%%%%%%%%%%%%%%%%%%%
\subsection{Proof of the Main Theorem.}

   It follows from lemma \ref{l334} that for any exceptional collection
   of bundles on the surface $S$ satisfying the conditions of the main
   theorem there is  a $hom$-collection
$$\tau=\OVEC F,k$$ equivalent to the original one such that the superrigid
bundle $F$ associated with $\tau$ is  included in  exact sequence
(\ref{orl}):  $$0\longrightarrow G\longrightarrow F\longrightarrow {\rm
Hom}(F,{\cal O}_e(-1))^*\otimes {\cal O}_e(-1)\longrightarrow 0 ,   $$
   where $G$ is a superrigid bundle with $\qquad {\rm
   Ext}^k(G,{\cal O}_e(-1))=0 \qquad\forall k=0,1,2$.  (Further we shall
 denote by ${\cal B}\OVEC F,k$ the superrigid bundle associated with a
$hom$-collection $\OVEC F,k$.)

  In particular, we see that $G|_e=s{\cal O}_e.$ Therefore there exists a
  superrigid bundle $G'$ on the surface $S'$ obtained from $S$ by blowing
  down the curve $e$$\quad (\sigma:S\longrightarrow S')$ such that
 $\sigma^*(G')=G.$

  Since $G'$ is superrigid, we see that there is an exceptional filtration
  of it:
  $$Gr(G')=(y_nG'_n,y_{n-1}G'_{n-1},\ldots,y_1G'_1).$$
  Using the induction on the number of blow up divisors on $S$, we can
  assume that the exceptional collection of bundles   $\VEC {G'},n$
  is constructible. That is it included in a full exceptional collection
  obtained from the basic collection
  $$\Bigl({\cal O}_S,{\cal O}_S(h),{\cal O}_S(2h),{\cal O}_{e_1}(-1),
  \ldots ,{\cal O}_{e_{d-1}}(-1)\Bigr)
$$  by mutations.
 (Note that $K_{S'}^2=K_{S}^2+1>1$. Therefore the constructibility of the
 collection $\VEC G',n$ does not depend on ranks of the sheaves  $G'_j$
(see theorem \ref{l318})).

  Let us remember that the base of the induction, i.e. the case of the
  projective plane, had been checked in the paper \cite{RP}.

  Since $\sigma^*(G')=G$, we obtain that the bundle $G$ has the exceptional
  filtration
  $$Gr(G)=(y_nG_n,y_{n-1}G_{n-1},\ldots,y_1G_1),$$
  where  $G_i=\sigma^*(G'_i)$. Moreover, the collection
   $\tau '=({\cal O}_e(-1),G_1, \ldots,G_n)$  is exceptional (the
   triviality of the groups ${\rm Ext}^k(G_i,{\cal O}_e(-1))$ follows from
    the fact that  $G_i|_e=s_i{\cal O}_e$).  Furthermore, the
  constructibility of the collection $\VEC {G'},n$ implies the
  constructibility of $\tau'$.

 Now to prove theorem \ref{l318} it is sufficient to show that the
 collection $\tau$ is included in an exceptional collection obtained
  from $\tau '$ by mutations.

 Let us illustrate the procedure of that inclusion in the projectivisation
 of $K_0(S)\otimes \QQ =K$. To each sheaf $E$ on $S$ assign a vector
 $[E]$ in $K$. It is obvious that vectors corresponding to sheaves from an
 exceptional collection are linearly independent. Recall that the
 nonsingular bilinear form $(\cdot,\cdot)$ on $K$ is defined. It corresponds
 to Euler characteristic of sheaves $\chi(E,F)$.
  Since all exceptional sheaves satisfy the equation $\chi(E,E)=1$,
  we see that the corresponding vectors are not proportional. Let us pass to
  the projectivisation of $K$. In this case, vectors corresponding to sheaves
  of an exceptional collection are projected to vertexes of some simplex.

  The key exact sequence implies that the vector $[F]$ gets into the
  simplex with the vertexes
   $[{\cal O}_e(-1)],[G_1],...,[G_n]$.
%%%%%%%%%%%%%%%%%%%%%%%%%%%%%%%%%%%%%%%%%%%%%%%%%%%%%%%%%%%%%%%%%%%%%%%
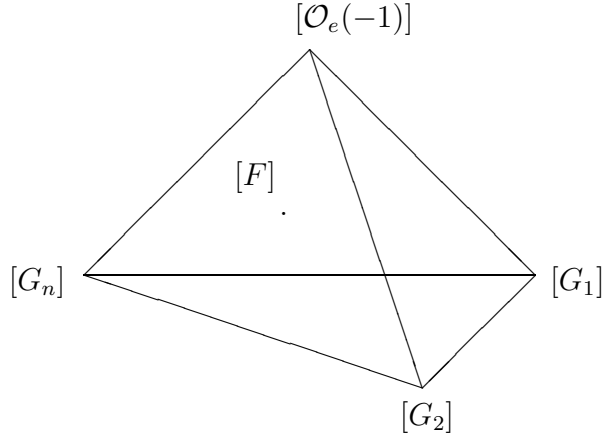
\begin{center}
\begin{picture}(100,85)(-50,-30)
 \put(0,35){\line(1,-1){30}}

 \put(0,35){\line(1,-3){15}}
 \put(0,35){\line(-1,-1){30}}
 \put(30,5){\line(-1,-1){15}}
 \put(30,5){\line(-1,0){60}}
 \put(15,-10){\line(-3,1){45}}
%%%%%%%%%%%%%%%%%%%%%%%%%%%%%%%%%%%%%%%%%%%%%%%%%%%%%%%%%%%%%
  \put(-2,38){$[{\cal O}_e(-1)]$}
  \put(32,3){$[G_1]$}
  \put(12,-15){$[G_2]$}
  \put(-40,3){$[G_n]$}
  \put(-4,13){.}
  \put(-10,17){$[F]$}

\end{picture}

\end{center}
%%%%%%%%%%%%%%%%%%%%%%%%%%%%%%%%%%%%%%%%%%%%%%%%%%%%%%%%%%%%%%%%%%

 Let us project the point $[F]$ to the edge $([{\cal O}_e(-1)],[G_1])$.
 Note that this projection corresponds to a superrigid sheaf, and the
 exceptional pair $(G'_0,G'_1)$ associated with it obtained by mutations
 of the pair $({\cal O}_e(-1),G_1)$. As a result, we get a lesser simplex
 containing $[F]$. Next let us project $[F]$ to the face
 $([G'_1],[G_2],\ldots,[G_n])$, etc... It remains to show that this
 process is finite.

Let us prove two lemmas about projections.

\TH{Lemma.}\label{l351} Let
  \begin{equation}\label{rpr1}
  0\longrightarrow G\longrightarrow F\longrightarrow E\longrightarrow 0
  \end{equation}
 be an exact sequence of superrigid sheaves on the surface $S$. Let
  $$Gr(E)=(y_kE_k,y_{k-1}E_{k-1},\ldots,y_1E_1),$$
$$ Gr(G)=(y_mG_m,y_{m-1}G_{m-1},\ldots,y_{k+1}G_{k+1}) $$
 be exceptional filtrations of $E$ and $G$ such that the collection
   $$(E_1,\ldots,E_k,G_{k+1},\ldots,G_m)$$
    is exceptional. Let us divide the
   filtration of the sheaf $G$ into two groups
  \begin{equation}\label{rpr2}
  0\longrightarrow G'\longrightarrow G\longrightarrow G''\longrightarrow 0,
  \end{equation}
  where  $G'$ and $G''$ are the sheaves with the exceptional filtrations
  $$Gr(G')=(y_mG_m,y_{m-1}G_{m-1},\ldots,y_{s+1}G_{s+1}),$$
$$  Gr(G'')=(y_sG_s,y_{s-1}G_{s-1},\ldots,y_{k+1}G_{k+1}).
   $$
   Then
\begin{enumerate}
\item
 $G'$ and $G''$ are superrigid;
\item
 ${\rm End}(G')\cong {\rm Hom}(G',F)$;
\item
 ${\rm Ext}^i(G',F)=0$ for $i>0$;
\item
 ${\rm Ext}^2(F,G')=0$;

\item
there is an exact sequence:
\end{enumerate}
  \begin{equation}\label{rpr3}
  0\longrightarrow G'\longrightarrow F\longrightarrow E'\longrightarrow 0,
  \end{equation}
 where $E'$ is a superrigid sheaf included in the exact triple
  \begin{equation}\label{rpr4}
  0\longrightarrow G''\longrightarrow E'\longrightarrow E\longrightarrow 0.
  \end{equation}
 Besides, ${\rm Ext}^i(G',E')=0\quad\forall i$.
\ETH

\PR By the definition of an exceptional collection, ${\rm Ext}^k(G_j,G_i)=0$
  for $j>i$ and all $k$. Therefore,
   $\forall k: \quad {\rm Ext}^k(G',G'')=0$
  (\ref{l124}). Hence it follows from lemma \ref{l222} that
  ${\rm Ext}^2(G'',G')=0$. We apply the Mukai lemma to exact sequence
 (\ref{rpr2}) to obtain that $G'$ and $G''$ are rigid. Since the collection
  $(G_{k+1},...,G_m)$ is exceptional, we see that
 ${\rm Ext}^2(G_i,G_j)=0$ for any pair $i,j$. This implies that
 ${\rm Ext}^2(G',G')={\rm Ext}^2(G'',G'')=0$. Thus the first lemma
  statement holds.

  We saw that ${\rm Ext}^k(G',G'')=0\quad\forall k.$  Whence, using exact
  triple   (\ref{rpr2}) and the fact that $G'$ is superrigid, we have
  $${\rm Hom}(G',G)\cong {\rm End}(G')\qquad\mbox{and}\qquad
  {\rm Ext}^i(G',G)=0\quad\mbox{for } i>0.
  $$
 Besides, in view of the definition of the sheaf $G'$ and the fact that
 the collection  $$(E_1,\ldots,E_k,G_{k+1},\ldots,G_m)$$
  is exceptional the
 following inequalities are valid.
$${\rm Ext}^i(G',E)=0\quad\forall i;\qquad{\rm Ext}^2(E,G')=
  {\rm Ext}^2(G,G')=0.
$$
  Consider two cohomology tables corresponding to sequence   (\ref{rpr1}).
$$\begin{array}{|c|ccccc|}
\hline
k&{\rm Ext}^k(G',G)&\rightarrow&
  {\rm Ext}^k(G',F)&\rightarrow&
  {\rm Ext}^k(G',E)\\
\hline
&{\rm End}(G') &&? &&0 \\
&0 &&? &&0 \\
&0 &&? &&0 \\
\hline
\end{array}$$

$$\begin{array}{|c|ccccc|}
\hline
k&{\rm Ext}^k(E,G')&\rightarrow&
  {\rm Ext}^k(F,G')&\rightarrow&
  {\rm Ext}^k(G,G')\\
\hline
&* &&? &&* \\
&* &&? &&* \\
&0 &&? &&0 \\
\hline
\end{array}   $$

 The lemma statements 2, 3 and 4  follow from this tables.

  Exact triples (\ref{rpr1}) and (\ref{rpr2}) give the following commutative
  diagram:
$$
  \begin{array}{ccccccccc}
  &&0&&&&&&\\
  &&\uparrow&&&&&&\\
  &&G''&&0&&0&&\\
  &&\uparrow&&\uparrow&&\uparrow&&\\
  0&\longrightarrow &G&\longrightarrow &F&\longrightarrow &E&
  \longrightarrow &0\\
  &&\uparrow&&\uparrow&&\uparrow&&\\
    0&\longrightarrow &G'&\longrightarrow &F&\longrightarrow &E'&
  \longrightarrow &0\\
  &&\uparrow&&\uparrow&&\uparrow&&\\
  &&0&&0&&G''&&\\
  &&&&&&\uparrow&&\\
  &&&&&&0&&
  \end{array}.
$$
  It yields exact sequences (\ref{rpr3}) and
  (\ref{rpr4}).

  Now to prove the lemma it remains to check that the sheaf
  $E'$ is superrigid and
 for all $i\quad$ ${\rm Ext}^i(G',E')=0$.
 All these facts follow from the following cohomology tables associated with
 sequence  (\ref{rpr3})
$$\begin{array}{|c|ccccc|}
\hline
k&{\rm Ext}^k(G',G')&\rightarrow&
  {\rm Ext}^k(G',F)&\rightarrow&
  {\rm Ext}^k(G',E')\\
\hline
&{\rm End}(G') &&{\rm Hom}(G',F) && ?\\
&0 &&0 &&? \\
&0 &&0 &&? \\
\hline
\end{array}
$$
$$
\begin{array}{|c|ccccc|}
\hline
k&{\rm Ext}^k(F,G')&\rightarrow&
  {\rm Ext}^k(F,F)&\rightarrow&
  {\rm Ext}^k(F,E')\\
\hline
&* &&* &&? \\
&* &&0 &&? \\
&0 &&0 &&? \\
\hline
\end{array}
$$
$$
\begin{array}{|c|ccccc|}
\hline
k&{\rm Ext}^k(E',E')&\rightarrow&
  {\rm Ext}^k(F,E')&\rightarrow&
  {\rm Ext}^k(G',E')\\
\hline
&? &&* &&0 \\
& ?&&0 && 0\\
&? &&0 &&0 \\
\hline
\end{array}.$$

This completes the proof.

\vspace{2ex}

  By the same argument the dual statement can be proved.

\TH{Lemma.}\label{l352} Under the conditions of the previous lemma let us
  divide the filtration of the sheaf $E$ into two groups:
  $$0\longrightarrow E'\longrightarrow E\longrightarrow E''\longrightarrow
  0,$$
  where $E'$ and $E''$ are sheaves with the exceptional filtrations
$$Gr(E')=(y_kE_k,y_{k-1}E_{k-1},\ldots,y_{s+1}E_{s+1}),\qquad
  Gr(E'')=(y_sE_s,y_{s-1}E_{s-1},\ldots,y_{1}E_{1}).
$$
 Then
\begin{enumerate}
\item
 $E'$ and $E''$ are superrigid;
\item
 ${\rm End}(E'')\cong {\rm Hom}(F,E'')$;
\item
 ${\rm Ext}^i(F,E'')=0$ for $i>0$;
\item
 ${\rm Ext}^2(E'',F)=0$;
\item
 There is an exact sequence:
$$  0\longrightarrow G'\longrightarrow F\longrightarrow E''\longrightarrow 0,
$$
 where $G'$ is a superrigid sheaf included in the exact triple:
$$  0\longrightarrow G\longrightarrow G'\longrightarrow E'\longrightarrow 0.
$$
 Besides, ${\rm Ext}^i(G',E'')=0\quad\forall i$.
\end{enumerate}
\ETH

\TH{Remark.}\label{l353} \begin{enumerate}
\item
 Lemma \ref{l351} is also valid provided $E=y_1{\cal O}_e(-1)$
 for the exceptional rational curve $e=e_d$;
\item
 lemma \ref{l352} is correct as well provided $E=y_1E_1\oplus y_2E_2$,
 where $E_1$ is an exceptional bundle and $E_2={\cal O}_e(-1)$;
\item
 the procedure described in \ref{l351} is called {\it the transfer of the
 collection  $(G_{k+1},\ldots,G_s)$ to the right\/}, and the similar
 procedure from \ref{l352} is
 {\it the transfer of the collection $(E_{s+1},\ldots,E_k)$ to the left\/}.
\end{enumerate}
\ETH

\vspace{2ex}

 Now let us prove a proposition concluding the proof of the main theorem.

\TH{Proposition.}\label{l354}
  Suppose a superrigid bundle $F={\cal B}\OVEC
F,k$ on the surface $S$ with $K_S^2>0$ is included in the exact sequence
  \begin{equation}\label{simp}
  0\longrightarrow G\longrightarrow F\longrightarrow E\longrightarrow 0,
  \end{equation}
 where $G$ is a superrigid bundle with an exceptional filtration
  $$Gr(G)=(y_nG_n,y_{n-1}G_{n-1},\ldots,y_sG_s),$$
 and $E$ is a superrigid sheaf. In addition we assume that
 $E$ is either locally free and
  $$Gr(E)=(y_{s-1}G_{s-1},y_{s-2}G_{s-2},\ldots,y_0G_0)$$
  or $E=y_0G_0=y_0{\cal O}_e(-1)$; but the collection $\OVEC G,n$
  is exceptional in any case. Then

\begin{enumerate}
\item
 $k\leq n$;
\item  the collection
 $\OVEC F,k$ is included in an exceptional collection obtained from
 $\OVEC G,n$ by mutations;
\item
 $\sum\limits_{i=0}^{k}r(F_i)\geq \sum\limits_{j=0}^{n}r(G_j)$ ;
\item
 If $E$ is locally free then the equality
 $\sum\limits_{i=0}^{k}r(F_i)= \sum\limits_{j=0}^{n}r(G_j)$ ;
  yields the equality $k=n$ . Moreover, in this case we have
 $F_i=G_i$ after some mutations of neighboring zero-pairs.
 \end{enumerate}
\ETH

\PR The proof is by induction on number of sheaves in the collection
   $$\OVEC G,n.$$ The case $n=1$ had been checked in the previous section.

\vspace{2ex}

{\sc Statement.} {\sl Without loss of generality it can be assumed that
 $E$ and $G$ is locally free.}

\vspace{2ex}

  Proof. Suppose   $E=y_0G_0=y_0{\cal O}_e(-1)$. Following remark \ref{l353},
  let us do the transfer of $G_1$ to the right. Namely, let us
  denote by $G'$ the bundle ${\cal B}(G_2,G_3,\ldots,G_n)$ and let us
  consider the exact sequences $$0\longrightarrow G'\longrightarrow
  F\longrightarrow E' \longrightarrow 0,$$ $$0\longrightarrow
y_1G_1\longrightarrow E'\longrightarrow y_0G_0 \longrightarrow 0.$$
 Taking into account lemma \ref{l351} and proposition \ref{l341},we
 obtain that $E'$ is a superrigid bundle such that $E'=x_0E'_0 \oplus
 x_1E'_1$ (or $E'=x_0E'_0$ ), where the exceptional pair $(E'_0,E'_1)$ (or
 $E'_0$) is obtained by mutations of the pair $(G_0,G_1)$. Moreover,
  $E'_0$ is locally free and $$r(E'_0)+r(E'_1)\ge r(G_0)+r(G_1) \qquad
  (r(E'_0)\ge r(G_0)+r(G_1)).$$ Let us show that the collection
  $(E'_0,E'_1,G_2,\ldots,G_n)$ is exceptional. From lemma \ref{l351} it
  follows that ${\rm Ext}^k(G',E')=0\quad\forall k$.  But,
  $E'=x_0E'_0\oplus x_1E'_1$ and $G'={\cal B}(G_2,G_3,\ldots,G_n)$.

  Provided $E'_i$ is locally free, the triviality of the groups
  ${\rm Ext}^k(G_j,E'_i)$ for $j=2,\ldots,n$ follows from lemma
   \ref{l257}. Let us check this property for the case
 $E'_1={\cal O}_e(-1)$. Since
  ${\rm Ext}^k(G',{\cal O}_e(-1))=0\quad\forall k$, we see that the
  restriction of   $G'$ to the curve $e$ is trivial.
  Therefore there exists a superrigid bundle $L$ on the surface $S'$
  obtained from $S$ by blowing down the curve $e$ ($\sigma:
  S\longrightarrow S'$) such that $\sigma^*(L)=G'$.

  Since $L$ is superrigid, we see that it has the exceptional filtration
 $$Gr(L)=(z_mL_m,z_{m-1}L_{m-1},\ldots,
  z_2L_2).$$ Besides,
  $Gr(G')=(z_m\sigma^*(L_m),z_{m-1}\sigma^*(L_{m-1}),\ldots,
  z_2\sigma^*(L_2))$ is the exceptional filtration of the bundle   $G'$.
  Now, by theorem \ref{l251}, $m=n$ and $G_i=\sigma^*(L_i)$.
  Thus the collection $(E'_0,E'_1,G_2,\ldots,G_n)$ is exceptional.

  Our statement is correct in the case $E'=y_0E'_0$.

  Assume that  $E'=x_0E'_0\oplus  x_1E'_1$  with positive $x_0,x_1$.
 Let us do the transfer of $E'_1$ to the left:
$$0\longrightarrow \tilde G \longrightarrow F\longrightarrow
x_0E'_0\longrightarrow 0,$$
$$0\longrightarrow G'\longrightarrow \tilde G \longrightarrow x_1E'_1
\longrightarrow 0  .
$$
  Using lemma \ref{l352} and the induction hypothesis, we obtain that
   $\tilde G $ is a superrigid sheaf with an exceptional filtration
  $Gr(\tilde G )=(x'_mG'_m,x'_{m-1}G'_{m-1},\ldots,x'_1G'_1)$. In addition,
  the collection $\VEC G',m$ is included in an exceptional collection
  obtained from $(E'_1,G_2,\ldots,G_n)$ by mutations and $$\sum r(G'_i)\geq
  \sum r(G_i)+ r(E'_1).$$

  Note that the sheaf $\tilde G $ has not torsion, as a subsheaf of a bundle.
  Since  $\tilde G $  is rigid, we see that it is locally free.
  It can be checked as above that the collection
 $(E'_0,G'_1,\ldots,G'_m)$ is exceptional. This completes the statement
 proof.

  We shall name {\it bounding\/} the collection $\OVEC G,n$ from
  the formulation of our proposition and all collections obtained from it
  by mutations.

  Now, consider exact sequence (\ref{simp}). We shall do the transfer of the
  bundle  $G_s$ to the right and to the left. Recall that in this procedure
 the sum of ranks of the bounding collections is not decreased.

  Since the sum of ranks of the bounding collections is less than or equals
  rank of the bundle $F$, we see that this process cannot be continued ad
  infinitum. Hence beginning with some moment the sum of ranks is a constant.
  Let us study this moment in the following statement.

  \vspace{2ex}

  {\sc Statement. } {\sl Assume that under the conditions of our
  proposition the sum of ranks of the bounding collection bundles does
  not change after the transfers of the bundle $G_s$ to the right and to
  the left;  then $k=n$ and $$\OVEC F,k=\OVEC G,n.$$ to within mutations
 of neighboring zero-pairs.}

 \vspace{2ex}

  Proof. After the transfer of the bundle $G_s$ to the right two exact
  sequence appear:
  $$0\longrightarrow {\cal B}(G_{s+1},\ldots,G_n)\longrightarrow
  F\longrightarrow {\cal B}\OVEC {G'},l \longrightarrow 0,$$
  $$0\longrightarrow y_sG_s\longrightarrow {\cal B}\OVEC {G'},l
  \longrightarrow {\cal B}\OVEC G,{s-1}\longrightarrow  0.$$
  Since $G_0,\ldots,G_{s-1},G_s$ are locally free, we obtain that by the
  induction  hypothesis it follows from the equality
  $$\sum\limits_{i=0}^{l}r(G'_i)=\sum\limits_{i=0}^{s}r(G_i)$$
  that  $l=s$ and $G_i=G'_i$ (to within mutations of neighboring zero-pairs).
  Therefore there is an exact sequence
  $$0\longrightarrow {\cal B}(G_{s+1},\ldots,G_n)\longrightarrow
    F\longrightarrow {\cal B}\OVEC G,s\longrightarrow 0.
  $$
   Moreover, $\OVEC G,s$ is the $hom$-collection. Whereby,
    $\mu_H(G_i)\geq \mu_H(G_j)$  for $s\ge i>j$.

  Now let us do the transfer of the bundle $G_s$ to the left (by
  assumption, the sum of ranks does not change as well):
  $$0\longrightarrow {\cal B}(G''_{s},\ldots,G''_m)\longrightarrow
  F\longrightarrow {\cal B}\OVEC {G},{s-1} \longrightarrow 0,$$
  $$0\longrightarrow {\cal B}(G_{s+1},\ldots,G_n)
  \longrightarrow {\cal B}(G''_{s},\ldots,G''_m)
  \longrightarrow  y_sG_s
  \longrightarrow 0.$$
  As before, by the induction hypothesis, we obtain that the collection
  $(G''_{s},\ldots,G''_m)$ consists with the collection $(G_s,\ldots,G_n)$
  to within mutations of neighboring zero-pairs. Hence
  $(G_s,\ldots,G_n)$ is the $hom$-collection and  $\mu_H(G_j)\le \mu_H(G_i)$.
  for $s\le j<i$.

   As a result we obtain that the all bounding collections $\OVEC G,n$
  is the $hom$-collection.
  Thus we can construct the  exceptional filtrations of the bundle
  $$F={\cal B}\OVEC F,k$$ from the exceptional filtrations of the
  bundles $E$ and $G$ in sequence (\ref{simp}).

Now the proof following from the uniqueness of the exceptional filtration.

%% FOLLOWING LINE CANNOT BE BROKEN BEFORE 80 CHAR
%%%%%%%%%%%%%%%%%%%%%%%%%%%%%%%%%%%%%%%%%%%%%%%%%%%%%%%%%%%%%%%%%%%%%%%%%%%%%%%%%%

\newpage
\tableofcontents
%

%% FOLLOWING LINE CANNOT BE BROKEN BEFORE 80 CHAR
%%%%%%%%%%%%%%%%%%%%%%%%%%%%%%%%%%%%%%%%%%%%%%%%%%%%%%%%%%%%%%%%%%%%%%%%%%%%%%%%%%%
\end{document}